\documentclass[12pt]{amsart}
\usepackage{amsmath,amsthm,amsfonts,amscd,amssymb,mathrsfs,latexsym}
\usepackage[all]{xy}

\theoremstyle{definition}
\newtheorem{theorem}{Theorem}[section]
\newtheorem{corollary}[theorem]{Corollary}
\newtheorem{lemma}[theorem]{Lemma}
\newtheorem{proposition}[theorem]{Proposition}
\newtheorem{remark}[theorem]{Remark}
\newtheorem{notation}[theorem]{Notation}
\newtheorem{example}[theorem]{Example}

\newtheorem{definition}[theorem]{Definition}
\newtheorem*{acknowledgments}{Acknowledgments}

\newcommand{\cA}{{\mathcal A}}
\newcommand{\cB}{{\mathcal B}}

\newcommand{\cD}{{\mathcal D}}
\newcommand{\cE}{{\mathcal E}}

\newcommand{\cG}{{\mathcal G}}

\newcommand{\cL}{{\mathcal L}}
\newcommand{\cO}{{\mathcal O}}
\newcommand{\cP}{{\mathcal P}}
\newcommand{\cX}{{\mathcal X}}
\newcommand{\cY}{{\mathcal Y}}
\newcommand{\cZ}{{\mathcal Z}}
\newcommand{\cU}{{\mathcal U}}

\newcommand{\sV}{{\mathscr V}}
\newcommand{\sM}{{\mathscr M}}
\newcommand{\sS}{{\mathscr S}}
\newcommand{\sC}{{\mathscr C}}
\newcommand{\sE}{{\mathscr E}}

\newcommand{\nl}{{\mathnormal l}}

\newcommand{\De}{{\mathbb D}}

\newcommand{\EA}{{\rm EA}}
\newcommand{\mul}{{\rm mul}}
\newcommand{\F}{{\rm F}}
\newcommand{\dist}{{\rm dist}}
\newcommand{\SU}{{\rm SU}}

\newcommand{\rre}{{\rm e}}
\newcommand{\rnu}{{\rm nu}}

\newcommand{\rrq}{{\rm q}}
\newcommand{\rH}{{\rm H}}
\newcommand{\sa}{{\rm sa}}

\newcommand{\tot}{\mathbin{\text{\footnotesize\textcircled{\raisebox{-.6pt}{\tiny \sf T}}}}}

\newcommand{\bd}{{\bf d}}

\newcommand{\id}{{\text{\rm id}}}

\newcommand{\pa}{{\|}}

\newcommand{\GNS}{{\rm GNS}}

\newcommand{\spn}{{\rm span}}

\newcommand{\fm}{{\rm fm}}

\newcommand{\Ce}{{\mathbb C}}
\renewcommand{\Re}{{\mathbb R}}
\newcommand{\Ze}{{\mathbb Z}}
\newcommand{\Ne}{{\mathbb N}}
\newcommand{\Te}{{\mathbb T}}

\allowdisplaybreaks

\begin{document}

\title[Ergodic actions of compact quantum groups]
{Compact quantum metric spaces and ergodic actions of compact quantum groups}

\author{Hanfeng Li}
\thanks{Partially supported by NSF Grant DMS-0701414.}
\address{Department of Mathematics \\
SUNY at Buffalo \\
Buffalo, NY 14260, USA} \email{hfli@math.buffalo.edu}
\date{September 8, 2008}

\begin{abstract}
We show that for any co-amenable compact quantum group
$A=C(\cG)$ there exists a unique compact Hausdorff topology on
the set $\EA(\cG)$ of isomorphism classes of ergodic actions of $\cG$
such that the following holds: for any continuous field of ergodic
actions of $\cG$ over a locally compact Hausdorff space $T$ the
map $T\rightarrow \EA(\cG)$ sending each $t$ in $T$ to the isomorphism
class of the fibre at $t$ is continuous if and only if the
function counting the multiplicity of $\gamma$ in each fibre is
continuous over $T$ for every equivalence class $\gamma$ of
irreducible unitary representations of  $\cG$. Generalizations for
arbitrary compact quantum groups are also obtained.
In the case $\cG$ is a compact group, the restriction of this
topology on the subset of isomorphism classes of ergodic actions
of full multiplicity coincides with the topology coming from the
work of Landstad and Wassermann.
Podle\'s spheres are shown to be continuous in the natural parameter
as ergodic actions of the quantum $\SU(2)$ group.
We also introduce a notion of regularity for quantum metrics on $\cG$, and show
how to construct a quantum metric from any ergodic action of $\cG$,
starting from a regular quantum metric on $\cG$.
Furthermore, we introduce a quantum Gromov-Hausdorff
distance between ergodic actions of $\cG$ when $\cG$ is separable and show that it induces the above topology.
\end{abstract}

\maketitle

\section{Introduction}
\label{intro:sec}

An ergodic action of a compact group $G$ on a unital $C^*$-algebra
$B$  is a  strongly continuous action of $G$ on $B$ such that
the fixed point algebra consists only of scalars.
For an irreducible representation of $G$ on a Hilbert space $H$, the conjugate action
of $G$ on the algebra $B(H)$ is ergodic. On the other hand, ergodic actions
of $G$ on commutative unital $C^*$-algebras correspond exactly to translations on
homogeneous spaces of $G$. Thus the theory of ergodic actions of $G$ connects both the representation
theory and the study of homogeneous spaces.
See \cite{HLS, Landstad, OPT, Wassermann89, Wassermann, Wassermann88} and references therein.

Olesen, Pedersen, and Takesaki classified faithful ergodic actions of an abelian  compact group
as skew-symmetric bicharacters on the dual group \cite{OPT}. Landstad and Wassermann generalized
their result independently to show that ergodic actions of full multiplicity of an arbitrary compact group
$G$ are classified by equivalence classes of dual cocycles \cite{Landstad, Wassermann}. However,
the general case is quite difficult--so far there is no classification of (faithful) ergodic actions of compact groups,
not to mention compact quantum groups.
In this paper
we are concerned with  topological properties of  the whole set $\EA(G)$ of isomorphism classes of ergodic actions of
a compact group $G$, and more generally, the set $\EA(\cG)$ of isomorphism classes of ergodic actions of
a compact quantum group $A=C(\cG)$.

As a consequence of their classification, Olesen, Pedersen, and Takesaki showed that the set
of isomorphism classes of faithful ergodic actions of an abelian  compact group
has a natural abelian compact group structure. From the work of Landstad and Wassermann,
the set $\EA(G)_{\fm}$ of ergodic actions of full multiplicity of an arbitrary compact group $G$
also carries a natural compact Hausdorff topology.

There are many ergodic actions not of full multiplicity, such as conjugation actions associated
to irreducible representations and actions corresponding to translations on homogeneous
spaces (unless $G$ is finite or the homogeneous space is $G$ itself).
In the physics literature concerning string theory and quantum field theory, people
talk about fuzzy spheres, the matrix algebras $M_n(\Ce)$, converging to the two-sphere
$S^2$ (see the introduction of \cite{Rieffel01} and references therein). One important feature of
this convergence is that each term carries an ergodic action of $\SU(2)$, which is used in the construction
of this approximation of $S^2$ by fuzzy spheres. Thus if one wants to give a concrete mathematical
foundation for this convergence, it is desirable to include the $\SU(2)$ symmetry. However, none of
these actions involved are of full multiplicity, and hence the topology of Landstad and Wassermann does not
apply here.

For compact quantum groups there are even more interesting
examples of ergodic actions \cite{Wang99}. Podle\'s introduced a
family of quantum spheres $S^2_{qt}$, parameterized by a compact
subset $T_q$ of the real line, as ergodic actions of the quantum
$\SU(2)$ group $\SU_q(2)$ satisfying certain spectral conditions
\cite{Podles87}. These quantum spheres carry interesting
non-commutative differential geometry \cite{DLPS, DS}. One also
expects that Podle\'s quantum spheres are continuous in the
natural parameter $t$ as ergodic actions of $\SU_q(2)$.

Continuous change of $C^*$-algebras is usually described qualitatively as continuous
fields of $C^*$-algebras over locally compact Hausdorff spaces \cite[Chapter 10]{Dixmier}.
There is no difficulty to formulate the equivariant version--continuous fields
of actions of compact groups \cite{Rieffel89} or even compact quantum
groups (see Section~\ref{cont field:sec} below). Thus if there is any natural topology on $\EA(\cG)$, the relation with
continuous fields of ergodic actions should be clarified.

One distinct feature of the theory of compact quantum groups is
that there is a full compact quantum group and a reduced compact
group associated to each compact quantum group $\cG$, which may
not be the same. A compact quantum group $\cG$ is called
co-amenable if the full and reduced
compact quantum groups coincide. This is the case for compact
groups and $\SU_q(N)$ (for $0<|q|<1$). Our result is simplified in
such case. Denote by $\hat{\cG}$ the set of equivalence classes of
irreducible unitary representations of $\cG$. For each ergodic
action of $\cG$, one can talk about the multiplicity of each
$\gamma \in \hat{\cG}$ in this action \cite{Podles95}, which is
known to be finite for the compact group case by \cite{HLS} and, for
the compact quantum group case by \cite{Boca}.

\begin{theorem} \label{main1:thm}
Let $\cG$ be a co-amenable compact quantum group. Then $\EA(\cG)$ has a
unique compact Hausdorff topology such that the following holds:
for any continuous field of ergodic actions of $\cG$ over a
locally compact Hausdorff space $T$ the map $T\rightarrow \EA(\cG)$
sending each $t\in T$ to the isomorphism class of the fibre at $t$
is continuous if and only if the function counting the
multiplicity of $\gamma$ in each fibre is continuous over $T$ for
each $\gamma$ in
$\hat{\cG}$.
\end{theorem}

In particular, fuzzy spheres converge to $S^2$ as ergodic actions
of $\SU(2)$ (see \cite[Example 10.12]{Li10}). Podle\'s quantum
spheres are also continuous as ergodic actions of $\SU_q(2)$:

\begin{theorem} \label{cont of Podles:thm}
Let $q$ be a real number with $0<|q|<1$, and let
$T_q$ be the parameter space of Podle\'s quantum spheres.
The map $T_q\rightarrow \EA(\SU_q(2))$ sending $t$ to the isomorphism class
of
$S^2_{qt}$ is continuous.
\end{theorem}

When $\cG$ is not co-amenable, the more appropriate object to study is a certain
quotient space of $\EA(\cG)$.
To each ergodic action of $\cG$, there is an associated
full ergodic action and an associated reduced ergodic action (see Section~\ref{f and r actions:sec} below),
which are always isomorphic when $\cG$ is co-amenable.
Two ergodic actions are said to be equivalent if the associated full (reduced resp.) actions are isomorphic.
Denote by $\EA^{\sim}(\cG)$ the quotient space of $\EA(\cG)$ modulo this equivalence relation.
We also have to deal with semi-continuous fields of ergodic actions in the general case.

\begin{theorem} \label{main2:thm}
Let $\cG$ be a compact quantum group. Then $\EA^{\sim}(\cG)$ has a
unique compact Hausdorff topology such that the following holds:
for any semi-continuous field of ergodic actions of $\cG$ over a
locally compact Hausdorff space $T$ the map $T\rightarrow \EA^{\sim}(\cG)$
sending each $t\in T$ to the equivalence class of the fibre at $t$
is continuous if and only if the function counting the
multiplicity of $\gamma$ in each fibre is continuous over $T$ for
each $\gamma$ in
$\hat{\cG}$.
\end{theorem}

Motivated partly by the need to give a mathematical foundation for
various approximations in the string theory literature, such as the approximation of $S^2$ by
fuzzy spheres in above, Rieffel initiated the theory of compact quantum  metric spaces
and quantum Gromov-Hausdorff distances \cite{Rieffel00, Rieffel03}. As the information
of the metric on a compact metric space $X$ is encoded in the Lipschitz seminorm on $C(X)$,
a quantum metric on (the non-commutative space corresponding to) a unital $C^*$-algebra $B$
is a (possibly $+\infty$-valued) seminorm on $B$ satisfying suitable conditions. The seminorm is called a Lip-norm.
Given a length function on a compact group $G$, Rieffel showed how to induce a quantum metric
on (the $C^*$-algebra carrying) any ergodic action of $G$ \cite{Rieffel98b}. We find that the
right generalizations of length functions for a compact quantum
group $A=C(\cG)$ are Lip-norms on $A$ being finite on the
algebra $\cA$ of regular functions, which we call regular
Lip-norms. Every separable co-amenable $A$ has a bi-invariant
regular Lip-norm (Corollary~\ref{biinvariant:coro}). Then we have
the following generalization of Rieffel's construction (see
Section~\ref{prelim:sec} below for more detail on the notation), answering
a question Rieffel raised at the end of Section 3 in \cite{Rieffel03}.

\begin{theorem} \label{induced Lip-norm:thm}
Suppose that $\cG$ is a co-amenable compact quantum group and
$L_{A}$ is a regular Lip-norm on $A$. Let
$\sigma:B\rightarrow B\otimes A$ be an ergodic action
of $\cG$ on a unital $C^*$-algebra $B$.
Define a (possibly $+\infty$-valued) seminorm
on $B$ via
\begin{eqnarray} \label{induced Lip-norm:eq}
L_{B}(b)=\sup_{\varphi\in S(B)}L_{A}(b*\varphi)
\end{eqnarray}
for all $b\in B$, where $S(B)$ denotes the state space of $B$ and $b*\varphi=(\varphi \otimes \id)(\sigma(b))$.
Then $L_{B}$
is finite on the algebra $\cB$ of regular functions
and is a Lip-norm on $B$ with $r_{B}\le 2r_{A}$, where
$r_{B}$ and $r_{A}$ are the radii of $B$ and $A$ respectively.
\end{theorem}

As an important step towards establishing  a mathematical
foundation for various convergence in the string theory
literature, such as the convergence of fuzzy spheres to $S^2$,
Rieffel introduced a quantum Gromov-Hausdorff distance
$\dist_{\rrq}$ between compact quantum metric spaces and showed,
among many properties of $\dist_{\rrq}$, that the fuzzy spheres
converge to $S^2$ under $\dist_{\rrq}$ when they are all endowed
with the quantum metrics induced from the ergodic actions of
$\SU(2)$ for a fixed length function on $\SU(2)$ \cite{Rieffel01}.
Two generalizations of $\dist_{\rrq}$ are introduced in
\cite{Kerr} and \cite{Li12} in order to distinguish the algebra
structures (see also \cite{KL}). However, none of these quantum
distances distinguishes the group symmetries. That is, there exist
non-isomorphic ergodic actions of a compact group such that
quantum distances between the compact quantum metric spaces
induced by these ergodic actions are zero (see
Example~\ref{tori:ex} below). One of the features of our quantum
distances in \cite{Li12, Li10} is that they can be adapted easily
to take care of other algebraic structures. Along the lines in
\cite{Li12, Li10}, we introduce a quantum distance $\dist_{\rre}$
(see Definition~\ref{distance:def} below)
between the compact quantum metric spaces coming from
ergodic action of $\cG$
as in Theorem~\ref{induced Lip-norm:thm}.
This distance distinguishes the ergodic actions:

\begin{theorem} \label{distance vs topology:thm}
Let $\cG$ be a co-amenable compact quantum group with a fixed
left-invariant regular Lip-norm $L_{A}$ on $A$.
Then $\dist_{\rre}$ is a metric on $\EA(\cG)$ inducing the topology in Theorem~\ref{main1:thm}.
\end{theorem}

The organization of this paper is as follows. In
section~\ref{prelim:sec} we recall some basic definitions and
facts about compact quantum groups, their actions, and compact
quantum metric spaces. Associated full and reduced actions are
discussed in Section~\ref{f and r actions:sec}. The topologies on
$\EA(\cG)$ and $\EA^{\sim}(\cG)$ are introduced in
Section~\ref{ergodic:sec}. We also prove that $\EA^{\sim}(\cG)$ is
compact Hausdorff there. In Section~\ref{cont field:sec} we
clarify the relation between semi-continuous fields of ergodic
actions and the topology introduced in Section~\ref{ergodic:sec}.
This completes the proofs of Theorems~\ref{main1:thm} and
\ref{main2:thm}. The continuity of Podle\'s quantum spheres is
discussed in Section~\ref{Podles:sec}. In Section~\ref{fm:sec} we
show that the topology of Landstad and Wassermann on $\EA(G)_{\fm}$ for a compact group $G$ is
simply the relative topology of $\EA(G)_{\fm}$ in $\EA(G)$.
Theorems~\ref{induced Lip-norm:thm} and \ref{distance vs
topology:thm} are proved in Sections~\ref{induced:sec} and
\ref{distance:sec} respectively.

\begin{acknowledgments} I am grateful to Florin Boca,
George Elliott, Marc Rieffel, Shuzhou Wang, and Wei Wu
for valuable discussions. I thank
Sergey Neshveyev for help on the proof of Lemma~\ref{approx e:lemma},
and thank Magnus Landstad for addressing a question on the topology of
$\EA(G)_{\fm}$ for compact groups. I also would like to thank the referee 
for several useful comments.
\end{acknowledgments}

\section{Preliminaries}
\label{prelim:sec}

In this section we collect some definitions and facts about compact quantum groups and
compact quantum metric spaces.

Throughout this paper we use $\otimes$ for the spatial tensor product of $C^*$-algebras,
and $\odot$ for the algebraic tensor product of vector spaces. $A=C(\cG)$ will be
a compact quantum group.

\subsection{Compact quantum groups and actions}
\label{actions:sub}

We recall first some definitions and facts about compact quantum groups.
See \cite{MV, Wor87, Wor95} for more detail.

A {\it compact quantum group} $(A, \Phi)$ is a unital $C^*$-algebra
$A$ and a unital $*$-homomorphism $\Phi:A\rightarrow
A\otimes A$ such that $(\id\otimes \Phi)\Phi=(\Phi\otimes
\id)\Phi$ and that both $\Phi(A)(1_{A}\otimes A)$ and
$\Phi(A)(A\otimes 1_{A})$ are dense in $A\otimes A$, where $\Phi(A)(1_A\otimes A)$ ($\Phi(A)(A\otimes 1_A)$ resp.)
denotes
the linear span of $\Phi(a)(1_A\otimes a')$ ($\Phi(a)(a'\otimes 1_A)$ resp.) for $a, a'\in A$.
We shall write $A$ as $C(\cG)$ and say that $\cG$ is the compact
quantum group. The Haar measure is the unique state $h$ of $A$
such that $(\id \otimes h)\Phi=(h \otimes \id)\Phi=h$.

A unitary representation of $\cG$ on a Hilbert space $H$ is a
unitary $u\in M(K(H)\otimes A)$ such that $(\id \otimes
\Phi)(u)=u_{12}u_{13}$, where $K(H)$ is the algebra of compact
operators, $M(K(H)\otimes A)$ is the multiplier algebra of
$K(H)\otimes A$, and we use the leg numbering notation
\cite[page 385]{PW}. When $H$ is finite dimensional, $u^c$
denotes the contragradient representation acting on the conjugate
Hilbert space of $H$. For unitary representations $v$ and $w$ of
$\cG$, the tensor product representation $v\tot w$ is
defined as $v_{13}w_{23}$ in the leg numbering notation. Denote by
$\hat{\cG}$ the set of equivalence classes of irreducible unitary
representations of $\cG$. For each $\gamma \in \hat{\cG}$ fix
$u^{\gamma}\in \gamma$ acting on $H_{\gamma}$ and an orthonormal
basis in $H_{\gamma}$. Each $H_{\gamma}$ is finite
dimensional. Denote by $d_{\gamma}$ the dimension of
$H_{\gamma}$. Then we may identify $B(H_{\gamma})$ with
$M_{d_{\gamma}}(\Ce)$, and hence $u^{\gamma}\in
M_{d_{\gamma}}(A)$. Denote by $A_{\gamma}$ the linear span of
$(u^{\gamma}_{ij})_{ij}$. Then
\begin{eqnarray} \label{action in A:eq}
A^*_{\gamma}=A_{\gamma^c}, \quad \Phi(A_{\gamma})\subseteq A_{\gamma}\odot A_{\gamma},
\end{eqnarray}
and
$\cA:=\oplus_{\gamma\in \hat{\cG}}A_{\gamma}$ is the algebra of
{\it regular functions} in $A$.
For any $1\le i, j\le d_{\gamma}$ denote by $\rho^{\gamma}_{ij}$ the unique element in
$A'$ such that
\begin{eqnarray} \label{rho:eq}
\rho^{\gamma}_{ij}(u^{\beta}_{sk})=\delta_{\gamma \beta}\delta_{is}\delta_{jk}
\end{eqnarray}
(the existence of such $\rho^{\gamma}_{ij}$ is guaranteed by \cite[Theorem 5.7]{Wor87}).
Moreover $\rho^{\gamma}_{ij}$ is of the form $h(\cdot a)$ for some $a\in A_{\gamma^c}$.
Denote $\sum_{1\le i\le d_{\gamma}}\rho^{\gamma}_{ii}$ by $\rho^{\gamma}$.
Denote the class of the trivial representations of $\cG$ by $\gamma_0$.

There exist a full compact quantum group $(A_u, \Phi_u)$ and
a reduced compact group $(A_r, \Phi_r)$ whose algebras of
regular functions and restrictions of comultiplications are the same as
$(\cA, \Phi|_{\cA})$. The quantum group $\cG$ is  said to be {\it co-amenable} if the canonical
surjective homomorphism $A_u\rightarrow A_r$ is an isomorphism
\cite[Definition 6.1]{Banica} \cite[Theorem 3.6]{BMT}.
There is a unique $*$-homomorphism $e:\cA\rightarrow \Ce$ such that
$(e\otimes \id)\Phi=(\id \otimes e)\Phi=\id$ on $\cA$, which is called the {\it counit}.
The quantum group $\cG$ is co-amenable exactly if it has bounded counit and faithful Haar measure
\cite[Theorem 2.2]{BMT}.

Next we recall some facts about actions of $\cG$. See \cite{Boca, Podles95} for detail.

\begin{definition}\cite[Definition 1.4]{Podles95} \label{action:def}
A (left) {\it action} of $\cG$ on a unital $C^*$-algebra $B$ is a unital
$*$-homomorphism $\sigma:B\rightarrow B\otimes A$ such that
\begin{enumerate}
\item $(\id\otimes \Phi)\sigma=(\sigma\otimes \id)\sigma$,

\item $\sigma(B)(1_{B}\otimes A)$ is dense in $B\otimes A$.

\end{enumerate}
The {\it fixed point algebra} of $\sigma$ is $B^{\sigma}=\{b\in B: \sigma(b)=b\otimes 1_{A}\}$.
The action $\sigma$ is {\it ergodic} if $B^{\sigma}=\Ce 1_{B}$.
\end{definition}

\begin{remark} \label{injective:remark}
When $A$ has bounded counit, the proof of \cite[Lemma 1.4.(a)]{EN}
shows that $(\id\otimes e)\sigma=\id$ on $B$ and that $\sigma$ is injective.
\end{remark}

Let $\sigma:B\rightarrow B\otimes A$ be an action
of $\cG$ on a unital $C^*$-algebra $B$.
For any $b\in B, \varphi\in B'$ and $\psi\in A'$  set
\begin{eqnarray} \label{slice:eq}
b*\varphi=(\varphi\otimes \id)(\sigma(b)), \quad
\psi*b=(\id \otimes \psi)(\sigma(b)).
\end{eqnarray}
Also set
$E^{\gamma}_{ij}, E^{\gamma}:B\rightarrow B$ by
\begin{eqnarray} \label{E_ij:eq}
E^{\gamma}_{ij}=(\id \otimes \rho^{\gamma}_{ij})\sigma, \quad E^{\gamma}=(\id \otimes \rho^{\gamma})\sigma.
\end{eqnarray}
Then \cite[page 98]{Boca}
\begin{eqnarray} \label{idempotent:eq}
E^{\beta}E^{\gamma}=\delta_{\beta\gamma}E^{\gamma}.
\end{eqnarray}
Set
\begin{eqnarray}
B_{\gamma}=E^{\gamma}(B),\quad \cB=\oplus_{\gamma\in \hat{\cG}}B_{\gamma}.
\end{eqnarray}
Then $\cB$ is a dense $*$-subalgebra of $B$ \cite[Theorem 1.5]{Podles95} \cite[Lemma 11, Proposition 14]{Boca}
(the ergodicity condition in \cite{Boca} is not used in Lemma 11 and Proposition 14 therein), which we shall
call the algebra of {\it regular functions} for $\sigma$.
Moreover,
\begin{eqnarray} \label{*-algebra:eq}
B^*_{\gamma}=B_{\gamma^c}, \quad B_{\alpha}B_{\beta}\subseteq \sum_{\gamma \preceq \alpha\tot \beta}B_{\gamma},\quad
\sigma(B_{\gamma})\subseteq B_{\gamma}\odot A_{\gamma}.
\end{eqnarray}
There exist a set $J_{\gamma}$
and a linear basis $\sS_{\gamma}=\{e_{\gamma ki}: k\in J_{\gamma}, 1\le i\le d_{\gamma}\}$
of $B_{\gamma}$ \cite[Theorem 1.5]{Podles95} such that
\begin{eqnarray} \label{action:eq}
\sigma(e_{\gamma ki})=\sum_{1\le j\le d_{\gamma}}e_{\gamma kj}\otimes u^{\gamma}_{ji}.
\end{eqnarray}
The {\it multiplicity} $\mul(B, \gamma)$ is defined as the cardinality
of $J_{\gamma}$, which  does not depend on the choice of
$\sS_{\gamma}$. Conversely, given a unital $*$-homomorphism
$\sigma:B\rightarrow B\otimes A$ for a unital $C^*$-algebra
$B$, if there exist a set $J'_{\gamma}$ and a set of linearly
independent elements $\sS'_{\gamma}=\{e_{\gamma ki}: k\in
J'_{\gamma}, 1\le i\le d_{\gamma}\}$ in $B$ satisfying
(\ref{action:eq}) for each $\gamma \in \hat{\cG}$ such that the
linear span of $\cup_{\gamma \in \hat{\cG}} \sS_{\gamma}$ is dense
in $B$, then $\sigma$ is an action of $\cG$ on $B$ and $\spn
\sS'_{\gamma}\subseteq B_{\gamma}$ for each $\gamma \in
\hat{\cG}$ \cite[Corollary 1.6]{Podles95}. In this event, if
$|J'_{\gamma}|$ or $\mul(B, \gamma)$ is finite, then
$B_{\gamma}=\spn \sS'_{\gamma}$.

We have that $B_{\gamma_0}=B^{\sigma}$ and that $E:=E^{\gamma_0}$ is a conditional expectation
from $B$ onto $B^{\sigma}$ \cite[Lemma 4]{Boca}. When $\sigma$ is ergodic, $E=\omega(\cdot)1_{B}$
for the unique $\sigma$-invariant state $\omega$ on $B$.

\subsection{Compact quantum metric
spaces}\label{cqms:sub}

In this subsection  we recall some facts about compact quantum
metric spaces \cite{Rieffel98b, Rieffel00, Rieffel03}.
Though Rieffel has set up his theory in the general framework of
order-unit spaces, we shall need it only for $C^*$-algebras. See
the discussion preceding Definition 2.1 in \cite{Rieffel00} for
the reason of requiring the reality condition (\ref{real:eq})
below.

\begin{definition}\cite[Definition 2.1]{Rieffel00} \label{C*CQM:def}
By a \emph{$C^*$-algebraic compact quantum metric space} we mean a
pair $(B, L)$ consisting of a unital $C^*$-algebra
$B$ and  a (possibly $+\infty$-valued) seminorm $L$ on
$B$ satisfying
 the \emph{reality condition}
\begin{eqnarray} \label{real:eq}
L(b)&=& L(b^*)
\end{eqnarray}
for all $b\in B$, such that $L$ vanishes on
$\Ce 1_{B}$ and the metric $\rho_L$ on the state space $S(B)$
defined by
\begin{eqnarray} \label{Lip to dist:eq}
\rho_L(\varphi, \psi)=\sup_{L(b)\le 1} |\varphi(b)-\psi(b)|
\end{eqnarray}
 induces the weak-$*$ topology.
The \emph{radius} $r_{B}$  of $(B, L)$ is defined to be the radius of $(S(B), \rho_L)$.
We say that $L$ is a \emph{Lip-norm}.
\end{definition}

Note that $L$ must in fact vanish precisely on $\Ce 1_{B}$
and take finite values on a dense subspace of $B$.

Let $B$ be a unital $C^*$-algebra and let $L$ be a
(possibly $+\infty$-valued) seminorm on $B$ vanishing on
$\Ce 1_{B}$. Then $L$ and $\pa \cdot \pa$ induce (semi)norms ${\tilde
L}$ and $\pa \cdot\pa^{\sim}$ respectively on the quotient space
$\tilde{B}=B/(\Ce 1_{B})$.

\begin{notation} \label{ball:notation}
Let
\begin{eqnarray*}
\cE(B):=\{b\in B_{\sa}:L(b)\le 1\}.
\end{eqnarray*}
For any $r\ge 0$, let
\begin{eqnarray*}
\cD_r(B):=\{b\in B_{\sa}: L(b)\le 1, \pa b\pa \le r\}.
\end{eqnarray*}
\end{notation}

Note that the definitions of $\cE(B)$ and $\cD_r(A)$ use $B_{\sa}$ instead of $B$.
The main criterion for when a seminorm $L$ is a Lip-norm is the following:

\begin{proposition}\cite[Proposition 1.6, Theorem 1.9]{Rieffel98b}
\label{criterion of Lip:prop}
Let $B$ be a unital $C^*$-algebra and let
$L$ be a (possibly $+\infty$-valued) seminorm on $B$
satisfying the reality condition (\ref{real:eq}).
Assume that $L$ takes finite values on a dense
subspace of $B$, and that $L$ vanishes on $\Ce 1_{B}$.
Then $L$ is a Lip-norm if and only if

\noindent \, \, \, \, \,\, \, \,
(1) there is a constant $K\ge 0$ such that $\pa \cdot\pa^{\sim}\le K
\tilde{L}$ on $\tilde{B}$;

\noindent
and \, \, \, \,
(2) for any $r\ge 0$, the  ball $\cD_r(B)$
is totally bounded in $B$ for $\pa \cdot \pa$;

\noindent \, \, \, \, \,
or
(2') for some $r> 0$, the  ball $\cD_r(B)$
is totally bounded in $B$ for $\pa \cdot \pa$.

In this event, $r_{B}$ is exactly the minimal $K$ such
that $\pa \cdot\pa^{\sim}\le K \tilde{L}$ on
$(\tilde{B})_{\sa}$.
\end{proposition}

\section{Full and reduced actions}
\label{f and r actions:sec}

In this section we discuss full and reduced actions associated
to actions of $\cG$. We will use the notation in subsection~\ref{actions:sub}
freely.
Throughout this section, $\sigma$ will be
an action of $\cG$ on a unital $C^*$-algebra $B$.

\begin{lemma} \label{faithful:lemma}
The conditional expectation $E=(\id \otimes h)\sigma$ is faithful on $\cB$. If $A$ is co-amenable, then $E$ is faithful on $B$.
\end{lemma}
\begin{proof} Suppose that $E(b)=0$ for some positive $b$ in $\cB$.
Then for any $\varphi\in S(B)$ we have
$0=\varphi(E(b))=h(b*\varphi)$. Observe that $b*\varphi$ is in
$\cA$ and is positive. By the faithfulness of $h$ on $\cA$
\cite[Theorem 4.2]{Wor87}, $b*\varphi=0$. Then $(\varphi\otimes
\psi)(\sigma(b))=\psi(b*\varphi)=0$ for all $\varphi\in S(B)$
and $\psi\in S(A)$. Since product states separate points of
$B\otimes A$ \cite[Lemma T.5.9 and Proposition T.5.14]{WO}, $\sigma(b)=0$. From (\ref{idempotent:eq}) one
sees that $b\in \{\phi*b:\phi \in A'\}$. Therefore $b=0$. The
second assertion is proved similarly, in view of
Remark~\ref{injective:remark}.
\end{proof}

For actions $\sigma_i:B_i\rightarrow B_i\otimes A$ of $\cG$ on $B_i$ for $i=1,2$,
a unital $*$-homomorphism $\theta:B_1\rightarrow B_2$ is said to be {\it equivariant}
(with respect to $\sigma_1$ and $\sigma_2$)
if $\sigma_2\circ \theta=(\theta\otimes \id)\circ \sigma_1$.

\begin{lemma} \label{equivariant:lemma}
Let $\theta:B_1\rightarrow B_2$ be a unital $*$-homomorphism
equivariant with respect
to actions $\sigma_1, \sigma_2$ of $\cG$ on $B_1$ and $B_2$. Then
\begin{eqnarray}
\label{equi:eq} E^{\gamma}\circ \theta=\theta \circ E^{\gamma}, \\
\label{equi spetral:eq} \theta((B_1)_{\gamma})\subseteq (B_2)_{\gamma}
\end{eqnarray}
for all $\gamma \in \hat{\cG}$. The map $\theta$ is surjective if and only
if $\theta(\cB_1)=\cB_2$. The map $\theta$ is injective on $\cB_1$ if and only if $\theta$ is
injective on $B^{\sigma_1}_1$.
\end{lemma}
\begin{proof} One has
\begin{eqnarray*}
E^{\gamma}\circ \theta &=&(\id \otimes \rho^{\gamma})\circ \sigma_2\circ \theta=
(\id \otimes \rho^{\gamma})\circ (\theta \otimes \id)\circ \sigma_1 \\
&=&
\theta \circ (\id \otimes \rho^{\gamma})\circ \sigma_1=\theta \circ E^{\gamma},
\end{eqnarray*}
which proves (\ref{equi:eq}). The formula (\ref{equi spetral:eq}) follows from (\ref{equi:eq}).

Since $\cB_2$ is dense in $B_2$, if $\theta(\cB_1)=\cB_2$, then
$\theta$ is surjective. Conversely, suppose that $\theta$ is
surjective. Applying both sides of (\ref{equi:eq}) to $B_1$ we
get $\theta((B_1)_{\gamma})=(B_2)_{\gamma}$ for each
$\gamma\in \hat{\cG}$. Thus $\theta(\cB_1)=\cB_2$.

Since $B^{\sigma_1}_1\subseteq \cB_1$, if $\theta$ is injective on $\cB_1$, then  $\theta$ is
injective on $B^{\sigma_1}_1$. Conversely, suppose that $\theta$ is injective on
$B^{\sigma_1}_1$. Let $b\in \cB_1\cap \ker \theta$.
Then $\theta(E(b^*b))=E(\theta(b^*b))=0$. Thus
$E(b^*b)\in \ker \theta$. By assumption we have $E(b^*b)=0$.
Then $b=0$ by Lemma~\ref{faithful:lemma}.
\end{proof}

\begin{proposition} \label{universal action:prop}
The $*$-algebra $\cB$ has a universal $C^*$-algebra $B_u$. The canonical $*$-homomorphism $\cB\rightarrow B_u$
is injective. Identify $\cB$ with its canonical image in $B_u$.
The unique $*$-homomorphism $\sigma_u:B_u\rightarrow B_u\otimes A$
extending $\cB\overset{\sigma|_{\cB}}\rightarrow  \cB\odot \cA\hookrightarrow B_u\otimes A$
is an action of $\cG$ on
$B_u$. Moreover, the unique $*$-homomorphism $\pi_u:B_u \rightarrow B$
extending the embedding $\cB\rightarrow B$
is equivariant, and
the algebra of regular functions for $\sigma_u$ is $\cB$.
\end{proposition}
\begin{proof}
Let $\gamma\in \hat{\cG}$ and let $\sS_{\gamma}$ be a linear basis
of $B_{\gamma}$ satisfying (\ref{action:eq}).
Set
$b_k=\sum^{d_{\gamma}}_{i=1}e_{\gamma ki}e^*_{\gamma ki}$. Then
\begin{eqnarray*}
\sigma(b_k)
=\sum_{1\le i,j, s\le d_{\gamma}}e_{\gamma kj}e^*_{\gamma ks}\otimes u^{\gamma}_{ji}{u^{\gamma}_{si}}^*
=\sum_{1\le j, s\le d_{\gamma}}e_{\gamma kj}e^*_{\gamma ks}\otimes \delta_{js}1_{A}
=b_k\otimes 1_{A}.
\end{eqnarray*}
Thus $b_k\in B^{\sigma}$.
Note that $B^{\sigma}$ is a $C^*$-subalgebra of $B$.
So $\pa \pi(\cdot)\pa \le \pa \cdot\pa $ on $B^{\sigma}$ for any $*$-representation
$\pi$ of $\cB$.
Consequently, $\pa \pi(e_{\gamma ki})\pa^2\le \pa b_k\pa$ for any $*$-representation
$\pi$ of $\cB$.
Thus for any $c\in \cB$ there is some $\lambda_c \in \Re$ such that $\pa \pi(c)\pa \le \lambda_c$
 for any $*$-representation
$\pi$ of $\cB$.  Therefore $\cB$ has a universal $C^*$-algebra $B_u$ with the canonical
$*$-homomorphism $\phi:\cB\rightarrow B_u$. Then there is a unique $*$-homomorphism
$\pi_u:B_u\rightarrow B$ such that $\pi_u\circ \phi$ is the canonical embedding $\iota:\cB\hookrightarrow B$.
Since $\iota$ is injective, so is $\phi$. Thus we may identify $\cB$ with $\phi(\cB)$.

Denote by $\sigma_u$ the unique $*$-homomorphism $B_u\rightarrow B_u\otimes A$
extending the $*$-homomorphism
$\cB\overset{\sigma|_{\cB}}\rightarrow  \cB\odot \cA\hookrightarrow B_u\otimes A$.
According to the characterization of actions of $\cG$ in terms of elements
satisfying (\ref{action:eq}) in subsection~\ref{actions:sub}, $\sigma_u$
is an action of $\cG$ on $B_u$ and $\cB\subseteq \cB_u$.
Since $\sigma\circ \pi_u$ and $(\pi_u\otimes \id)\circ \sigma_u$
coincide on $\cB$, they
also coincide on $B_u$.
Thus $\pi_u$ is equivariant.

Since $B^{\sigma}$ is closed and $E(\cB)=B^{\sigma}$, $B^{\sigma_u}_u=E(B_u)=B^{\sigma}$.
Thus $\pi_u$ is injective on $B^{\sigma_u}_u$. By Lemma~\ref{equivariant:lemma}
the map $\theta$ is injective on $\cB_u$.
Let $b\in \cB_u$. Then $\pi_u(b)\in \cB$ by Lemma~\ref{equivariant:lemma}, and
hence $\pi_u(b-\pi_u(b))=0$. Therefore $b=\pi_u(b)\in \cB$. This proves
$\cB_u=\cB$ as desired.
\end{proof}

We refer the reader to \cite{Lance} for basics on Hilbert $C^*$-modules.
Since $E: B\rightarrow B^{\sigma}$ is a conditional expectation,
$B$ is a right semi-inner-product $B^{\sigma}$-module
with the inner product $\left<\cdot, \cdot\right>_{B^{\sigma}}$ given by
$\left<x, y\right>_{B^{\sigma}}=E(x^*y)$ \cite[page 7]{Lance}.
Denote by $H_{B}$ the completion, and by
$\pi_r$ the associated representation of $B$ on $H_{B}$.
Denote $\pi_r(B)$ by $B_r$.

\begin{proposition} \label{reduced:prop}
There exists a unique $*$-homomorphism $\sigma_r:B_r\rightarrow B_r\otimes A$
such that
\begin{eqnarray} \label{reduced:eq}
\sigma_r\circ \pi_r=(\pi_r\otimes \id)\circ \sigma.
\end{eqnarray}
The homomorphism $\sigma_r$ is injective and is an action of $\cG$ on $B_r$.
The map $\pi_r$ is equivariant, and is injective on $\cB$.
The algebra of regular functions for $\sigma_r$ is $\pi_r(\cB)$.
\end{proposition}
\begin{proof} The uniqueness of such $\sigma_r$ follows from the surjectivity of $\pi_r$.
Consider the right Hilbert $(B^{\sigma}\otimes A)$-module $H_{B}\otimes A$.
Denote by $B(H_B\otimes A)$
the $C^*$-algebra of adjointable operators of the Hilbert $(B^{\sigma}\otimes A)$-module $H_{B}\otimes A$.
Then $B_r\otimes A\subseteq B(H_B\otimes A)$.
The argument in the proof of \cite[Lemma 3]{Boca} shows that
there is a unitary $U\in B(H_{B}\otimes A)$ satisfying
$U(b\otimes a)=((\pi_r\otimes \id)(\sigma(b)))(1_{B}\otimes a)$ for all $a\in A, b\in B$.
It follows that $U(\pi_r(b)\otimes 1_{A})=((\pi_r\otimes \id)(\sigma(b)))U$
for all $b\in B$.
Thus $U(B_r\otimes 1_{A})U^{-1}\subseteq B_r\otimes A$.
Define $\sigma_r:B_r\rightarrow B_r\otimes A$ by $\sigma_r(b')=
U(b'\otimes 1_{A})U^{-1}$. Then (\ref{reduced:eq}) follows.
Clearly  $\sigma_r$ is injective. Since $\sigma$ is an action
of $\cG$ on $B$, it follows easily that $\sigma_r$ is an action of $\cG$ on $B_r$.
 The equivariance of $\pi_r$ follows from (\ref{reduced:eq}).
By Lemma~\ref{equivariant:lemma},
$\pi_r(\cB)=\cB_r$. It is clear that $\pi_r$ is injective on $B^{\sigma}$. Thus by Lemma~\ref{equivariant:lemma}
the map $\pi_r$ is injective on $\cB$.
\end{proof}

\begin{definition} \label{universal and reduced:def}
We call the action $(B_u, \sigma_u)$ in
Proposition~\ref{universal action:prop} the {\it full action} associated
to $(B, \sigma)$, and call  the action $(B_r, \sigma_r)$ in
Proposition~\ref{reduced:prop} the {\it reduced action} associated
to $(B, \sigma)$. The action $(B, \sigma)$ is said to be {\it
full} ({\it reduced}, {\it co-amenable}, resp.) if $\pi_u$
($\pi_r$, both $\pi_u$ and $\pi_r$, resp.) is an isomorphism.
\end{definition}

\begin{example} \label{actions:eg}
\begin{enumerate}
\item When $B$ is finite dimensional, the action $(B, \sigma)$ is co-amenable. This applies to
the actions constructed in \cite{Wang98} and the adjoint action on $B(H)$ associated to
any finite-dimensional representation of $\cG$ on $H$ \cite[notation after Theorem 2.5]{Wang99}.
\item Consider the Cuntz algebra $\cO_n$ \cite{Cuntz} for an integer $n\ge 2$, that is, the universal $C^*$-algebra
generated by isometries $S_1,{\cdots}, S_n$ satisfying
$\sum^n_{j=1}S_jS^*_j=1$.
Since $\cO_n$ is simple, any action of a compact quantum group on $\cO_n$ is reduced.
Given a compact quantum group $A=C(\cG)$ and an $n$-dimensional unitary representation
$u=(u_{ij})_{ij}$ of $\cG$, one has an action $\sigma$ of $\cG$ on $\cO_n$ determined
by $\sigma(S_i)=\sum^n_{j=1}S_j\otimes u_{ji}$ for all $1\le i\le n$ \cite[Theorem 1]{KNW}.
The regular subalgebra $\cB$ for this action $\sigma$ contains $S_1, {\cdots}, S_n$, thus
$\sigma$ is full and hence is co-amenable, because of the universal property of $\cO_n$.
This kind of actions has been considered for $\cG$ being $\SU_q(2)$ \cite{KNW, Marciniak}, $\SU_q(N)$ \cite{Paolucci},
and $A_u(Q)$ \cite[Section 5]{Wang99}.
\item For the action $\Phi:A\rightarrow A\otimes A$ of $\cG$ on $A$, the $C^*$-algebra
for the associated full action is the $C^*$-algebra of the full quantum group \cite[Section 3]{BMT},
while the $C^*$-algebra for the associated reduced action is the $C^*$-algebra of the reduced quantum group
\cite[Section 2]{BMT}. Thus the action $(A, \Phi)$ is full (reduced resp.) exactly
if $\cG$ is a full (reduced resp.) compact quantum group.
\end{enumerate}
\end{example}


\begin{remark} \label{equivalence:remark}
Having isomorphic $(\cB, \sigma|_{\cB})$ is an equivalence relation
between actions of $\cG$ on unital $C^*$-algebras.
Two actions are equivalent in this sense
exactly if they have isomorphic full actions, exactly if they have
isomorphic reduced actions. If $(A_1, \Phi_1)$ is another compact quantum group
with $(\cA_1, \Phi_1|_{\cA_1})$ isomorphic to $(\cA, \Phi|_{\cA})$, then
$A_1$ has also a natural action on $B_u$. Thus the class of the equivalence classes
of actions of $\cG$ depends only on $(\cA, \Phi|_{\cA})$.
\end{remark}

\begin{proposition} \label{universal:prop}
The following are equivalent:
\begin{enumerate}
\item $\cG$ is co-amenable,

\item every action of $\cG$ on a unital $C^*$-algebra is co-amenable,

\item every ergodic action of $\cG$ on a unital $C^*$-algebra is co-amenable.
\end{enumerate}
\end{proposition}
\begin{proof}
(1)$\Rightarrow$(2). Let $\sigma$ be an action of $\cG$ on a unital
$C^*$-algebra $B$. Then $(B_r, \sigma_r)$ is also
the reduced action associated to $(B_u, \sigma_u)$.
By Lemma~\ref{faithful:lemma} $E$ is faithful on $B_u$.
Thus the canonical homomorphism $B_u\rightarrow B_r$ is
injective, and hence is an isomorphism. Therefore $(B, \sigma)$ is co-amenable.

(2)$\Rightarrow$(3). This is trivial.

(3)$\Rightarrow$(1). This follows from Example~\ref{actions:eg}(3).
\end{proof}

\section{Ergodic actions}
\label{ergodic:sec}

In this section we introduce a topology on the set of isomorphism
classes of ergodic actions of $\cG$ in Definition~\ref{topology on
EA:def} and prove Theorem~\ref{topology on EA:thm}. At the end of
this section we also discuss the behavior of this topology under
taking Cartesian products of compact quantum groups.

\begin{notation} \label{EA:notation}
Denote by $\EA(\cG)$ the set of isomorphism classes of ergodic actions
of $\cG$.
Denote by $\EA^{\sim}(\cG)$
the quotient space of $\EA(\cG)$
modulo the equivalence relation in Remark~\ref{equivalence:remark}.
\end{notation}

What we shall do is to define a topology on $\EA^{\sim}(\cG)$, then pull it back to
a topology on $\EA(\cG)$.
For each $\gamma \in \hat{\cG}$, let $M_{\gamma}$ be the quantum dimension defined
after Theorem 5.4 in \cite{Wor87}. One knows that $M_{\gamma}$ is a positive number no less
than $d_{\gamma}$ and that $M_{\gamma_0}=1$.
Set $N_{\gamma}$ to be the largest integer no bigger than $M^2_{\gamma}/d_{\gamma}$.
Let $(B, \sigma)$ be an ergodic action of $\cG$.
According to \cite[Theorem 17]{Boca}, one has $\mul(B, \gamma)\le N_{\gamma}$
for each $\gamma \in \hat{\cG}$ (the assumption on the injectivity
of $\sigma$ in \cite{Boca} is not used in the proof of Theorem 17 therein; this
can be also seen by passing to the associated reduced action in Proposition~\ref{reduced:prop}
for which $\sigma_r$ is always injective).

The pair $(\cB, \sigma|_{\cB})$ consists of the $*$-algebra $\cB$ and
the action $\sigma|_{\cB}:\cB\rightarrow \cB\odot \cA$.
For each $\gamma \in \hat{\cG}$,
one has a linear basis $\sS_{\gamma}$ of $B_{\gamma}$ satisfying (\ref{action:eq}),
where we take $J_{\gamma}$ to be $\{1, {\cdots}, \mul(B, \gamma)\}$.
If we choose such a basis $\sS_{\gamma}$ for each $\gamma \in \hat{\cG}$, then the action
$\sigma|_{\cB}$
is fixed by (\ref{action:eq}) and the pair
$(\cB, \sigma|_{\cB})$ is determined by the $*$-algebra structure on $\cB$ which in turn can
be determined by the coefficients appearing in the multiplication
and $*$-operation rules on these basis elements.  In order to
reduce the set of possible coefficients appearing this way, we
put one more restriction on $\sS_{\gamma}$.
By the argument on \cite[page 103]{Boca}, one can
require $\sS_{\gamma}$ to be an orthonormal basis of
$B_{\gamma}$ with respect to the inner product $\left<x, y\right>=\omega(x^*y)$,
that is,
\begin{eqnarray} \label{ON:eq}
\omega(e^*_{\gamma sj}e_{\gamma ki})=\delta_{sk}\delta_{ji}.
\end{eqnarray}
We can always choose $e_{\gamma_011}=1_{B}$.
We shall call a basis $\sS_{\gamma}$ satisfying all these conditions
a {\it standard
basis} of $B_{\gamma}$, and call
the union $\sS$ of a standard basis for each $B_{\gamma}$ a standard basis of $B$.

\begin{notation} \label{sets:notation}
Set
\begin{eqnarray*}
\hat{\cG}^{\flat}=\hat{\cG}\setminus \{\gamma_0\},\quad \sM=\{(\alpha, \beta, \gamma)\in
\hat{\cG}^{\flat}\times \hat{\cG}^{\flat}\times \hat{\cG}:\gamma\preceq \alpha\tot \beta\}.
\end{eqnarray*}
For each $\gamma \in \hat{\cG}$, set
\begin{eqnarray*}
X_{\gamma}&=&\{(\gamma, k, i): 1\le k\le N_{\gamma},  1\le i\le d_{\gamma}\}, \\
X'_{\gamma}&=&\{(\gamma, k, i): 1\le k\le \mul(B, \gamma),  1\le i\le d_{\gamma}\}.
\end{eqnarray*}
Denote by $x_0$ the unique element $(\gamma_0, 1, 1)$
in $X_{\gamma_0}$. Set
\begin{eqnarray*}
Y=\cup_{(\alpha, \beta, \gamma)\in \sM}X_{\alpha}\times X_{\beta}\times X_{\gamma}, & &\quad
Z= \cup_{\gamma\in \hat{\cG}^{\flat}} X_{\gamma}\times X_{\gamma^c}, \\
Y'=\cup_{(\alpha, \beta, \gamma)\in \sM}X'_{\alpha}\times X'_{\beta}\times X'_{\gamma}, & &\quad
Z'= \cup_{\gamma\in \hat{\cG}^{\flat}} X'_{\gamma}\times X'_{\gamma^c}.
\end{eqnarray*}
\end{notation}

Fix a standard basis $\sS$ of $\cB$.
Since we have chosen $e_{x_0}$ to be $1_{B}$, the algebra
structure of $\cB$ is determined by
the linear expansion of $e_{x_1}e_{x_2}$ for all
$x_1\in X'_{\alpha}, x_2\in X'_{\beta}, \alpha, \beta \in \hat{\cG}^{\flat}$.
By (\ref{*-algebra:eq}) we have
$e_{x_1}e_{x_2}\in \sum_{\gamma \preceq \alpha\tot \beta}B_{\gamma}$.
Thus the coefficients of the expansion of $e_{x_1}e_{x_2}$
under $\sS$ for all such $x_1, x_2$ determine a scalar function on $Y'$,
that is,
there exists a unique element $f\in \Ce^{Y'}$ such that for any $x_1\in X'_{\alpha}, x_2\in X'_{\beta}, \alpha, \beta\in
\hat{\cG}^{\flat}$,
\begin{eqnarray} \label{coefficient1:eq}
e_{x_1}e_{x_2}=\sum_{(x_1, x_2, x_3)\in Y'}f(x_1, x_2, x_3)e_{x_3}.
\end{eqnarray}
Similarly, the $*$-structure of $\cB$ is determined by
the linear expansion of $e_{x_1}$ for all
$x_1\in X'_{\gamma}, \gamma \in \hat{\cG}^{\flat}$.
By (\ref{*-algebra:eq}) we have
$e^*_{x_1}\in B_{\gamma^c}$. Thus there exists a unique element $g\in \Ce^{Z'}$ such that for any $x_1\in X'_{\gamma},
\gamma \in \hat{\cG}^{\flat}$,
\begin{eqnarray} \label{coefficient3:eq}
e^*_{x_1}=\sum_{(x_1, x_2)\in Z'}g(x_1, x_2)e_{x_2}.
\end{eqnarray}
Then $(f, g)$ determines the isomorphism class of $(\cB,
\sigma|_{\cB})$ and hence determines the equivalence class of
$(B, \sigma)$ in $\EA^{\sim}(\cG)$. Note that $(f, g)$ does not
determine the isomorphism class of $(B, \sigma)$ in $\EA(\cG)$ unless
$(B, \sigma)$ is co-amenable. Since we are going to consider all
ergodic actions of $\cG$ in a uniform way, we extend $f$ and $g$
to functions on $Y$ and $Z$ respectively by
\begin{eqnarray} \label{coefficient2:eq}
f|_{Y\setminus Y'}=0, \quad g|_{Z\setminus Z'}=0.
\end{eqnarray}
We shall say that $(f, g)$ is the element in $\Ce^Y\times \Ce^Z$ associated to
$\sS$.

Denote by $\cP$ the set of $(f, g)$ in $\Ce^Y\times \Ce^Z$ associated to
various bases of ergodic actions of $\cG$. We say that $(f_1, g_1)$ and
$(f_2, g_2)$ in $\cP$ are {\it equivalent}
if they are associated to standard bases of $(B_1, \sigma_1)$ and
$(B_2, \sigma_2)$ respectively such that $(\cB_1, \sigma_1|_{\cB_1})$
and $(\cB_2, \sigma_2|_{\cB_2})$ are isomorphic.
Then this is an equivalence relation on $\cP$ and
we can identify the quotient space of $\cP$ modulo this equivalence relation
with $\EA^{\sim}(\cG)$ naturally.

\begin{definition} \label{topology on EA:def}
Endow $\Ce^Y\times \Ce^Z$ with the product topology.
Define the topology on  $\cP$ as the relative topology, and
define the topology on  $\EA^{\sim}(\cG)$ as the quotient topology from
$\cP\rightarrow \EA^{\sim}(\cG)$.
Also define the topology on $\EA(\cG)$ via setting the
open subsets in $\EA(\cG)$ as inverse image of open subsets in
$\EA^{\sim}(\cG)$ under the quotient map $\EA(\cG)\rightarrow \EA^{\sim}(\cG)$.
\end{definition}

\begin{theorem} \label{topology on EA:thm}
Both $\cP$
and $\EA^{\sim}(\cG)$ are compact Hausdorff spaces.
The space $\EA(\cG)$ is also compact, but it
is Hausdorff if and only if $\cG$ is co-amenable.
Both quotient maps $\cP\rightarrow \EA^{\sim}(\cG)$ and
$\EA(\cG)\rightarrow \EA^{\sim}(\cG)$ are open.
\end{theorem}

\begin{remark} \label{indep of ON:remark}
The equation (\ref{action:eq}) depends on the identification
of $B(H_{\gamma})$ with $M_{d_{\gamma}}(\Ce)$, which in turn depends
on the choice of an orthonormal basis of $H_{\gamma}$. Then $\cP$ also depends
on such choice. However, using Lemma~\ref{change of basis:lemma} below one can
show directly that the quotient topology on $\EA^{\sim}(\cG)$ does not
depend on such choice. This will also follow from Corollary~\ref{unique topology:coro} below.
\end{remark}

In order to prove Theorem~\ref{topology on EA:thm}, we need
to characterize $\cP$ and its equivalence relation more explicitly.
We start with characterizing $\cP$, that is,
we consider which elements of $\Ce^Y\times \Ce^Z$  come from standard bases of ergodic
actions of $\cG$. For this purpose, we take $f(y)$ and $g(z)$ for $y\in Y, z\in Z$ as
variables and try to find algebraic conditions they should satisfy in order to construct
$(\cB, \sigma|_{\cB})$.
Set
\begin{eqnarray*}
X=\cup_{\gamma \in \hat{\cG}^{\flat}}X_{\gamma}, \quad \mbox{ and }\quad X_0=\cup_{\gamma \in \hat{\cG}}X_{\gamma}.
\end{eqnarray*}
Let $\sV$ be a vector space with basis $\{v_x:x\in X_0\}$.
We hope to construct $\cB$ out of $\sV$ such that $v_x$ becomes $e_x$.
Corresponding to (\ref{coefficient1:eq})-(\ref{coefficient2:eq}) we want to make $\sV$ into
a $*$-algebra with identity $v_{x_0}$ satisfying
\begin{eqnarray} \label{coefficient11:eq}
v_{x_1}v_{x_2}=\sum_{(x_1, x_2, x_3)\in Y}f(x_1, x_2, x_3)v_{x_3}
\end{eqnarray}
for any $x_1, x_2\in  X$, and
\begin{eqnarray} \label{coefficient31:eq}
v^*_{x_1}=\sum_{(x_1, x_2)\in Z}g(x_1, x_2)v_{x_2},
\end{eqnarray}
for any $x_1\in X$.
Corresponding to (\ref{action:eq}), we also want a unital $*$-homomorphism $\sigma_\sV: \sV\rightarrow \sV\odot \cA$ satisfying
\begin{eqnarray} \label{coefficient51:eq}
\sigma_{\sV}(v_{\gamma ki})=\sum_{1\le j\le d_{\gamma}}v_{\gamma kj}\otimes u^{\gamma}_{ji}
\end{eqnarray}
for $(\gamma, k, i)\in X$.
Thus consider the equations
\begin{eqnarray*}
(v_{x_1}v_{x_2})v_{x_3}=v_{x_1}(v_{x_2}v_{x_3}), & &\quad (v^*_{x_1})^*=v_{x_1}, \quad (v_{x_1}v_{x_2})^*=v^*_{x_2}v^*_{x_1}, \\
\sigma_\sV(v_{x_1}v_{x_2})=\sigma_\sV(v_{x_1})\sigma_\sV(v_{x_2}), & &\quad (\sigma_\sV(v_{x_1}))^*=\sigma_\sV(v^*_{x_1})
\end{eqnarray*}
for all $x_1, x_2, x_2\in X$. Expanding both sides of these equations
formally using (\ref{coefficient11:eq})-(\ref{coefficient51:eq})
and identifying  the corresponding coefficients,
we get a set $\sE_1$ of polynomial equations in the variables $f(y), g(z)$ and their conjugates for $y\in Y, z\in Z$.
For any $(f, g)\in \Ce^Y\times \Ce^Z$ satisfying $\sE_1$,
we have a conjugate-linear map $*:\sV\rightarrow \sV$ specified by (\ref{coefficient31:eq}).
Set $I_{f, g}$ to be the kernel of $*$, and set $\sV_{f, g}=\sV/I_{f, g}$.
Denote the quotient map $\sV\rightarrow \sV_{f, g}$
by $\phi_{f, g}$, and denote $\phi_{f, g}(v_x)$ by $\nu_x$
for $x\in X_0$.
Then the formulas
\begin{eqnarray}
\label{coefficient12:eq} \nu_{x_1}\nu_{x_2}&=&\sum_{(x_1, x_2, x_3)\in Y}f(x_1, x_2, x_3)\nu_{x_3}, \\
\label{coefficient32:eq} \nu_{x_1}^*&=&\sum_{(x_1, x_2)\in Z}g(x_1, x_2)\nu_{x_2},\\
\label{coefficient52:eq} \sigma_{f, g}(\nu_{\gamma ki})&=& \sum_{1\le j\le d_{\gamma}}\nu_{\gamma kj}\otimes u^{\gamma}_{ji}
\end{eqnarray}
corresponding to (\ref{coefficient11:eq})-(\ref{coefficient51:eq})
determine a unital $*$-algebra structure of $\sV_{f, g}$ with the identity $\nu_{x_0}$ and
a unital $*$-homomorphism $\sigma_{f, g}:\sV_{f, g}\rightarrow \sV_{f, g}\odot \cA$.

In order to make sure that $(f, g)$ is associated to some standard basis of
some ergodic action of $\cG$, we need to
also take care of (\ref{ON:eq}). Note that $\omega|_{\cB}$ is simply to take the coefficient at $1_{B}$.
For any $\gamma\in \hat{\cG}^{\flat}$ and any
$x_1, x_2\in X_{\gamma}$, expand $v^*_{x_2}v_{x_1}$ formally using (\ref{coefficient31:eq}) and (\ref{coefficient11:eq})
and denote by $\F_{x_1, x_2}$
the coefficient at $1_\sV$.
Then we want the existence of a non-negative  integer $m_{\gamma, f, g}\le N_{\gamma}$ for
each $\gamma\in \hat{\cG}^{\flat}$, which one expects to be
$\mul(B, \gamma)$,
such that the value of $\F_{\gamma sj, \gamma ki}$
at $(f, g)$ is $\delta_{sk}\delta_{ji}$ or $0$ depending on
$s, k\le m_{\gamma, f, g}$ or not. This condition can be expressed
as the set $\sE_2$ of the equations
$\F_{x_1, x_2}=0$ for all $x_1, x_2\in X_{\gamma}$ with
$x_1\neq x_2$,
the equations $\F_{\gamma si, \gamma si}=\F_{\gamma sj, \gamma sj}$ for all
$1\le i, j\le d_{\gamma}, 1\le s\le N_{\gamma}$,
and the equations $\F_{\gamma s1, \gamma s1}\F_{\gamma k1, \gamma k1}=\F_{\gamma s1, \gamma s1}$
for all $1\le k\le s\le N_{\gamma}$ (and for all $\gamma \in \hat{\cG}^{\flat}$).
We also need to take care of (\ref{coefficient2:eq}). Thus denote by $\sE_3$
the set of equations
\begin{eqnarray*}
f(x_1, x_2, x_3)=f(x_1, x_2, x_3)\F_{x_1, x_1}=f(x_1, x_2, x_3)\F_{x_2, x_2}=f(x_1, x_2, x_3)\F_{x_3, x_3}
\end{eqnarray*}
for all $(x_1, x_2, x_3)\in Y$ (the last equation is vacuous when $x_3=x_0$), and
the equations
\begin{eqnarray*}
g(x_1, x_2)=g(x_1, x_2)\F_{x_1, x_1}=g(x_1, x_2)\F_{x_2, x_2}
\end{eqnarray*}
for all $(x_1, x_2)\in Z$.

\begin{notation} \label{equation:notation}
Denote by $\sE$ the union of $\sE_1, \sE_2$ and $\sE_3$.
\end{notation}

Clearly every element in $\cP$ satisfies $\sE$. This proves
part of the following characterization of $\cP$:

\begin{proposition} \label{P:prop}
$\cP$ is exactly the set of elements in $\Ce^Y\times \Ce^Z$ satisfying $\sE$.
\end{proposition}

Let $(f, g)\in \Ce^Y\times \Ce^Z$ satisfy $\sE$.  Set
\begin{eqnarray*}
X_{f, g}=\{(\gamma, s, i)\in X: 1\le s\le m_{\gamma, f, g}\},
\end{eqnarray*}
which one expects to parameterize $\sS\setminus \{e_{x_0}\}$.
Since $(f, g)$ satisfies $\sE_3$, $\spn \{v_x:x\in X\setminus
X_{f, g}\}\subseteq I_{f, g}$. Thus $\nu_x$'s for $x\in X_{f,
g}\cup \{x_0\}$ span $\sV_{f, g}$. Clearly $\sV_{f, g}$ is the
direct sum of $\Ce \nu_{x_0}$ and $\spn \{\nu_x:x\in X_{\gamma}\}$
for all $\gamma \in \hat{\cG}^{\flat}$. Thus it makes sense to
talk about the coefficient of $\nu$ at $\nu_{x_0}$ for any $\nu\in
\sV_{f, g}$. This defines a linear functional $\varphi_{f, g}$ on
$\sV_{f, g}$, which one expects to be $\omega$. Clearly
\begin{eqnarray} \label{invariant state:eq}
\varphi_{f, g}(\cdot)\nu_{x_0}=(\id \otimes h)\sigma_{f, g}(\cdot)
\end{eqnarray}
on $\sV_{f, g}$.

\begin{lemma} \label{exists universal:lemma}
Let $(f, g)\in \Ce^Y\times \Ce^Z$ satisfy $\sE$.
Then $\sV_{f, g}$ has a universal $C^*$-algebra $B_{f, g}$.
The canonical $*$-homomorphism $\sV_{f, g}\rightarrow B_{f, g}$ is injective.
Identifying $\sV_{f, g}$ with its canonical image
in $B_{f, g}$ one has
\begin{eqnarray} \label{bound2:eq}
\pa \nu_x\pa \le \sqrt{\pa F_{\gamma}\pa M_{\gamma}}
\end{eqnarray}
for any $x=(\gamma, k, i)\in X_{f, g}$, where
$F_{\gamma}$ denotes the element in $M_{d_{\gamma}}(\Ce)$ defined after Theorem 5.4 in \cite{Wor87}.
The set $\sS:=\{\nu_x:x\in X_{f, g}\cup\{x_0\}\}$
is a linear basis of $\sV_{f, g}$.
\end{lemma}
\begin{proof}
We show first that for each $\nu \in \sV_{f, g}$ there  exists some $c_{\nu}\in \Re$
such that $\pa \pi(\nu)\pa \le c_{\nu}$ for any $*$-representation $\pi$ of $\sV_{f, g}$.
Recalling that $\sS$ spans $\sV_{f, g}$,
it suffices to prove the claim for $\nu=\nu_x$ for every $x\in X_{f, g}$.
Say $x=(\gamma, k, i)$. Set
\begin{eqnarray*}
W_{\gamma k}&=&F^{-\frac{1}{2}}_{\gamma}(\nu_{\gamma k1}, {\cdots}, \nu_{\gamma  kd_{\gamma}})^T\in M_{d_{\gamma}\times 1}
(\sV_{f, g}),\\
W'_{\gamma k}&=&(W_{\gamma k}, 0, {\cdots}, 0)\in M_{d_{\gamma}\times d_{\gamma}}
(\sV_{f, g}).
\end{eqnarray*}
Note that $\{\nu \in \sV_{f, g}:\sigma_{f, g}(\nu)=\nu\otimes 1_{A}\}=\Ce\nu_{x_0}$.
The argument in \cite[page 103]{Boca} shows that
\begin{eqnarray} \label{bound3:eq}
W^*_{\gamma k}W_{\gamma s}=\delta_{ks}M_{\gamma}\nu_{x_0}
\end{eqnarray}
for all $1\le k, s\le m_{\gamma, f, g}$. Thus for any $*$-representation $\pi$ of $\sV_{f, g}$
we have
$\pa \pi(W'_{\gamma k})\pa \le \sqrt{M_{\gamma}}$ and hence
\begin{eqnarray} \label{bound:eq}
\pa \pi(\nu_x)\pa \le
\pa F^{\frac{1}{2}}_{\gamma}\pa \sqrt{M_{\gamma}}=\sqrt{\pa F_{\gamma}\pa M_{\gamma}}.
\end{eqnarray}

Next we show that $\sV_{f, g}$ does have a $*$-representation.
By \cite[Theorem 5.7]{Wor87} one has $h(A^*_{\alpha}A_{\beta})=0$ for any
$\alpha\neq \beta \in \hat{\cG}$. Using (\ref{invariant state:eq}) one sees that
\begin{eqnarray} \label{orthogonal:eq}
\varphi_{f, g}(\nu_{x_2}^*\nu_{x_1})=0
\end{eqnarray}
for all $x_1\in X_{\alpha}, x_2\in X_{\beta},\alpha \neq \beta$.
Using (\ref{orthogonal:eq}) and the assumption that $(f, g)$ satisfies $\sE_2$, one observes
that
\begin{eqnarray} \label{positive:eq}
\varphi_{f, g}(\nu^*\nu)=\sum_{x\in X_{f, g}\cup \{x_0\}}|\lambda_x|^2\ge 0
\end{eqnarray}
for any $\nu=\sum_{x\in X_0}\lambda_x\nu_x \in \sV_{f, g}$. Denote
by $H$ the Hilbert space completion of $\sV_{f, g}$ with respect
to the inner product $\left<\nu_1, \nu_2\right>=\varphi_{f,
g}(\nu^*_1\nu_2)$, and by $H^{(d_{\gamma})}$ the direct sum of
$d_{\gamma}$ copies of $H$. By (\ref{bound3:eq}) the
multiplication by $W'_{\gamma s}$ extends to a bounded operator on
$H^{(d_{\gamma})}$. Then so does the multiplication by
$((\nu_{\gamma k1}, {\cdots}, \nu_{\gamma  kd_{\gamma}})^T, 0,
{\cdots}, 0)\in M_{d_{\gamma}\times d_{\gamma}}(\sV_{f, g})$.
Consequently, the multiplication by $\nu_x$ for $x=(\gamma, k,
i)\in X_{f, g}$ extends to a bounded operator on $H$. Since
$\sS$ spans $\sV_{f, g}$, the multiplication of $\sV_{f, g}$
extends to a $*$-representation $\pi$ of $\sV_{f, g}$ on $H$.

Now we conclude that $\sV_{f, g}$ has a universal $C^*$-algebra
$B_{f, g}$. It follows from (\ref{positive:eq}) that $\pi\circ
\phi_{f, g}$ is injective on $\spn\{v_x: x\in X_{f, g}\cup
\{x_0\}\}$, where $\phi_{f, g}:\sV\rightarrow \sV_{f, g}$ is the
quotient map. Thus $\sS$ is a linear basis of $\sV_{f, g}$, and
the canonical $*$-homomorphism $\sV_{f, g}\rightarrow B_{f,g}$
must be injective. The inequality (\ref{bound2:eq}) follows from
(\ref{bound:eq}).
\end{proof}

For $(f, g)$ as in Lemma~\ref{exists universal:lemma},
by the universality of $B_{f, g}$, the $*$-homomorphism
$\sV_{f, g}\overset{\sigma_{f, g}}\rightarrow \sV_{f, g}\odot \cA\hookrightarrow
B_{f, g}\otimes A$ extends uniquely to a (unital) $*$-homomorphism $B_{f, g}\rightarrow B_{f, g}\otimes A$,
which we still denote by $\sigma_{f, g}$.

\begin{proposition} \label{coefficient to ergodic:prop}
Let $(f, g)$ be as in Lemma~\ref{exists universal:lemma}. Then $\sigma_{f, g}$ is an ergodic action of
$\cG$ on $B_{f, g}$. The algebra of regular functions for
this action is $\sV_{f, g}$. The set $\sS$
is a standard basis
of $\sV_{f, g}$. The element in $\cP$ associated to this basis is exactly $(f, g)$.
\end{proposition}
\begin{proof}
By Lemma~\ref{exists universal:lemma}, $\sS$ is a basis of
$\sV_{f, g}$. By (\ref{coefficient52:eq}) and the characterization
of actions of $\cG$ in terms of elements satisfying
(\ref{action:eq}) in subsection~\ref{actions:sub}, $\sigma_{f, g}$
is an ergodic action of $\cG$ on $B_{f, g}$, and $(B_{f,
g})_{\gamma}=\spn\{\nu_x:x=(\gamma, k, i)\in X_{f, g}\}$,
$\mul(B_{f, g}, \gamma)=m_{\gamma, f, g}$ for all $\gamma \in
\hat{\cG}^{\flat}$. Thus $\sV_{f,g}$ is  the algebra of regular
functions. Denote by $\omega$ the unique $\cG$-invariant state on
$B_{f, g}$. By (\ref{invariant state:eq}) $\omega$ extends
$\varphi_{f, g}$. Since $(f, g)$ satisfies $\sE_2$, we have
$\omega(\nu^*_x\nu_y)=\delta_{xy}$ for any $x=(\gamma, k, i),
y=(\gamma, s, j)\in X_{f, g}$. Thus $\sS$ is a standard basis of
$\cB_{f, g}$. Clearly the element in $\cP$ associated to this
basis is exactly $(f, g)$.
\end{proof}

Now Proposition~\ref{P:prop} follows from Proposition~\ref{coefficient to ergodic:prop}.

We are ready to prove the compactness of $\cP$.

\begin{lemma} \label{compact:lemma}
Let $(f, g)\in \cP$. Then
\begin{eqnarray} \label{compact:eq1}
|f(x_1, x_2, x_3)|\le \sqrt{\pa F_{\alpha}\pa M_{\alpha}}
\end{eqnarray}
for any $(x_1, x_2, x_3)\in Y, x_1\in X_{\alpha}$.
And
\begin{eqnarray} \label{compact:eq2}
|g(x_1, x_2)|\le  \sqrt{\pa F_{\alpha}\pa M_{\alpha}}
\end{eqnarray}
for any $(x_1, x_2)\in Z, x_1\in X_{\alpha}$.
The space $\cP$ is compact.
\end{lemma}
\begin{proof}
Say, $(f, g)$ is associated to a standard basis $\sS$ for an
ergodic action $(B, \sigma)$ of $\cG$. Let $(H_{B}, \pi_r)$
be the $\GNS$ representation associated to the unique
$\sigma$-invariant state $\omega$ of $B$. Then $B_{\alpha}$
and $B_{\beta}$ are orthogonal to each other in $H_{B}$ for
distinct $\alpha, \beta \in \hat{\cG}$ \cite[Corollary 12]{Boca}.
In view of (\ref{ON:eq}), $\sS$ is an orthonormal basis of
$H_{B}$. We may identify $\cB$ with $\sV_{f, g}$ naturally via
$e_x\leftrightarrow \nu_x$. Then there is a $*$-homomorphism from
$B_{f, g}$ in Lemma~\ref{exists universal:lemma} to $B$
extending this identification. Thus by (\ref{bound2:eq}) we have
$\pa e_x\pa \le \sqrt{\pa F_{\alpha}\pa M_{\alpha}}$ for any $x\in
X'_{\alpha}, \alpha\in \hat{\cG}$. For any $(x_1, x_2, x_3)\in Y',
x_1\in X'_{\alpha}$, by (\ref{coefficient1:eq}),
\begin{eqnarray*}
|f(x_1, x_2, x_3)|=|\left<e_{x_3}, e_{x_1}e_{x_2}\right>|\le \pa
e_{x_1}\pa \le \sqrt{\pa F_{\alpha}\pa M_{\alpha}}.
\end{eqnarray*}
If $y\in Y\setminus Y'$, then $f(y)=0$ by (\ref{coefficient2:eq}).
This proves (\ref{compact:eq1}). The inequality (\ref{compact:eq2})
is proved similarly.

By Proposition~\ref{P:prop} the space $\cP$ is closed
in $\Ce^Y\times \Ce^Z$. It follows
from (\ref{compact:eq1}) and (\ref{compact:eq2}) that $\cP$ is compact.
\end{proof}

Next we characterize the equivalence relation on $\cP$. For this
purpose, we need to consider the relation between two standard
bases of $\cB$. The argument in the proof of \cite[Theorem
1.5]{Podles95} shows the first two assertions of the following
lemma:

\begin{lemma} \label{change of basis:lemma}
Let $\gamma\in \hat{\cG}$.
If $b_i\in B, 1\le i\le d_{\gamma}$ satisfy
\begin{eqnarray} \label{action:eq2}
\sigma(b_i)=\sum_{1\le j\le d_{\gamma}}b_j\otimes u^{\gamma}_{ji}
\end{eqnarray}
for all $1\le i\le d_{\gamma}$, then $b_i=E^{\gamma}_{i1}(b_1)$
(see (\ref{E_ij:eq})) for all $1\le i\le d_{\gamma}$. Conversely,
given $b\in E^{\gamma}_{11}(B)$, if we set
$b_i=E^{\gamma}_{i1}(b)$, then $b_i\in B_{\gamma}, 1\le i\le
d_{\gamma}$ satisfy (\ref{action:eq2}), and $b_1=b$. For any
$b_1, {\cdots}, b_{d_{\gamma}}$ ($b'_1, {\cdots}, b'_{d_{\gamma}}$
resp.)  in $B$ satisfying (\ref{action:eq2}) ((\ref{action:eq2})
with $b_i$ replaced by $b'_i$ resp.) we have
\begin{eqnarray} \label{ON:eq1}
\omega(b^*_jb'_i)=\delta_{ji}\omega(b^*_1b'_1)
\end{eqnarray}
for all $1\le i, j\le d_{\gamma}$.
\end{lemma}
\begin{proof}
We just need to prove (\ref{ON:eq1}). By \cite[Theorem 5.7]{Wor87} we have
\begin{eqnarray} \label{or relation:eq}
h({u^{\gamma}_{lj}}^*u^{\gamma}_{ni})=\frac{1}{M_{\gamma}}f_{-1}(u^{\gamma}_{nl})\delta_{ji},
\end{eqnarray}
where $f_{-1}$ is the linear functional on $\cA$ defined in \cite[Theorem 5.6]{Wor87}.
Thus
\begin{eqnarray*}
\omega(b^*_jb'_i)1_{B}&=&(\id \otimes h)(\sigma(b^*_jb'_i))
\overset{(\ref{action:eq2})}= (\id \otimes h)(\sum_{1\le l, n\le d_{\gamma}}b^*_{l}b'_{n}\otimes
{u^{\gamma}_{lj}}^*u^{\gamma}_{ni})\\
&\overset{(\ref{or relation:eq})}=&
 \sum_{1\le l, n\le d_{\gamma}}b^*_{l}b'_{n}\frac{1}{M_{\gamma}}f_{-1}(u^{\gamma}_{nl})\delta_{ji}.
\end{eqnarray*}
Therefore
\begin{eqnarray*}
\omega(b^*_jb'_i)1_{B}=\delta_{ji}\omega(b^*_1b'_1)1_{B},
\end{eqnarray*}
which proves (\ref{ON:eq1}).
\end{proof}

By Lemma~\ref{change of basis:lemma}, for $\gamma\in \hat{\cG}^{\flat}$, there is a $1$-$1$ correspondence between
standard bases of $B_{\gamma}$ and
orthonormal bases
of $E^{\gamma}_{11}(B)$ with respect to the inner product $\left<b, b'\right>=\omega(b^*b')$.
It also follows from Lemma~\ref{change of basis:lemma} that $\dim (E^{\gamma}_{11}(B))=\mul(B, \gamma)$.
Denote by $\cU_n$ the unitary group of $M_n(\Ce)$. Then $\prod_{\gamma \in \hat{\cG}^{\flat}}\cU_{\mul(B, \gamma)}$ has a
right free transitive action
on the set of standard bases of $\cB$ via acting on
the set of orthonormal bases of $E^{\gamma}_{11}(B)$
for each $\gamma \in \hat{\cG}^{\flat}$.
For $n\le m$ identify $\cU_n$ with the subgroup of $\cU_m$ consisting
of elements with $1_{m-n}$ at the lower-right corner.
Denote $\prod_{\gamma \in \hat{\cG}^{\flat}}\cU_{N_{\gamma}}$ by $\cU$, equipped with the product topology.
Then $\cU$ has a natural partial right (not necessarily free) action $\tau$ on $\cP$, that is, $\xi \in \cU$ acts at $t\in \cP$ exactly
if $\xi\in \prod_{\gamma \in \hat{\cG}^{\flat}}\cU_{m_{\gamma, t}}$, where $m_{\gamma, t}$ was defined in the paragraph before
Notation~\ref{equation:notation},
and the image $t{\cdot}\xi$ is the element in
$\cP$ associated to the standard basis $\sS{\cdot}\xi$ of $\cB_t$, where $\sS{\cdot}\xi$ is the image of the action of
$\xi$ at the standard basis $\sS$ of $\cB_t$ in Proposition~\ref{coefficient to ergodic:prop}.
Clearly the orbits of this partial action are exactly the fibres of the quotient map
$\cP\rightarrow \EA^{\sim}(\cG)$, equivalently, exactly the equivalence classes in $\cP$ introduced
before Definition~\ref{topology on EA:def}.
Thus we may identify $\EA^{\sim}(\cG)$ with the quotient space $\cP/\cU$.

\begin{lemma} \label{open:lemma}
The quotient map $\cP\rightarrow \cP/\cU$ is open. The quotient
topology on $\cP/\cU$ is compact Hausdorff.
\end{lemma}
\begin{proof} Denote by $\pi$ the quotient map $\cP\rightarrow \cP/\cU$.
To show the openness of $\pi$, it suffices to show that $\pi^{-1}(\pi(V))$
is open for every open subset $V$ of $\cP$. Let $t\in V$ and $\xi\in \cU$ such
that $t{\cdot}\xi$ is defined. Say $\xi=(\xi_{\gamma})_{\gamma\in \hat{\cG}^{\flat}}$.
Let $J$ be a finite subset of $\hat{\cG}$. Replacing $\xi_{\gamma}$ by $1_{N_{\gamma}}$
for $\gamma \in \hat{\cG}^{\flat}\setminus J$ we get an element $\xi'\in \cU$. Notice that
when $t'\in \cP$ is close enough to $t{\cdot}\xi$, $t'{\cdot}(\xi')^{-1}$ is defined.
Moreover, the restrictions of $ t'{\cdot}(\xi')^{-1}$ on $(X_{J}\times X_{J}\times X_{J})\cap Y$
and $(X_{J}\times X_{J})\cap Z$ converge to the restrictions of $t$ as $t'$ converges to $t{\cdot}\xi$,
where $X_{J}=\cup_{\gamma\in J}X_{\gamma}$.
Clearly we can find a large enough finite subset $J$ of $\hat{\cG}$ such that when $t'$ is close enough
to $t\cdot \xi$, the element $t'\cdot (\xi')^{-1}$ is in $V$.
Then $t'=(t'{\cdot}(\xi')^{-1}){\cdot} \xi'$ is in $\pi^{-1}(\pi(V))$. Therefore $\pi^{-1}(\pi(V))$ is open,
and hence $\pi$ is open.

Denote by $\De$ the domain of $\tau$, i.e., the subset of
$\cP\times \cU$ consisting of elements  $(t, \xi)$ for which
$t{\cdot}\xi$ is defined. From the equations in $\sE_2$ it is
clear that $\mathbb D$ is closed in
$\cP\times \cU$. By Lemma~\ref{compact:lemma} the space $\cP$ is
compact. Since $\cU$ is also compact, so is $\De$. It is also
clear that $\tau$ is continuous in the sense that the map
$\De\rightarrow \cP$ sending $(t, \xi)$ to $t{\cdot}\xi$ is
continuous. Thus the set $\{(t, t')\in \cP\times
\cP:\pi(t)=\pi(t')\}= \{(t, t{\cdot}\xi)\in \cP\times \cP: (t,
\xi)\in \De\}$ is closed in $\cP\times \cP$. Since $\pi$ is open,
a standard argument shows that the quotient topology on $\cP/\cU$
is compact Hausdorff.
\end{proof}

Since $\EA^{\sim}(\cG)$ is Hausdorff by Lemma~\ref{open:lemma}, $\EA(\cG)$
is Hausdorff exactly if the quotient map $\EA(\cG)\rightarrow
\EA^{\sim}(\cG)$ is a bijection, exactly if $\cG$ is co-amenable by
Proposition~\ref{universal:prop}. Then Theorem~\ref{topology on
EA:thm} follows from Lemmas~\ref{compact:lemma} and
\ref{open:lemma}.

Notice that the function $t\mapsto m_{\gamma, t}$ is continuous on $\cP$ for each $\gamma \in \hat{\cG}^{\flat}$.
Thus we have

\begin{proposition} \label{mul:prop}
The multiplicity function
$\mul(\cdot, \gamma)$ is continuous on both $\EA(\cG)$ and
$\EA^{\sim}(\cG)$ for each $\gamma \in \hat{\cG}$.
\end{proposition}


To end this section, we discuss the behavior of $\EA(\cG)$ when we take
Cartesian products of compact quantum groups. Let
$\{A_{\lambda}=C(\cG_{\lambda})\}_{\lambda \in \Lambda}$ be a family of compact
quantum groups indexed by a set $\Lambda$. Then
$\otimes_{\lambda}A_{\lambda}$ has a unique compact quantum
group structure such that the embeddings
$A_{\mu}\hookrightarrow \otimes_{\lambda}A_{\lambda}$ for 
$\mu\in \Lambda$ are
all morphisms between compact quantum groups \cite[Theorem 1.4,
Proposition 2.6]{Wang95}, which we shall denote by
$C(\prod_{\lambda}\cG_{\lambda})$. The Haar measure of
$\otimes_{\lambda}A_{\lambda}$ is the tensor product
$\otimes_{\lambda}h_{\lambda}$ of the Haar measures $h_{\lambda}$
of $A_{\lambda}$ \cite[Proposition 2.7]{Wang95}.

If $\sigma_{\lambda}:B_{\lambda}\rightarrow B_{\lambda}\otimes
A_{\lambda}$ is an action of $\cG_{\lambda}$ on a unital $C^*$-algebra $B_{\lambda}$
for each $\lambda$,
then the unique $*$-homomorphism
$\otimes_{\lambda}\sigma_{\lambda}:
\otimes_{\lambda}B_{\lambda}\rightarrow
(\otimes_{\lambda}B_{\lambda})\otimes
(\otimes_{\lambda}A_{\lambda})$ extending all
$\sigma_{\lambda}$'s is easily seen to be an action of
$\prod_{\lambda}\cG_{\lambda}$. Using the canonical conditional
expectation $\otimes_{\lambda}B_{\lambda}\rightarrow
(\otimes_{\lambda}B_{\lambda})^{\otimes_{\lambda}\sigma_{\lambda}}$,
one checks easily that
$(\otimes_{\lambda}B_{\lambda})^{\otimes_{\lambda}\sigma_{\lambda}}
=\otimes_{\lambda}B_{\lambda}^{\sigma_{\lambda}}$. In
particular, $\otimes_{\lambda}\sigma_{\lambda}$ is ergodic if and
only if every  $\sigma_{\lambda}$ is.

\begin{proposition} \label{product:prop}
Let $\{A_{\lambda}=C(\cG_{\lambda})\}_{\lambda\in \Lambda}$ be a family of compact
quantum groups indexed by a set $\Lambda$. The map
$\prod_{\lambda}\EA(\cG_{\lambda})\rightarrow
\EA(\prod_{\lambda}\cG_{\lambda})$ sending the isomorphism classes
of $(B_{\lambda}, \sigma_{\lambda})$'s to the isomorphism class
of $(\otimes_{\lambda}B_{\lambda},
\otimes_{\lambda}\sigma_{\lambda})$
descends to a map $\prod_{\lambda}\EA^{\sim}(\cG_{\lambda})\rightarrow \EA^{\sim}(\prod_{\lambda}\cG_{\lambda})$, that
is, there exists a (unique) map $\prod_{\lambda}\EA^{\sim}(\cG_{\lambda})\rightarrow \EA^{\sim}(\prod_{\lambda}\cG_{\lambda})$
such that the diagram
\begin{equation}\label{descend:diagram}
\xymatrix{
\prod_{\lambda}\EA(\cG_{\lambda}) \ar[d] \ar[r] & \EA(\prod_{\lambda}\cG_{\lambda}) \ar[d] \\
\prod_{\lambda}\EA^{\sim}(\cG_{\lambda})  \ar[r] & \EA^{\sim}(\prod_{\lambda}\cG_{\lambda})
}
\end{equation}
commutes. Moreover, both of these maps are injective and continuous,
where both $ \prod_{\lambda}\EA(\cG_{\lambda})$ and $ \prod_{\lambda}\EA^{\sim}(\cG_{\lambda})$
are endowed with the product topology.
\end{proposition}
\begin{proof}
Denote by $\prod^{\sim}_{\lambda}\widehat{\cG_{\lambda}}$ the
subset of $\prod_{\lambda}\widehat{\cG_{\lambda}}$ consisting of
elements whose all but finitely many components are classes of
trivial representations. For any $\gamma\in
\prod^{\sim}_{\lambda}\widehat{\cG_{\lambda}}$, say
$\gamma_{\lambda_1}, {\cdots}, \gamma_{\lambda_n}$ are the
nontrivial components of $\gamma$, the element
$u^{\gamma_{\lambda_1}}_{1(n+1)}
u^{\gamma_{\lambda_2}}_{2(n+2)}\cdots 
u^{\gamma_{\lambda_n}}_{n(2n)}$ (in the leg numbering notation) is
an irreducible unitary representation of
$\prod_{\lambda}\cG_{\lambda}$. Moreover, this map
$\prod^{\sim}_{\lambda}\widehat{\cG_{\lambda}}\rightarrow
\widehat{\prod_{\lambda}\cG_{\lambda}}$ is bijective \cite[Theorem
2.11]{Wang95}, and hence we may identify these two sets. Fixing an
orthonormal basis of $H_{\gamma_{\lambda}}$ we take the tensor
products of the bases of $H_{\gamma_{\lambda_1}}, {\cdots},
H_{\gamma_{\lambda_n}}$ as an orthonormal basis of
$H_{\gamma}$.

Let $\sS_{\lambda}$ be a standard basis of $B_{\lambda}$. Say,
it consists of a standard basis $\sS_{\alpha_{\lambda}}$ of $(B_{\lambda})_{\alpha_{\lambda}}$
for each $\alpha_{\lambda}\in \widehat{\cG_{\lambda}}$.
Denote by $\omega_{\lambda}$ the $\sigma_{\lambda}$-invariant state on $B_{\lambda}$.
Then $\otimes_{\lambda}\omega_{\lambda}$ is the $\prod_{\lambda}\sigma_{\lambda}$-invariant state
of $ \otimes_{\lambda}B_{\lambda}$.
Using the characterization of ergodic actions in terms of elements satisfying (\ref{action:eq}) in subsection~\ref{actions:sub},
one sees that the algebra of regular functions for
$\otimes_{\lambda}\sigma_{\lambda}$ is $\odot_{\lambda}\cB_{\lambda}$ and that
the tensor products of $\sS_{\gamma_{\lambda_1}}, {\cdots}, \sS_{\gamma_{\lambda_n}}$ is
a standard basis of $(\otimes_{\lambda}B_{\lambda})_{\gamma}$.
This shows the existence of the map $\prod_{\lambda}\EA^{\sim}(\cG_{\lambda})\rightarrow \EA^{\sim}(\prod_{\lambda}\cG_{\lambda})$
making (\ref{descend:diagram}) commute. Taking the union
of the above standard basis of $(\otimes_{\lambda}B_{\lambda})_{\gamma}$,
we also get
a standard basis of $ \otimes_{\lambda}B_{\lambda}$, which we shall denote by $\prod^{\sim}_{\lambda}\sS_{\lambda}$.
For any fixed $\lambda_0$, if we take all $\gamma \in \prod^{\sim}_{\lambda}\widehat{\cG_{\lambda}}$
whose components are trivial at all $\lambda\neq \lambda_0$ and take the sum of
the corresponding spectral subspaces of $\otimes_{\lambda}B_{\lambda}$, we
get $\cB_{\lambda_0}$. Taking norm closure, we get $B_{\lambda_0}$.
This proves the injectivity of the maps  $\prod_{\lambda}\EA(\cG_{\lambda})\rightarrow \EA(\prod_{\lambda}\cG_{\lambda})$
and
$\prod_{\lambda}\EA^{\sim}(\cG_{\lambda})\rightarrow \EA^{\sim}(\prod_{\lambda}\cG_{\lambda})$.

Clearly the map $\prod_{\lambda}\cP(\cG_{\lambda})\rightarrow \cP(\prod_{\lambda}\cG_{\lambda})$
sending $(t_{\lambda})_{\lambda \in \Lambda}$ to the element of $\cP(\prod_{\lambda}\cG_{\lambda})$
associated to the standard basis $\prod^{\sim}_{\lambda}\sS_{t_{\lambda}}$ is continuous,
where $\sS_{t_{\lambda}}$ is the standard basis of $B_{t_{\lambda}}$ in Proposition~\ref{coefficient to ergodic:prop}.
Note that  the diagram
\begin{equation}\label{P:diagram}
\xymatrix{
\prod_{\lambda}\cP(\cG_{\lambda}) \ar[d] \ar[r] & \cP(\prod_{\lambda}\cG_{\lambda}) \ar[d] \\
\prod_{\lambda}\EA^{\sim}(\cG_{\lambda})  \ar[r] & \EA^{\sim}(\prod_{\lambda}\cG_{\lambda})
}
\end{equation}
commutes, where the left vertical map is the product map.
By Theorem~\ref{topology on EA:thm} the map
$\cP(\cG_{\lambda})\rightarrow \EA^{\sim}(\cG_{\lambda})$ is open for each $\lambda$.
Thus the product map $\prod_{\lambda}\cP(\cG_{\lambda})\rightarrow \prod_{\lambda}\EA^{\sim}(\cG_{\lambda})$
is open.
It follows from the commutativity
of the diagram (\ref{P:diagram})
that the map $\prod_{\lambda}\EA^{\sim}(\cG_{\lambda})\rightarrow \EA^{\sim}(\prod_{\lambda}\cG_{\lambda})$ is continuous.
Then the continuity of the map $\prod_{\lambda}\EA(\cG_{\lambda})\rightarrow \EA(\prod_{\lambda}\cG_{\lambda})$
follows from the commutativity of the diagram (\ref{descend:diagram}).
\end{proof}

\section{Semi-continuous fields of ergodic actions}
\label{cont field:sec}

In this section we prove Theorems~\ref{cont field over P:thm} and
\ref{cont vs mul:thm}, from which we deduce
Theorems~\ref{main1:thm} and \ref{main2:thm}.

We start with discussion of semi-continuous fields of
$C^*$-algebras.

\begin{notation} \label{gloable sections:notation}
For a field $\{C_t\}_{t\in T}$ of $C^*$-algebras over
a locally compact Hausdorff space $T$, we
denote by $\prod_{t}C_t$ the $C^*$-algebra of
bounded cross-section (for the supremum norm),
and by $\prod^{\sim}_{t}C_t$ the $C^*$-algebra
of bounded cross-sections vanishing at infinity on $T$.
\end{notation}

Note that both $\prod_t C_t$ and $\prod_t^{\sim}C_t$ are
Banach modules over the $C^*$-algebra $C_{\infty}(T)$ of
continuous $\Ce$-valued functions on $T$ vanishing at infinity. We use
Rieffel's definition of semi-continuous fields of $C^*$-algebras
\cite[Definition 1.1]{Rieffel89}. We find that it is convenient to
extend the definition slightly.

\begin{definition} \label{cont fields:def}
Let $\{C_t\}_{t\in T}$ be a field of $C^*$-algebras over
a locally compact Hausdorff space $T$, and let $C$ be a $C^*$-subalgebra
of $\prod_t^{\sim}C_t$.
We say that $(\{C_t\}_{t\in T}, C)$ is a {\it topological}
field
of $C^*$-algebras if
\begin{enumerate}
\item the evaluation map $\pi_t$ from $C$ to $C_t$ is surjective for each $t\in T$,
\item $C$ is a $C_{\infty}(T)$-submodule of $\prod^{\sim}_tC_t$.
\end{enumerate}
We say that $(\{C_t\}_{t\in T}, C)$ is {\it upper semi-continuous}
({\it lower semi-continuous}, {\it continuous}, resp.) if furthermore
for each $c\in C$ the function $t\mapsto \pa \pi_t(c)\pa$ is
upper semi-continuous (lower semi-continuous, continuous, resp.).
In such case we say that $(\{C_t\}_{t\in T}, C)$ is
{\it semi-continuous}.
\end{definition}

\begin{remark} \label{unique field:remark}
If we have two upper semi-continuous fields of
$C^*$-algebras $(\{C_t\}_{t\in T}, C_1)$ and
$(\{C_t\}_{t\in T},C_2)$ over $T$ with the same fibres and
$C_1\subseteq C_2$, then $C_1=C_2$ \cite[Proposition 2.3]{DG}.
This is not true for lower semi-continuous fields
of $C^*$-algebras. For example, let $T$ be a compact Hausdorff space
and let $H$ be a Hilbert space. Take $C_t=B(H)$ for each $t$.
Set $C_1$ to be the set of all cross-sections $c$ such that
$t\mapsto \pi_t(c)$ is norm continuous, while set $C_2$ to be the set
of all norm-bounded cross-sections $c$ such that both
$t\mapsto \pi_t(c)$ and $t\mapsto (\pi_t(c))^*$  are continuous with respect to
 the strong operator topology in $B(H)$.
Then $C_1\subsetneqq C_2$ when $T$ is the one-point compactification of
$\Ne$ and $H$ is infinite-dimensional.
\end{remark}

\begin{definition} \label{homo:def}
By a  {\it homomorphism} $\varphi$ between two topological fields
of $C^*$-algebras $(\{C_t\}_{t\in T}, C)$ and $(\{B_t\}_{t\in T}, B)$ over a
locally compact Hausdorff space $T$ we mean a $*$-homomorphism
$\varphi_t:C_t\rightarrow B_t$ for each $t\in T$ such that the
pointwise $*$-homomorphism
$\prod_t\varphi_t:\prod_tC_t\rightarrow \prod_tB_t$ sends
$C$ into $B$.
\end{definition}

\begin{lemma} \label{generate cont field:lemma}
Let $\{C_t\}_{t\in T}$ be a field of $C^*$-algebras over
a locally compact Hausdorff space $T$, and let $\sC$ be a linear subspace
of $\prod_tC_t$.
 Then a section $c'\in \prod^{\sim}_tC_t$ is in
 $C:=\overline{C_{\infty}(T)\sC}$ if and only if for any $t_0\in T$ and
$\varepsilon>0$, there exist a neighborhood $U$ of $t_0$ and  $c\in \sC$ such that $\pa \pi_t(c-c')\pa <\varepsilon$
throughout $U$.  If furthermore $\pi_t(\sC)$ is dense in
$C_t$ for each $t$
and $\sC\sC, \sC^*\subseteq C$,
then $(\{C_t\}_{t\in T}, C)$ is a topological
field
of $C^*$-algebras over $T$, which we shall call the {\it topological
field generated by $\sC$}. If furthermore the function $t\mapsto \pa \pi_t(c)\pa$ is
upper semi-continuous (lower semi-continuous, continuous, resp.) for each $c\in C$,
then $(\{C_t\}_{t\in T}, C)$ is
upper semi-continuous (lower semi-continuous, continuous, resp.).
\end{lemma}
\begin{proof} The ``only if'' part is obvious. The ``if'' part follows
from a partition-of-unity argument. The second and the third assertions follow easily.
\end{proof}

Let $(\{C_t\}_{t\in T}, C)$ be a topological
field of $C^*$-algebras over a locally compact Hausdorff space $T$. If $\Theta$ is another locally
compact Hausdorff space and $p:\Theta\rightarrow T$ is a
continuous map, then we have the pull-back field
$\{C_{p(\theta)}\}_{\theta \in \Theta}$ of $C^*$-algebras over $\Theta$. There is a
natural $*$-homomorphism $p^*:\prod_tC_t\rightarrow
\prod_{\theta}C_{p(\theta)}$ sending $c$ to
$\{\pi_{p(\theta)}(c)\}_{\theta \in \Theta}$. We will call the topological field
generated by $p^*(C)$ in Lemma~\ref{generate cont field:lemma}
the {\it pull-back} of $(\{C_t\}_{t\in T}, C)$ under $p$.
In particular, if $\Theta$ is a closed
or open subset of $T$ and $p$ is the embedding, we get {\it the
restriction of $(\{C_t\}_{t\in T}, C)$ on $\Theta$}. Clearly the
pull-back and restriction of homomorphisms between topological
fields are also homomorphisms.

\begin{lemma} \label{cont field of unit:lemma}
 Let $(\{C_t\}_{t\in T}, C)$ be a semi-continuous
field of unital $C^*$-algebras over a locally compact Hausdorff space
$T$ such that the section
$\{f(t)1_{C_t}\}_{t\in T}$ is in $C$ for each $f\in C_{\infty}(T)$.
Then for any bounded function $g$ on $T$ vanishing at infinity,
the section $\{g(t)1_{C_t}\}_{t\in T}$ is in $C$ if and only if
$g\in C_{\infty}(T)$.
\end{lemma}
\begin{proof}
Via restricting to compact subsets of $T$, we may assume that $T$
is compact. The ``if'' part is given by assumption. To prove the
``only if'' part, it suffices to show that when the section
$\{g(t)1_{C_t}\}_{t\in T}$ is in $C$ and $g(t_0)=0$ for some
$t_0\in T$, we have $g(t)\to 0$ as $t\to t_0$. Replacing $g$
by $g^*g$, we may assume that $g$ is nonnegative. When the
field is upper semi-continuous, the function $t\mapsto \pa
g(t)1_{C_t}\pa=g(t)$ is upper semi-continuous at $t_0$ and
hence $g(t)\to 0$ as $t\to t_0$. When the field is lower
semi-continuous, the function $t\mapsto \pa(\pa g\pa
-g(t))1_{C_t}\pa=\pa g\pa -g(t)$ is lower semi-continuous
at $t_0$ and hence $g(t)\to 0$ as $t\to t_0$.
\end{proof}


\begin{lemma} \label{tensor product field:lemma}
Let $(\{C_t\}_{t\in T}, C)$ be a topological field of $C^*$-algebras over a locally compact Hausdorff space $T$. 
Let $D$ be a $C^*$-algebra. Then there is a natural injective $*$-homomorphism 
$\varphi: C\otimes D\rightarrow \prod^{\sim}_{t}(C_t\otimes D)$ determined by $\pi'_s(\varphi(c\otimes d))=\pi_s(c)\otimes d$
for all $c\in C$, $d\in D$, and $s\in T$, where $\pi_s$ and $\pi'_s$ denote
the coordinate maps $\prod_{t}C_t\rightarrow C_s$ and
$\prod_{t}(C_t\otimes D)\rightarrow C_s\otimes D$ respectively. 
Identifying $C\otimes D$ with $\varphi(C\otimes D)$, the pair
$(\{C_t\otimes D\}_{t\in T}, C\otimes D)$ is also a
topological field of $C^*$-algebras over $T$.
\end{lemma}
\begin{proof} 
For each $s\in T$ we have the $*$-homomorphism $\pi_s\otimes \id: (\prod_tC_t)\otimes D\rightarrow C_s\otimes D$.
Then we have the product $*$-homomorphism $(\prod_t C_t)\otimes D\rightarrow \prod_t(C_t\otimes D)$.
Denote by $\varphi$ the restriction of this homomorphism to $C\otimes D$.
We have $\pi'_s(\varphi(c\otimes d))=\pi_s(c)\otimes d$ for all $c\in C$, $d\in D$, and $s\in T$.
Clearly this identity also determines $\varphi$.

To show that $\varphi$ is injective, we may assume that $C_s$ is contained in the
algebra of bounded linear operators on $H_s$ for some Hilbert space $H_t$ for each $s\in T$,
and $D$ is contained in the algebra of bounded linear operators on $K$ for some Hilbert space
$K$. Denote
the Hilbert space direct sum $\oplus_t H_t$ by $H_T$. 
Then $\prod_tC_t$ can be 
identified with the algebra of bounded linear operators $c$ on $H_T$ satisfying
that $c$ preserves $H_s$ for each $s\in T$ and the restriction of $c$ on $H_s$ is in
$C_s$ for each $s\in T$. Now $C\otimes D$ is naturally a $C^*$-algebra of bounded
linear operators on the Hilbert space tensor product $H_T\otimes K=\oplus_{t\in T}(H_t\otimes K)$. 
It is easily checked that for every $g\in C\otimes D$,
$g$ preserves $H_s\otimes K$ for each $s\in T$, the restriction of $c$ on $H_s$ is equal to $\pi'_s(\varphi(g))$
for each $s\in T$, and the function $t\mapsto \pa \pi'_t(\varphi(g))\pa$ on $T$ vanishes at infinity (check this
for $g'\in C\odot D$ first, then approximate $g\in C\otimes D$ by $g'\in C\odot D$). It follows that $\varphi$ is
injective and maps $C\otimes D$ into $\prod^{\sim}_t(C_t\otimes D)$. 

Clearly the restriction of $\pi'_s$ on $\varphi(C\otimes D)$ is onto $C_s\otimes D$ for each $s\in T$. 
Since $C$ is a $C_{\infty}(T)$-module, $\varphi(C\odot D)$ is easily seen to be a $C_{\infty}$-submodule
of $\prod_t(C_t\otimes D)$. It follows that $\varphi(C\otimes D)$ is a $C_{\infty}$-submodule
of $\prod_t(C_t\otimes D)$. Thus the pair
$(\{C_t\otimes D\}_{t\in T}, \varphi(C\otimes D))$ is a
topological field of $C^*$-algebras over $T$.
\end{proof}

From now on, for a topological field $(\{C_t\}_{t\in T}, C)$ 
of $C^*$-algebras over a locally compact Hausdorff space $T$
and a $C^*$-algebra $D$, we shall take $(\{C_t\otimes D\}_{t\in T}, C\otimes D)$ 
to be the
topological field of $C^*$-algebras over $T$ in Lemma~\ref{tensor product field:lemma}.

In general, for a continuous field $(\{C_t\}_{t\in T}, C)$ of $C^*$-algebras over a compact metrizable space
$T$ and a $C^*$-algebra $D$, the topological field $(\{C_t\otimes D\}_{t\in T}, C\otimes D)$ of $C^*$-algebras may fail
to be continuous \cite[Theorem A]{KW}. The following lemma tells us that if a field $(\{C_t\}_{t\in T}, C)$
over  a locally compact Hausdorff space $T$ 
can be {\it subtrivialized} in the sense that there is a $C^*$-algebra $B$ containing each $C_t$ as a $C^*$-subalgebra
so that the elements of $C$ are exactly the continuous maps $T\rightarrow B$ vanishing at $\infty$ whose images at each $t$ are in $C_t$,
then the field $(\{C_t\otimes D\}_{t\in T}, C\otimes D)$ can also be subtrivialized and hence is continuous.

\begin{lemma} \label{tensor product:lemma}
Let $(\{C_t\}_{t\in T}, C)$ be a topological field of $C^*$-algebras over a locally compact Hausdorff space $T$. 
Suppose that there is a $C^*$-algebra $B$ containing each $C_t$ as a $C^*$-subalgebra
so that the elements of $C$ are exactly the continuous maps $T\rightarrow B$ vanishing at $\infty$ whose images at each $t$ are in $C_t$.
Let $D$ be a $C^*$-algebra, and identify $C\otimes D$ with a $C^*$-subalgebra of
$\prod^{\sim}_t(C_t\otimes D)$ as in Lemma~\ref{tensor product field:lemma}. Then
elements of $C\otimes D$ are exactly 
the continuous maps $T\rightarrow B\otimes D$ vanishing at $\infty$ whose images at each $t$ are
in $C_t\otimes D$.
\end{lemma}
\begin{proof}

Denote by $W$ the continuous maps $T\rightarrow B\otimes D$ vanishing at $\infty$ whose images at each $t$ are
in $C_t\otimes D$. This is a $C^*$-subalgebra of $\prod^{\sim}_t(C_t\otimes D)$. 

Denote by $\pi'_t$ the coordinate map $C\otimes D\rightarrow C_t\otimes D$ for each $t\in T$. Then
$\pi'_t(f\otimes d)=f(t)\otimes d$ for all $t\in T$, $f\in C$, and $d\in D$. It is
easy to check that $C\odot D\subseteq W$. Thus $C\otimes D\subseteq W$. 

Let $w\in W$ and let $\varepsilon>0$. For any $s\in T$, we can find some $\sum_j b_j\otimes d_j\in C_s\odot D$
satisfying $\pa w(s)- \sum_j b_j\otimes d_j\pa<\varepsilon$. 
Take $f_j\in C$ with
$f_j(s)=b_j$. Then $\pa w(t)-(\sum_j f_j\otimes d_j)(t)\pa <\varepsilon$ for $t=s$ and hence
for all $t$ in some neighborhood of $s$ by continuity. Note that both $C$ and $W$ are Banach modules
over $C_{\infty}(T)$.
Now a standard partition of unity argument
shows that we can find some $g\in C\odot D$ with $\pa w-g\pa <\varepsilon$. Thus $C\otimes D$ is
dense in $W$ and hence $C\otimes D=W$.
\end{proof}

Next we discuss semi-continuous fields of ergodic actions of
$\cG$.
The following definition
is a natural generalization of
Rieffel's definition of upper semi-continuous fields of actions of
locally compact groups \cite[Definition 3.1]{Rieffel89}. 

\begin{definition} \label{continuous field of actions:def}
By a {\it topological field
of actions} of $\cG$ on unital $C^*$-algebras we mean
a topological field $(\{B_t\}_{t\in T}, B)$ of unital $C^*$-algebras over a locally compact Hausdorff space $T$,
and an action $\sigma_t$ of $\cG$ on $B_t$
for each $t\in T$ such that
the section $\{f(t)1_{B_t}\}_{t\in T}$ is in $B$ for
each $f\in C_{\infty}(T)$ and
$\{\sigma_t\}_{t\in T}$ is a homomorphism
from $(\{B_t\}_{t\in T}, B)$ to $(\{B_t\otimes A\}_{t\in T}, B\otimes A)$.
If the field $(\{B_t\}_{t\in T}, B)$ is actually upper semi-continuous (lower semi-continuous, continuous, resp.),
then we will say that the field of actions is {\it upper semi-continuous} ({\it lower semi-continuous}, {\it continuous}, resp.).
If each $\sigma_t$ is ergodic, we say that this is a {\it field of ergodic
actions}.
\end{definition}

Clearly the pull-back of a topological (upper semi-continuous, lower semi-continuous, continuous, resp.)
field of
actions of $\cG$ on unital $C^*$-algebras is a topological (upper semi-continuous, lower semi-continuous, continuous, resp.)
field of actions of $\cG$.


\begin{lemma} \label{passing to red field:lemma}
Let $(\{(B_t, \sigma_t)\}_{t\in T}, B)$ be a semi-continuous field of
ergodic actions of $\cG$ over a locally compact Hausdorff space
$T$. Then for any $b\in B$ the function $t\mapsto
\omega_t(\pi_t(b))$ is continuous on $T$, where $\omega_t$ is the
unique $\sigma_t$-invariant state on $B_t$. Denote by $(B_{t,
r}, \sigma_{t, r})$ the reduced action associated to $(B_t,
\sigma_t)$ and by $\pi_{t, r}$ the canonical $*$-homomorphism
$B_t\rightarrow B_{t, r}$. Denote by $\pi_r$ the
$*$-homomorphism $\prod_tB_t\rightarrow \prod_tB_{t, r}$ given
pointwisely by $\pi_{t, r}$. Then $(\{(B_{t, r}, \sigma_{t, r})\}_{t\in T},
\pi_r(B))$ is a lower semi-continuous field of ergodic actions
of $\cG$ over $T$.
\end{lemma}
\begin{proof}
We prove the continuity of the function  $t\mapsto \omega_t(\pi_t(b))$ first.
Via taking restrictions to compact subsets of $T$ we may assume that $T$ is compact.
The cross-section $t\mapsto \omega_t(\pi_t(b))1_{B_t}$ is simply $((\id\otimes h)\circ(\prod_t\sigma_t))(b)$, which is in
$B$. Thus the function $t\mapsto \omega_t(\pi_t(b))$ is continuous by Lemma~\ref{cont field of unit:lemma}.

Next we show that $(\{(B_{t, r}, \sigma_{t, r})\}_{t\in T}, \pi_r(B))$ is a
lower semi-continuous field of actions. Clearly $\pi_r(B)$ is a
$C^*$-subalgebra and $C_{\infty}(T)$-submodule of
$\prod^{\sim}_tB_{t, r}$, and the evaluation map
$\pi_t:\pi_r(B)\rightarrow B_{t, r}$ is surjective for each
$t$. Since $\prod_t\sigma_{t, r}\circ \prod_t\pi_{t,
r}=\prod_t(\pi_{t, r}\otimes \id) \circ \prod_t\sigma_t$, one sees
that
$\prod_t\sigma_{t, r}$
sends $\pi_r(B)$ into $\pi_r(B)\otimes A$. We are left to
show that the function $t\mapsto \pa \pi_t(\pi_r(b))\pa$ is lower
semi-continuous for each $b\in B$. Note that for any $b\in B$
and $t\in T$, the norm of $\pi_t(\pi_r(b))$ is the smallest number
$K$ such that $\omega_t(\pi_t(b^*_1b^*bb_1))^{\frac{1}{2}}\le
K\omega_t(\pi_t(b^*_1b_1))^{\frac{1}{2}}$ for all $b_1\in B$. It
follows easily that the function $t\mapsto \pa \pi_t(\pi_r(b))\pa$
is lower semi-continuous over $T$ for each $b\in B$. This
completes the proof of Lemma~\ref{passing to red field:lemma}.
\end{proof}

It is well-known that there is a continuous field of ergodic
actions of the $n$-dimensional torus $\Te^d$
over the compact space of isomorphism classes
of faithful ergodic actions of $\Te^d$ such that the isomorphism
class of the fibre at each point is exactly the point (see
\cite[Theorem 1.1]{AP} for
a proof for the case $n=2$; the proof for the higher-dimensional
case is similar). We have not
been able to extend this to arbitrary compact quantum groups. What
we find is that there are two natural semi-continuous fields of
ergodic actions of $\cG$ over $\cP$ such that the equivalence
class of the fibre at each $t\in \cP$ is the image of $t$ under
the quotient map $\cP\rightarrow \EA^{\sim}(\cG)$ defined before
Definition~\ref{topology on EA:def}. By Propositions~\ref{P:prop}
and \ref{coefficient to ergodic:prop}, for each $t\in \cP$, the
pair $(\sV_t, \sigma_t)$ defined after the formula
(\ref{coefficient52:eq})
is isomorphic to the regular part of
some ergodic action of $\cG$. By Propositions~\ref{universal
action:prop} and \ref{reduced:prop} there exist (unique up to
isomorphisms) a full action $(B_{t, u},\sigma_{t, u})$ and a
reduced action $(B_{t, r}, \sigma_{t, r})$ of $\cG$ whose
regular parts are exactly $(\sV_t, \sigma_t)$. In fact, one can
take $(B_t, \sigma_t)$ in Proposition~\ref{coefficient to
ergodic:prop} as $(B_{t, u}, \sigma_{t, u})$. Recall the
the quotient map $\phi_t: \sV\rightarrow \sV_t$ defined before (\ref{coefficient12:eq})
for each $t\in
\cP$.

\begin{theorem} \label{cont field over P:thm}
The set of cross-sections $\{\phi_t(v)\}_{t\in \cP}$ over $\cP$ for $v\in
\sV$ is in $\prod_tB_{t, u}$ ($\prod_tB_{t, r}$ resp.). It
generates an upper (lower resp.) semi-continuous field $(\{B_{t,
u}\}_{t\in \cP}, B_u)$ ($(\{B_{t, r}\}_{t\in \cP}, B_r)$ resp.) of $C^*$-algebras
over $\cP$. Moreover, the field $(\{(B_{t, u}, \sigma_{t, u})\}_{t\in \cP},
B_u)$ ($(\{(B_{t, r}, \sigma_{t, r})\}_{t\in \cP}, B_r)$ resp.) is an
upper (lower resp.) semi-continuous field of full (reduced resp.)
ergodic actions of $\cG$. If $\cG$ is co-amenable, then these two
fields coincide and are continuous.
\end{theorem}
\begin{proof}
Consider generators $w_x$ for $x\in X_0$,
 $\theta(y)$ for $y\in Y$ and $\zeta(z)$ for $z\in Z$ subject to the following relations:
\begin{enumerate}
\item $w_{x_0}$ is the identity,

\item the equations (\ref{coefficient12:eq}) and (\ref{coefficient32:eq})
with $\nu_x, f(y), g(z)$ replaced by $w_x$, $\theta(y)$, $\zeta(z)$ respectively,

\item the equations in $\sE$ with $f(y), \overline{f(y)}, g(z), \overline{g(z)}$ replaced
by $\theta(y)$, $\theta(y)^*$, $\zeta(z)$, $\zeta(z)^*$ respectively,

\item $\theta(y)$ and $\zeta(z)$ are in the center.
\end{enumerate}
These relations have $*$-representations since $B_{f, g, u}$ for
any $(f, g)\in \cP$ has generators satisfying these conditions.
Consider an irreducible representation $\pi$ of these relations.
Because of (4), $\pi(\theta(y))$ and $\pi(\zeta(z))$ have to be
scalars. Say $\pi(\theta(y))=f(y)$ and $\theta(\zeta(z))=g(z)$.
Then $(f, g)\in \Ce^Y\times \Ce^Z$ satisfies the equations in
$\sE$ because of (3). Thus the inequalities (\ref{compact:eq1})
and (\ref{compact:eq2}) hold with $|f(y)|$ and $|g(z)|$ replaced
by $\pa \pi(\theta(y))\pa$ and $\pa \pi(\zeta(z))\pa$
respectively. Also, there is a $*$-homomorphism from $\sV_{f, g}$
to the $C^*$-algebra generated by $\pi(w_x), \pi(\theta(y)),
\pi(\zeta(z))$ sending $\nu_x$ to $\pi(w_x)$.  Thus
(\ref{bound2:eq}) holds with $\nu_x$ replaced by $\pi(w_x)$.
Consequently, above generators and relations do have a universal
$C^*$-algebra $B_u$.

In particular, there is a surjective $*$-homomorphism $\pi_{f,
g}:B_u\rightarrow B_{f, g, u}$ for each $(f, g)\in \cP$
sending $w_x$, $\theta(y)$, $\zeta(z)$ to $\phi_{f, g}(v_x)$,
$f(y)\phi_{f, g}(v_{x_0})$, $g(z)\phi_{f, g}(v_{x_0})$
respectively. These $*$-homomorphisms $\pi_t$'s for $t\in \cP$
combine to a $*$-homomorphism $\pi:B_u\rightarrow \prod_tB_{t,
u}$. In above we have seen that every irreducible
$*$-representation of $B_u$ factors through $\pi_t$ for some
$t\in \cP$. Thus $\pi$ is faithful and we may identify $B_u$
with $\pi(B_u)$. Since $B_{f, g, u}$ is the universal
$C^*$-algebra of $\sV_{f, g}$, one sees easily that $\ker(\pi_{f,
g})$ is generated by $\theta(y)-f(y)w_{x_0}$ and
$\zeta(z)-g(z)w_{x_0}$.  Since $\theta(y)-f'(y)w_{x_0}\rightarrow
\theta(y)-f(y)w_{x_0}$ and $\zeta(z)-g'(z)w_{x_0}\rightarrow
\zeta(z)-g(z)w_{x_0}$ as $(f', g')\rightarrow (f, g)$, the
function $t\mapsto \pa \pi_t(b)\pa$ is upper semi-continuous on
$\cP$ for each $b\in B_u$. Thanks to the Stone-Weierstrass
theorem, the unital $C^*$-subalgebra of $B_u$ generated by
$\theta(y)$ and $\zeta(z)$ is exactly $C(\cP)$. Thus $B_u$ is a
$C(\cP)$-submodule of $\prod_tB_{t, u}$. Therefore $(\{B_{t,
u}\}_{t\in \cP}, B_u)$ is an upper semi-continuous field of $C^*$-algebras
over $\cP$. Clearly it is generated by the sections
$\{\phi_t(v)\}_{t\in \cP}$ for $v\in \sV$.

The formula (\ref{coefficient32:eq}) tells us that $\prod_t\sigma_{t, u}$ sends the section $\{\phi_t(v_x)\}_{t\in \cP}$
into $B_u\otimes A$ for each $x\in X_0$. Since $B_u$ is generated by such sections and $C(\cP)$,
$\prod_t(\sigma_{t, u})$ sends $B_u$ into $B_u\otimes A$. Thus
$(\{(B_{t, u}, \sigma_{t, u})\}_{t\in \cP}, B_u)$ is an upper semi-continuous
field of ergodic actions of $\cG$.

The assertions about the reduced actions follow from Lemma~\ref{passing to red field:lemma}.
The assertion about the case $\cG$ is co-amenable follows from Proposition~\ref{universal:prop}.
\end{proof}

\begin{theorem} \label{cont vs mul:thm}
Let $(\{(B_t, \sigma_t)\}_{t\in T}, B)$ be a semi-continuous field
of ergodic actions of $\cG$  over a locally compact Hausdorff space $T$.
Let $t_0\in T$. Then the following are equivalent:
\begin{enumerate}
\item the map $T\rightarrow \EA(\cG)$ sending each $t$ to
the isomorphism class of
$(B_t, \sigma_t)$ is continuous at $t_0$,

\item  the map $T\rightarrow \EA^{\sim}(\cG)$ sending each $t$ to
the equivalence class of  $(B_t, \sigma_t)$
is continuous at $t_0$,

\item $\limsup_{t\rightarrow t_0}\mul(B_t, \gamma)\le \mul(B_{t_0}, \gamma)$ for
all $\gamma \in \hat{\cG}$,

\item $\lim_{t\rightarrow t_0}\mul(B_t, \gamma)= \mul(B_{t_0}, \gamma)$ for
all $\gamma \in \hat{\cG}$.
\end{enumerate}
\end{theorem}

\begin{lemma} \label{lifting ON:lemma}
Let the notation be as in Theorem~\ref{cont vs mul:thm}.
Let $\gamma \in \hat{\cG}$, and
let $c_{\gamma si}, 1\le i\le \mul(B_{t_0}, \gamma), 1\le i\le d_{\gamma}$ be
a standard basis of $(B_{t_0})_{\gamma}$. Then there is a linear map
$\varphi_t:(B_{t_0})_{\gamma}\rightarrow (B_{t})_{\gamma}$ for all $t\in T$
such that the section $t\mapsto \varphi_t(c)$ is in $B$ for every $c\in (B_{t_0})_{\gamma}$,
that $\varphi_{t_0}=\id$, and that $\varphi_t(c_{\gamma si}), 1\le s\le \mul(B_{t_0}, \gamma),
1\le i\le d_{\gamma}$
satisfy (\ref{action:eq}) and (\ref{ON:eq}) (with $e_{\gamma si}$ and $\omega$ replaced by $\varphi_t(c_{\gamma si})$
and the unique $\sigma_t$-invariant state $\omega_t$ respectively)
throughout a neighborhood of $t_0$.
\end{lemma}
\begin{proof} 
We may assume that $T$ is compact.
Denote by $\sigma$ the restriction of $\prod_t\sigma_t$ on $B$.
Recall the map $E^{\gamma}_{ij}$ defined via (\ref{E_ij:eq}).
Then $E^{\gamma}_{ij}$ is also defined on $B$ for the unital $*$-homomorphism
$\sigma:B\rightarrow B\otimes A$.
Set $m=\mul(B_{t_0}, \gamma)$ and $S=\{c_{\gamma s1}: 1\le s\le m\}$.
For each $c\in S$ take  $b\in B$ with $\pi_{t_0}(b)=c$.
Then $\pi_t(E^{\gamma}_{11}(b))=E^{\gamma}_{11}(\pi_t(b))$ is in $E^{\gamma}_{11}(B_t)$
for each $t\in T$.
By Lemma~\ref{change of basis:lemma}
$S$ is a linear basis of $E^{\gamma}_{11}(B_{t_0})$.
Set $\psi_t$ to be the linear map $E^{\gamma}_{11}(B_{t_0})\rightarrow E^{\gamma}_{11}(B_t)$
sending each $c\in S$ to $\pi_t(E^{\gamma}_{11}(b))$.
By Lemma~\ref{change of basis:lemma} we have $\psi_{t_0}=\id$.
By Lemma~\ref{passing to red field:lemma}
the function $t\mapsto \omega_t(\pi_t(b'))$ is continuous on $T$ for
any $b'\in B$, where $\omega_t$ is the unique $\sigma_t$-invariant state
on $B_t$. Consequently, for any $c_1, c_2\in S$, we have
\begin{eqnarray*}
\omega_t({\psi_t(c_1)}^*\psi_t(c_2))\rightarrow \omega_{t_0}({\psi_{t_0}(c_1)}^*\psi_{t_0}(c_2))=
\omega_{t_0}(c_1c_2)=\delta_{c_1c_2}
\end{eqnarray*}
as $t\rightarrow t_0$.
Shrinking $T$ if necessary, we may assume that the matrix
$Q_t=(\omega_t({\psi_t(c_{\gamma k1})}^*\psi_t(c_{\gamma s1})))_{ks}\in
M_m(\Ce)$ is invertible for all $t\in T$.  Set $(c_{t,1}, {\cdots}, c_{t, m})=
(\psi_t(c_{\gamma 11}), {\cdots}, \psi_t(c_{\gamma m1}))Q^{-\frac{1}{2}}_t$.
Then $c_{t, k}\in  E^{\gamma}_{11}(B_t)$ and $\omega_t(c^*_{t, k}c_{t, s})=\delta_{ks}$
for all $t\in T$. Note that the section $t\mapsto c_{t, s}$ is in $B$ for
each $1\le s\le m$. Thus  the section $t\mapsto E^{\gamma}_{i1}(c_{t, s})$ is in $B$
for all $1\le s\le m, 1\le i\le d_{\gamma}$. Set $\varphi_t$ to be the
linear map $(B_{t_0})_{\gamma}\rightarrow B_t$
sending $c_{\gamma si}$ to $E^{\gamma}_{i1}(c_{t, s})$.
Then the section $t\mapsto \varphi_t(c)$ is in $B$ for every $c\in (B_{t_0})_{\gamma}$.
By Lemma~\ref{change of basis:lemma} these maps have the other desired properties.
\end{proof}

\begin{remark} \label{upper:remark}
Using Remark~\ref{unique field:remark} one can show easily that for an
upper semi-continuous field $(\{(B_t, \sigma_t)\}_{t\in T}, B)$ of
actions of $\cG$ over a compact Hausdorff space $T$, the
$*$-homomorphism $(\prod_t\sigma_t)|_{B}:B\rightarrow
B\otimes A$ is an action of $\cG$ on $B$. Using the
well-known fact that upper semi-continuous fields of unital
$C^*$-algebras over a compact Hausdorff space $T$ satisfying the
hypothesis in Lemma~\ref{cont field of unit:lemma} correspond
exactly to unital $C^*$-algebras containing $C(T)$ in the centers,
one can show further that upper semi-continuous fields of ergodic
actions of $\cG$ over $T$ correspond exactly to actions of $\cG$
on unital $C^*$-algebras whose fixed point algebras are $C(T)$ and
are in the centers.
\end{remark}

As a corollary of Lemma~\ref{lifting ON:lemma} we get

\begin{lemma} \label{mul lower cont:lemma}
Let the notation be as in Theorem~\ref{cont vs mul:thm}.
The function $t\mapsto \mul(B_t, \gamma)$ is lower semi-continuous on $T$
for each $\gamma \in \hat{\cG}$.
\end{lemma}

We are ready to prove Theorem~\ref{cont vs mul:thm}.

\begin{proof}[Proof of Theorem~\ref{cont vs mul:thm}]
(1)$\iff$(2) follows from the definition of the topology on $\EA(\cG)$.
(2)$\Rightarrow$(3) follows from Proposition~\ref{mul:prop}.
(3)$\Rightarrow$(4) follows from Lemma~\ref{mul lower cont:lemma}.
We are left to show (4)$\Rightarrow$(2). Assume (4).
Fix a standard basis $\sS$
of $\cB_{t_0}$, consisting of a standard basis $\sS_{\gamma}$ of $(B_{t_0})_{\gamma}$
for each $\gamma\in \hat{\cG}$.
Let $J$ be a finite subset of $\hat{\cG}$.
Then $\mul(B_t, \gamma)=\mul(B_{t_0}, \gamma)$ for each $\gamma\in J$ throughout
some neighborhood $U$ of $t_0$.
By Lemma~\ref{lifting ON:lemma}, shrinking $U$ if necessary,
 we can find a linear map $\varphi_t:(B_{t_0})_{J}\rightarrow (B_{t})_{J}$
for all $t\in T$, where $(B_{t})_{J}=\sum_{\gamma\in J}(B_t)_{\gamma}$,
such that the section $t\mapsto \varphi_t(c)$ is in $B$ for every $c\in (B_{t_0})_{J}$,
that $\varphi_{t_0}=\id$, and
that $\varphi_t(\sS_{\gamma})$ is
a standard basis of $(B_t)_{\gamma}$ for all $\gamma \in J$ and $t\in U$.
For each $t\in U$, extend these bases of $(B_t)_{\gamma}$ for $\gamma \in J$ to a standard basis
$\sS_t$ of $\cB_t$. Set $(f_t, g_t)$ to
be the element in $\cP$ associated to $\sS_t$ via (\ref{coefficient1:eq})-(\ref{coefficient2:eq}).
Suppose that $\alpha, \beta \in J\setminus \{\gamma_0\}$.
By Lemma~\ref{passing to red field:lemma}
the function $t\mapsto \omega_t(\pi_t(b))$ is continuous
for each $b\in B$, where $\omega_t$ is the unique $\sigma_t$-invariant state
on $B_t$.
Then one sees easily that the function $t\mapsto f_t(x_1, x_2, x_3)$
is continuous over $U$ for any $x_1\in X_{\alpha}, x_2\in X_{\beta}, x_3\in X_{\gamma}, \gamma \in J$.
Similarly, if $\alpha, \bar{\alpha}\in J\setminus \{\gamma_0\}$, then the function $t\mapsto
g_t(x_1, x_2)$
is continuous over $U$ for any $x_1\in X_{\alpha}, x_2\in X_{\bar{\alpha}}$.
Since $J$ is an arbitrary finite subset of $\hat{\cG}$, this means that for any neighborhood $W$ of $(f_{t_0}, g_{t_0})$
in $\cP$, we can find a neighborhood $V$ of $t_0$ in $T$ and choose a standard basis of
$\cB_t$ for each $t\in V$ such that the associated element in $\cP$ is in $W$.
Therefore (2) holds.
\end{proof}

Now Theorems~\ref{main1:thm} and \ref{main2:thm} follow from
Theorems~\ref{topology on EA:thm}, \ref{cont vs mul:thm} and
\ref{cont field over P:thm}. In fact we have a stronger assertion:

\begin{corollary} \label{unique topology:coro}
The topology on $\EA^{\sim}(\cG)$ defined in Definition~\ref{topology
on EA:def} is the unique Hausdorff topology on  $\EA^{\sim}(\cG)$ such
that the implication (4)$\Rightarrow$(2) in Theorem~\ref{cont vs
mul:thm} holds for all upper semi-continuous (lower
semi-continuous resp.) fields of ergodic actions of $\cG$ over
compact Hausdorff spaces. If $\cG$ is co-amenable, then the
topology on $\EA(\cG)$ defined in Definition~\ref{topology on EA:def}
is the unique Hausdorff topology on  $\EA(\cG)$ such that the
implication (4)$\Rightarrow$(2) in Theorem~\ref{cont vs mul:thm}
holds for all continuous fields of ergodic actions of $\cG$ over
compact Hausdorff spaces.
\end{corollary}

When $A$ is separable and co-amenable, one can describe the
topology on $\EA(\cG)$ more explicitly in terms of continuous fields of
actions:

\begin{theorem} \label{separable:thm}
Suppose that $A$ is separable and co-amenable. Then both $\cP$
and $\EA(\cG)$ are metrizable. The isomorphism classes of a sequence
$\{(B_n, \sigma_n)\}_{n\in \Ne}$ of ergodic actions of $\cG$
converge  to that of $(B_{\infty}, \sigma_{\infty})$ in $\EA(\cG)$ if
and only if there exists a continuous field of ergodic actions of
$\cG$ over the one-point compactification $\Ne\cup \{\infty\}$ of
$\Ne$ with fibre $(B_n, \sigma_n)$ at $n$ for $1\le n\le \infty$
and $\lim_{n\rightarrow \infty}\mul(B_n,
\gamma)=\mul(B_{\infty}, \gamma)$ for all $\gamma\in \hat{\cG}$.
\end{theorem}
\begin{proof} Denote by $(\pi_{A}, H_{A})$ the $\GNS$
representation of $A$ associated to $h$. Since $A$ is
separable, so is $H_{A}$. Note that the subspaces
$A_{\gamma}$ are nonzero and orthogonal to each other in
$H_{A}$ for $\gamma \in \hat{\cG}$.
Thus $\hat{\cG}$ is countable. Then $Y$ and $Z$ are both
countable. Therefore $\cP$ and $\EA(\cG)$ are metrizable. The ``if''
part follows from Theorem~\ref{cont vs mul:thm}. Suppose that
the isomorphism class of $(B_n, \sigma_n)$ converges to that of
$(B_{\infty}, \sigma_{\infty})$ in $\EA(\cG)$ as $n\to \infty$. By
Proposition~\ref{mul:prop} we have $\lim_{n\rightarrow
\infty}\mul(B_n, \gamma)=\mul(B_{\infty}, \gamma)$ for all
$\gamma\in \hat{\cG}$. Also the map $\xi:\Ne\cup \{\infty\}
\rightarrow \EA(\cG)$ sending $1\le n\le \infty$ to the isomorphism
class of $(B_n, \sigma_n)$ is continuous. By
Theorem~\ref{topology on EA:thm} the quotient map $\cP\rightarrow
\EA(\cG)$ is open. Thus $\xi$ lifts up to a continuous map
$\eta:\Ne\cup \{\infty\}\rightarrow \cP$. The pull-back of the
continuous field
of ergodic actions of $\cG$ over $\cP$ in Theorem~\ref{cont field
over P:thm} via $\eta$ is a continuous field of ergodic actions of
$\cG$ over $ \Ne\cup \{\infty\}$ with the desired fibres. This
proves the ``only if'' part.
\end{proof}

\section{Podle\'s spheres}
\label{Podles:sec}

In this section we prove Theorem~\ref{cont of Podles:thm}.

Fix $q\in [-1, 1]$. The quantum $\SU(2)$ group
$A=C(\SU_q(2))$ \cite{VS, WorSU2} is defined as the universal $C^*$-algebra generated by
$\alpha$ and $\beta$ subject to the condition that
\begin{eqnarray*}
u=\left( \begin{matrix} \alpha & -q\beta^* \\ \beta &  \alpha^* \end{matrix}\right)
\end{eqnarray*}
is a unitary in $M_2(A)$. The comultiplication $\Phi:A\rightarrow A$ is defined in such a way that
$u$ is a representation of $A$.


Below we assume $0< |q|< 1$. The quantum group $\SU_q(2)$ is co-amenable
\cite{Nagy}\cite[Corollary 6.2]{Banica}\cite[Theorem 2.12]{BMT}.
Let
$$T_q=\{c(1), c(2), {\cdots}\}\cup [0, 1],$$
where
\begin{eqnarray*}
c(n)=-q^{2n}/(1+q^{2n})^2.
\end{eqnarray*}
For $t\in T_q$ with $t\le 0$, Podle\'s quantum sphere $C(S^2_{qt})$ \cite{Podles87} is defined
as the universal $C^*$-algebra generated by $a_t, b_t$ subject to the relations
\begin{eqnarray} \label{relation finite:eq}
a^*_t=a_t, \quad & & b^*_tb_t=a_t-a^2_t+t,  \\
\nonumber b_ta_t=q^2a_tb_t,\quad & & b_tb^*_t=q^2a_t-q^4a^2_t+t.
\end{eqnarray}
For $t\in T_q$ with $t\ge 0$, $C(S^2_{qt})$ is defined as the universal $C^*$-algebra generated by $a_t, b_t$
subject to the relations
\begin{eqnarray} \label{relation infinit:eq}
a^*_t=a_t, \quad b^*_tb_t=(1-t^2)a_t-a^2_t+t^2, \\
\nonumber b_ta_t=q^2a_tb_t,
\quad b_tb^*_t=(1-t^2)q^2a_t-q^4a^2_t+t^2.
\end{eqnarray}
The action $\sigma_t:C(S^2_{qt})\rightarrow C(S^2_{qt})\otimes A$ is determined by
\begin{eqnarray}   \label{action finite:eq}
\sigma_t(a_t)&=&a_t\otimes 1_{A}+c_t\otimes \beta^*\beta +b^*_t\otimes \alpha^*\beta
+b_t\otimes\beta^*\alpha,\\
\nonumber \sigma_t(b_t)&=&-qb^*_t\otimes \beta^2+c_t\otimes \alpha\beta+b_t\otimes\alpha^2,
\end{eqnarray}
where $c_t$ is $1_{C(S^2_{qt})}-(1+q^2)a_t$ or $(1-t^2)1_{C(S^2_{qt})}-(1+q^2)a_t$
depending on $t\le 0$ or $t\ge 0$.
As in \cite{HMS}, here we reparametrize the family for $0\le c\le \infty$ in \cite{Podles87}
for the parameters $0\le t\le 1$ by $t=2\sqrt{c}/(1+\sqrt{1+4c})$
(and $c=(t^{-1}-t)^{-2}$), and rescale the generators $A, B$ in \cite{Podles87}
by $a_t=(1-t^2)A$, $b_t=(1-t^2)B$ for $0\le t<1$.

\begin{proposition} \label{cont field of Podles:prop}
There is a unique continuous field
of $C^*$-algebras
over $T_q$ with fibre $C(S^2_{qt})$ at each $t\in T_q$ such that
the sections $t\mapsto a_t$ and $t\mapsto b_t$ are in the algebra $B$ of continuous sections.
Moreover, the field $\{\sigma_t\}_{t\in T_q}$ of ergodic actions of $\SU_q(2)$ is
continuous.
\end{proposition}
\begin{proof}
The uniqueness is clear.
We start to show that there exits an upper semi-continuous field $(\{C(S^2_{qt})\}_{t\in T_q}, B)$ of
$C^*$-algebras over $T_q$ such that the sections $t\mapsto a_t$ and $t\mapsto b_t$ are in $B$.
For this purpose, by Lemma~\ref{generate cont field:lemma}
it suffices to show that the function $\eta_p:t\mapsto \pa p(a_t, b_t, a^*_t, b^*_t)\pa$
is upper semi-continuous over $T_q$ for any noncommutative polynomial $p$ in four variables.
Denote by $T'_q$ the set of the non-positive numbers in $T_q$.
We prove the upper semi-continuity of $\eta_p$ over $T'_q$ first.

We claim that there exists a universal $C^*$-algebra generated by
$a, b, x$ subject to the relations
\begin{enumerate}
\item the equations in (\ref{relation finite:eq}) with $a_t, b_t, t$ replaced by
$a, b, x$ respectively,

\item the inequality $\pa x\pa \le |c(1)|$,

\item $x=x^*$ is in the center.
\end{enumerate}
Clearly $C(S^2_{qt})$ for $t\in T'_q$ has generators satisfying these conditions.
Let $a, b, x$ be bounded linear operators on a Hilbert space satisfying these relations.
We have
\begin{eqnarray} \label{bound norm:eq}
& &(1+q^2)(b^*b+q^{-2}bb^*)\\
\nonumber &\overset{(\ref{relation finite:eq})}=& (1+q^2)(a-a^2+x)+(1+q^2)(a-q^2a^2+q^{-2}x)\\
\nonumber &=&-(1-(1+q^2)a)^2+(1+(1+q^2)^2q^{-2}x).
\end{eqnarray}
Thus
\begin{eqnarray*}
\pa (1-(1+q^2)a)^2\pa,\, \pa (1+q^2)b^*b\pa  &\le &\pa  1+(1+q^2)^2q^{-2}x\pa \\
&\le &1+(1+q^2)^2q^{-2}|c(1)|=2,
\end{eqnarray*}
and hence
\begin{eqnarray*}
\pa a\pa &\le &(1+q^2)^{-1}(1+2^{\frac{1}{2}}), \\
\pa b\pa &\le &(1+q^2)^{-\frac{1}{2}}2^{\frac{1}{2}}.
\end{eqnarray*}
Therefore there does exist a universal $C^*$-algebra $C$ generated by
$a, b, x$ subject to these relations.
An argument similar to that in the proof of Theorem~\ref{cont field over P:thm}
shows that $\eta_p$
is upper semi-continuous over $T'_q$.

%
The upper semi-continuity
of $\eta_p$ over $[0, 1]$ is proved similarly,
replacing (\ref{bound norm:eq}) by
\begin{eqnarray} \label{bound norm:eq1}
& &(1+q^2)(b^*b+q^{-2}bb^*)\\
\nonumber &=&-((1-x^2)-(1+q^2)a)^2+((1-x^2)^2+(1+q^2)^2q^{-2}x^2).
\end{eqnarray}

This proves the existence of the desired upper semi-continuous
field of $C^*$-algebras over $T_q$. Note that $B$ is generated
as a $C^*$-algebra by $C(T_q)$ and the sections $t\mapsto a_t$ and
$t\mapsto b_t$. From (\ref{action finite:eq}) one sees immediately
that $(\{(C(S^2_{qt}), \sigma_t)\}_{t\in T_q}, B)$ is an upper
semi-continuous field of ergodic actions of $\SU_q(2)$. Since
$\SU_q(2)$ is co-amenable, by Proposition~\ref{universal:prop} and
Lemma~\ref{passing to red field:lemma} this is actually a
continuous field of actions.
\end{proof}

We are ready to prove Theorem~\ref{cont of Podles:thm}.

\begin{proof}[Proof of Theorem~\ref{cont of Podles:thm}]
It is customary to index $\widehat{\SU_q(2)}$ by $0, \frac{1}{2}, 1, 1+\frac{1}{2}, {\cdots}$ \cite[remark after the proof
of Theorem 5.8]{WorSU2}. Say $\widehat{\SU_q(2)}=\{\bd_0, \bd_{\frac{1}{2}}, \bd_1, \bd_{1+\frac{1}{2}}, {\cdots}\}$.
Then $\mul(C(S^2_{qt}), \bd_k)=1$, $\mul(C(S^2_{qt}), \bd_{k+\frac{1}{2}})=0$ for $k=0, 1, 2, {\cdots}$
when $t\ge 0$.  And $\mul(C(S^2_{qt}), \bd_l)=1$ or $0$ depending on $l\in \{0, 1, {\cdots}, n-1\}$
or not when $t=c(n)$ \cite[the note after Proposition 2.5]{Podles95}. Thus the multiplicity function
$t\mapsto \mul(C(S^2_{qt}), \gamma)$ is continuous over $T_q$ for
any $\gamma \in \widehat{\SU_q(2)}$. Then
Theorem~\ref{cont of Podles:thm}
follows from Proposition~\ref{cont field of Podles:prop} and Theorem~\ref{cont vs mul:thm}.
\end{proof}

\section{Ergodic actions of full multiplicity of compact groups}
\label{fm:sec}

In this section we show that the topology of Landstad and Wassermann on the
set $\EA(G)_{\fm}$ of isomorphism classes of ergodic actions of full multiplicity
of a compact group $G$ coincides with the relative topology of $\EA_{\fm}$
in $\EA(G)$.

Throughout this section we let $\cG=G$ be a compact Hausdorff group.
An ergodic action $(B, \sigma')$ of $G$ is said to be of {\it full multiplicity}
if $\mul(B, \gamma)=d_{\gamma}$ for all $\gamma \in \hat{G}$.
Denote by $\EA(G)_{\fm}$ the set of isomorphism classes of ergodic actions of full multiplicity
of $G$. By Proposition~\ref{mul:prop} $\EA(G)_{\fm}$ is a closed subset of $\EA(G)$.

Landstad \cite{Landstad} and Wassermann \cite{Wassermann} showed independently that
$\EA(G)_{\fm}$ can be identified with the set of  equivalence classes of dual cocycles.
Let us recall the notation in \cite{Landstad}. Denote by $\cL(G)$
the von Neumann algebra generated by the left regular representation of $G$
on $L^2(G)$.
One has a natural decomposition $L^2(G)\cong \oplus_{\gamma\in \hat{G}}H_{\gamma}$
as unitary representations of $G$.
Then $\cL(G)=\prod_{\gamma \in \hat{G}}B(H_{\gamma})$ under this decomposition.
Denote by $1_{\gamma_0}$ the identity of $B(H_{\gamma_0})$ for the trivial representation
$\gamma_0$.
One has the normal $*$-homomorphism $\delta:\cL(G)\rightarrow \cL(G)\otimes \cL(G)$
(tensor product of von Neumann algebras)
and the normal $*$-anti-isomorphism $\nu:\cL(G)\rightarrow \cL(G)$ determined by
\begin{eqnarray*}
 \delta(x)=x\otimes x \quad \mbox{ and } \quad \nu(x)=x^{-1} \quad \mbox{ for } x\in G.
\end{eqnarray*}
Denote by $\sigma$ the flip automorphism of $\cL(G)\otimes \cL(G)$ determined by $\sigma(a\otimes b)=b\otimes a$
for all $a, b\in \cL(G)$.
One also has Takesaki's unitary $W$ in $B(L^2(G))\otimes \cL(G)$ defined by
$ (Wf)(x, y)=f(x, xy)$ for $f\in C(G\times G)$ and $x, y\in G$.
A {\it normalized dual cocycle} \cite[page 376]{Landstad}
is a unitary $w\in \cL(G)\otimes \cL(G)$ satisfying
\begin{eqnarray*}
& &(w\otimes I)((\delta\otimes \id)(w))=(I\otimes w)((\id \otimes \delta)(w)), \\
& &(\nu\otimes \nu) (w)=\sigma(w^*), \, \qquad \qquad w(I\otimes 1_{\gamma_0})=I\otimes 1_{\gamma_0}, \\
& &w(1_{\gamma_0}\otimes I)=1_{\gamma_0}\otimes I, \quad   \qquad \qquad w\delta(1_{\gamma_0})=\delta(1_{\gamma_0}),\\
& &(\id\otimes \nu)(w \sigma(w^*))=\sigma(w)w^*, \quad  (\id\otimes \nu)(w W^*)=Ww^*.
\end{eqnarray*}
Denote by $C^2$ the set of all normalized dual cocycles.
Also denote by $H$ the group of unitaries $\xi$ in $\cL(G)$ satisfying $\xi=\nu(\xi^*)$ and $\xi 1_{\gamma_0}=1_{\gamma_0}$
(on page 376 of \cite{Landstad} only the condition $\xi=\nu(\xi^*)$ is mentioned, but in order for
$\alpha_{\xi}(w)$ below to satisfy $\alpha_{\xi}(w)(I\otimes 1_{\gamma_0})=I\otimes 1_{\gamma_0}$,
one has to require $\xi 1_{\gamma_0}=1_{\gamma_0}$;
this can be seen using the formula $\delta(x)(I\otimes 1_{\gamma_0})=x\otimes 1_{\gamma_0}$ for all $x\in \cL(G)$).
Then $H$ has a left action $\alpha$ on $C^2$ via $\alpha_{\xi}(w)=(\xi\otimes \xi)w \delta(\xi^*)$.
The result of Landstad and Wassermann says that $\EA(G)_{\fm}$ can be identified with $C^2/H$ \cite[Remark 3.13]{Landstad}
in a natural way.

Note that the unitary groups of $\cL(G)\otimes \cL(G)$
and $\cL(G)$ are both compact Hausdorff groups with the weak topology.
Clearly $C^2$ and $H$ are closed subsets of the unitary groups of $\cL(G)\otimes \cL(G)$
and $\cL(G)$ respectively.
Thus
$C^2$ is a compact Hausdorff space and $H$ is a compact Hausdorff group,
with the relative topologies. It is also clear that the action $\alpha$ is continuous.
Therefore $C^2/H$ equipped with the quotient topology is a compact Hausdorff space.

In order to show that the quotient topology on $C^2/H$ coincides with the relative topology
of $\EA(G)_{\fm}$ in $\EA(G)$, we need to recall the map $C^2\rightarrow \EA(G)_{\fm}$ constructed in the proof
of \cite[Theorem 3.9]{Landstad}.
Let $w \in C^2$. Set $U=w W^*\in B(L^2(G))\otimes \cL(G)$. Recall that for each $\gamma \in \hat{G}$
we fixed
an orthonormal basis of $H_{\gamma}$ and identified $B(H_{\gamma})$ with $M_{d_{\gamma}}(\Ce)$.
Let $e^{\gamma}_{ij}, 1\le i, j\le d_{\gamma}$ be the matrix units of $M_{d_{\gamma}}(\Ce)$ as usual.
Then we may write $U$ as $\sum_{\gamma\in \hat{G}}\sum_{1\le i, j\le d_{\gamma}}b_{\gamma ij}\otimes e^{\gamma}_{ij}$
for $b_{\gamma ij}\in B(L^2(G))$. The  conjugation of the right regular representation of $G$ on $L^2(G)$
restricts on an ergodic action $\alpha$ of $G$ on the
$C^*$-algebra $B$ generated by $b_{\gamma ij}$ for all $\gamma \in \hat{G}, 1\le i, j\le d_{\gamma}$.
The isomorphism class of $\alpha$ is the image of $w$ under the map $C^2\rightarrow \EA(G)_{\fm}$.
Furthermore, each $U_{\gamma}=\sum_{1\le i, j\le d_{\gamma}}b_{\gamma ij}\otimes e^{\gamma}_{ij}$
is a {\it unitary $u^{\gamma}$-eigenoperator} meaning that $U_{\gamma}$ is a unitary in
$B\otimes B(H_{\gamma})$ satisfying
\begin{eqnarray} \label{eigen:eq}
(\alpha_x\otimes \id)(U_{\gamma})=U_{\gamma}(1_{B}\otimes u^{\gamma}(x))
\end{eqnarray}
for all $x\in G$. If we let $\sigma:B\rightarrow B\otimes C(G)=C(G, B)$ be the $*$-homomorphism
associated to $\alpha$, i.e., $(\sigma(b))(x)=\alpha_x(b)$, then
(\ref{eigen:eq}) simply means $(\sigma\otimes \id)(U_{\gamma})=(U_{\gamma})_{13}(\tau(u^{\gamma}))_{23}$,
where $(U_{\gamma})_{13}$ and $(\tau(u^{\gamma}))_{23}$ are in the leg numbering notation
and $\tau:B(H_{\gamma})\otimes C(G)\rightarrow C(G)\otimes B(H_{\gamma})$ is the flip.
It follows that (\ref{eigen:eq}) is equivalent to (\ref{action:eq}) with
$e_{\gamma ki}$ replaced by $b_{\gamma ki}$.
Then $\sum_{1\le j\le d_{\gamma}}b_{\gamma ij}b^*_{\gamma kj}$ is easily seen to be $\sigma$-invariant and
hence is in $\Ce 1_{B}$. One checks easily that $U_{\gamma}U^*_{\gamma}=1_{B}\otimes 1_{B(H_{\gamma})}$ means
that
\begin{eqnarray} \label{on:eq}
\omega(\sum_{1\le j\le d_{\gamma}}b_{\gamma ij}b^*_{\gamma kj})=\delta_{ki}
\end{eqnarray}
for all $1\le k, i\le d_{\gamma}$, where $\omega$ is the unique $\alpha$-invariant state on $B$.
Since $G$ is a compact group, $\omega$ is a trace \cite[Theorem 4.1]{HLS}.
From Lemma~\ref{change of basis:lemma} one sees that (\ref{on:eq}) is equivalent to
(\ref{ON:eq}) with $e_{\gamma ki}$ replaced by $d^{1/2}_{\gamma}b_{\gamma ki}$.
Using $W(I\otimes 1_{\gamma_0})=w(I\otimes 1_{\gamma_0})=I\otimes 1_{\gamma_0}$ one gets
$U(I\otimes 1_{\gamma_0})=I\otimes 1_{\gamma_0}$. Thus $b_{\gamma_0 11}=1_{B}$.
Therefore $d^{1/2}_{\gamma}b_{\gamma ij}$ for $\gamma \in \hat{G}, 1\le i, j\le d_{\gamma}$
is a standard basis of $B$. Denote by $\psi(w)$ the associated element in $\cP$. Then the diagram
\begin{equation*} 
\xymatrix{
C^2 \ar[d]     \ar[rr]^{\psi} && \cP \ar[d] \\
C^2/H  \ar[r] & \EA(G)_{\fm} \ar@{^{(}->}[r]  & \EA(G)=\EA^{\sim}(G)
}
\end{equation*}
commutes, where we identify $\EA(G)$ with $\EA^{\sim}(G)$ since $G$ is co-amenable.
It was showed in the proof of \cite[Theorem 3.9]{Landstad} that one has
\begin{eqnarray*}
U_{12}U_{13}=(I\otimes w)((\id \otimes \delta)(U)) \quad \mbox{ and } \quad (\id\otimes \nu)(U)=U^*,
\end{eqnarray*}
where $U_{12}$ and $U_{13}$ are in the leg numbering notation.
 It follows that the map
$\psi$ is continuous. Consequently, the relative topology on $\EA(G)_{\fm}$ in $\EA(G)$ coincides with
the quotient topology coming from $C^2\rightarrow \EA(G)_{\fm}$.

\section{Induced Lip-norm}
\label{induced:sec}

In this section we prove Theorem~\ref{induced Lip-norm:thm}.

We recall first Rieffel's construction of Lip-norms from ergodic actions of compact groups.
Let $G$ be a compact group. A {\it length function} on $G$ is a continuous function
$\nl:G\rightarrow \Re_{+}$ such that
\begin{eqnarray*}
\mathnormal{l}(xy)&\le &\mathnormal{l}(x)+\mathnormal{l}(y)\mbox{ for
  all } x, y\in G \\
\mathnormal{l}(x^{-1})&=& \mathnormal{l}(x) \mbox{ for all }x\in G \\
\mathnormal{l}(x)&=& 0 \mbox{ if and only if } x=e_G.
\end{eqnarray*}
Given an ergodic action $\alpha$ of $G$ on a unital $C^*$-algebra $B$, Rieffel
showed that the seminorm
$L_{B}$ on $B$ defined by
\begin{eqnarray} \label{Rieffel induce:eq}
L_{B}(b)=\sup \{\frac{\pa \alpha_x(b)-b\pa}{\mathnormal{l}(x)}: x\in G, x\neq e_G\}
\end{eqnarray}
is a Lip-norm \cite[Theorem 2.3]{Rieffel98b}.

Note that there is a $1$-$1$ correspondence between length functions on $G$ and
left-invariant metrics on $G$ inducing the topology of $G$, via $\rho(x, y)=\mathnormal{l}(x^{-1}y)$
and $\nl(x)=\rho(x, e_G)$. Since a quantum metric on
(the non-commutative space corresponding to) a unital $C^*$-algebra is a Lip-norm
on this $C^*$-algebra,  a length function for a compact quantum group $A=C(\cG)$
should be a Lip-norm $L_{A}$ on $A$ satisfying certain compatibility condition with
the group structure. The proof of \cite[Proposition 2.2]{Rieffel98b} shows that
$L_{B}$ in above is finite on any $\alpha$-invariant finite-dimensional subspace of $B$, and
hence is finite on $\cB$.  If one applies this observation to the action of $G$ on $C(G)$ corresponding
to the right translation of $G$ on itself, then we see that the Lipschitz seminorm $L_{C(G)}$ on $C(G)$ associated to
the above metric $\rho$ via
\begin{eqnarray*}
L_{C(G)}(a)=\sup_{x\neq y}\frac{|a(x)-a(y)|}{\rho(x, y)}=\sup_{x\neq e_G}\sup_{y}\frac{|a(yx)-a(y)|}{\nl(x)}
\end{eqnarray*}
is finite on the algebra of regular functions in $C(G)$.
This leads to the following definition:

\begin{definition} \label{regular Lip-norm:def}
We say that a Lip-norm $L_{A}$ on a compact quantum group $A=C(\cG)$ is {\it regular} if $L_{A}$ is finite on
the algebra $\cA$ of regular functions.
\end{definition}

It turns out that a regular Lip-norm is sufficient for us to induce Lip-norms on $C^*$-algebras carrying
ergodic actions
of co-amenable compact quantum groups. We leave the discussion of the left and right invariance of
$L_{A}$ to the end of this section.

\begin{remark} \label{existence of regular:remark}
If a unital $C^*$-algebra $B$ has a Lip-norm, then $S(B)$ with the weak-$*$ topology is metrizable and hence
$B$ is separable. Conversely, if $B$ is a separable unital $C^*$-algebra, then for any countable
subset $W$ of $B$, there exist Lip-norms on $B$ being finite on $W$ \cite[Proposition 1.1]{Rieffel02}.
When $A=C(\cG)$ is separable, $\cA$ is a countable-dimensional vector space,
and hence $A$ has regular Lip-norms.
\end{remark}

\begin{example} \label{regular:example}
Let $\Gamma$ be a discrete group. Then the reduced group $C^*$-algebra $C^*_r(\Gamma)$ is a compact quantum group
with $\Phi(g)=g\otimes g$ for $g\in \Gamma$.
Its algebra of regular functions is $\Ce \Gamma$. Let $\nl$ be a length function
on $\Gamma$. Denote by $D$ the (possibly unbounded) linear operator of pointwise multiplication by $\nl$ on $\ell^2(\Gamma)$.
One may consider the seminorm $L$ defined on $\Ce \Gamma$ as $L(a)=\pa [D, a]\pa$ and
extend it to $C^*_r(\Gamma)$ via setting $L=\infty$ on $C^*_r(\Gamma)\setminus \Ce\Gamma$.
The seminorm $L$ so defined is always finite on $\Ce \Gamma$, and hence is regular if it is
a Lip-norm. This is the case for $\Gamma=\Ze^d$ when $\nl$ is a word-length, or the restriction
to $\Ze^d$ of a norm on $\Re^d$ \cite[Theorem 0.1]{Rieffel02}, and for $\Gamma$ being a hyperbolic group when
$\nl$ is a word-length \cite[Corollary 4.4]{Rieffel-Ozawa}.
\end{example}

Now we try to extend (\ref{Rieffel induce:eq}) to ergodic actions of compact quantum groups.
Let $\sigma:B\rightarrow B\otimes C(G)=C(G, B)$ be the $*$-homomorphism associated to $\alpha$, i.e.,
$(\sigma(b))(x)=\alpha_x(b)$ for $b\in B$ and $x\in G$.
For any $b\in B_{\sa}$, we have
\begin{eqnarray*}
L_{B}(b)&=&\sup_{x\neq e_G}\sup_{y} \frac{\pa \alpha_{yx}(b)-\alpha_y(b)\pa}{\nl(x)}\\
&=&\sup_{x\neq e_G}\sup_{y}\sup_{\varphi \in S(B)}\frac{|\varphi(\alpha_{yx}(b))-\varphi(\alpha_y(b))|}{\nl(x)}\\
&=&\sup_{\varphi\in S(B)} L_{C(G)}(b*\varphi),
\end{eqnarray*}
where $S(B)$ denotes the state space of $B$.
Note that for quantum metrics, only the restriction of $L_{B}$ on $B_{\sa}$ is essential. Thus
the above formula leads to our definition of the (possibly $+\infty$-valued)
seminorm $L_{B}$ on $B$ in (\ref{induced Lip-norm:eq})
for any ergodic action $\sigma:B\rightarrow B\otimes A$
of a compact quantum group $A=C(\cG)$ equipped with a regular Lip-norm $L_{A}$.

Throughout the rest of this section we assume that $L_{A}$ is a regular Lip-norm on $A$.
%

\begin{lemma} \label{Haar vs radius:lemma}
We have
\begin{eqnarray} \label{Harr vs radius:eq}
\pa a-h(a)1_{A}\pa \le 2r_{A}L_{A}(a)
\end{eqnarray}
for all $a\in A_{\sa}$.
\end{lemma}
\begin{proof} By Proposition~\ref{criterion of Lip:prop}
we can find $a'\in \Ce 1_{A}$ such that $\pa a-a'\pa \le r_{A}L_{A}(a)$.
Then $\pa h(a)1_{A}-a'\pa =|h(a-a')|\le \pa a-a'\pa\le r_{A}L_{A}(a)$. Thus
$\pa a-h(a)1_{A}\pa \le \pa a-a'\pa +\pa a'-h(a)1_{A}\pa \le 2r_{A}L_{A}(a)$.
\end{proof}

\begin{lemma} \label{radius:lemma}
Let
$L_{B}$ be the seminorm on  a unital $C^*$-algebra $B$ defined via (\ref{induced Lip-norm:eq})
for an action $\sigma:B\rightarrow B\otimes A$
of $\cG$ on $B$.
Assume that $A$ has bounded counit $e$.
Then for any $b\in B_{\sa}$ we have $\pa b-E(b)\pa \le 2r_{A}L_{B}(b)$,
where $E:B\rightarrow B^{\sigma}$ is the canonical conditional expectation.
\end{lemma}
\begin{proof} Let $\varphi \in S(B)$. Note that $h(b*\varphi)=\varphi(E(b))$.
We have
\begin{eqnarray*}
\pa b*\varphi-\varphi(E(b))1_{A}\pa &=&\pa b*\varphi-h(b*\varphi)1_{A}\pa \\
&\overset{(\ref{Harr vs radius:eq})}\le & 2r_{A}L_{A}(b*\varphi)
\overset{(\ref{induced Lip-norm:eq})}\le 2r_{A}L_{B}(b).
\end{eqnarray*}
Thus
\begin{eqnarray*}
\sup_{\varphi\in S(B)}\pa (b-E(b))*\varphi\pa =\sup_{\varphi\in S(B)}\pa b*\varphi-\varphi(E(b))1_{A}\pa \le 2r_{A}L_{B}(b).
\end{eqnarray*}
Therefore by Remark~\ref{injective:remark} we have
\begin{eqnarray*}
\pa b-E(b)\pa &=&\pa e*(b-E(b))\pa =\sup_{\varphi\in S(B)}|\varphi( e*(b-E(b)))| \\
&=& \sup_{\varphi\in S(B)}|e((b-E(b))*\varphi)|
\le \sup_{\varphi\in S(B)}\pa (b-E(b))*\varphi\pa \\
&\le &2r_{A}L_{B}(b)
\end{eqnarray*}
as desired.
\end{proof}

For any $J\subseteq \hat{\cG}$ denote $\sum_{\gamma \in J}A_{\gamma}$ and $\sum_{\gamma \in J}B_{\gamma}$
by $A_{J}$ and $B_{J}$ respectively.

\begin{lemma} \label{approx e:lemma}
Assume that $A$ has faithful Haar measure.
For any $\varepsilon>0$ and $\phi\in S(A)$ there exist $\psi\in S(A)$ and
a finite subset $J\subseteq \hat{\cG}$ such that $\psi$ vanishes on
$A_{\gamma}$ for all $\gamma \in \hat{\cG}\setminus J$ and
\begin{eqnarray} \label{approx e:eq1}
|(\phi-\psi)(a)|\le \varepsilon L_{A}(a)
\end{eqnarray}
for all $a\in A_{\sa}$.
\end{lemma}
\begin{proof}
Denote by $W$ the set of states of $A$ consisting of convex combinations of states of the form $h(a^* (\cdot) a)$ for
$a\in \cA$ with $h(a^*a)=1$. Let $\psi \in W$. Clearly there exists a finite subset
$F\subseteq \hat{\cG}$ such that if $h(A^*_{F}a'A_{F})=0$ for some $a'\in A$
then $\psi(a')=0$.
By the faithfulness of $h$ on $\cA$ and
the Peter-Weyl theory \cite[Theorems 4.2 and 5.7]{Wor87}, for any $a'\in \cA$ and any finite subset
$J'\subseteq \hat{\cG}$,
$h(A^*_{J'}a')\neq 0$ if and only if $h(a'A^*_{J'})\neq 0$.
Denote by $F'$ the set of equivalence classes of irreducible unitary subrepresentations of the tensor products
$u^{\alpha}\tot u^{\beta}$ of all $\alpha \in F$ and $\beta\in  F^c=\{\gamma^c:\gamma \in F\}$.
Denote $(F')^c$ by $J$.
Suppose that $\psi$ does not vanish on $A_{\gamma}$
for some $\gamma\in \hat{\cG}$. Then
$h(A^*_{F}A_{\gamma}A_{F})\neq 0$.
Thus
$$h(A_{\gamma}A_{F'})\overset{(\ref{*-algebra:eq})}\supseteq h(A_{\gamma}A_{F}A_{F^c})=h(A_{\gamma}A_{F}A^*_{F})\supsetneqq \{0\}.$$
Since $h(A_{\alpha}A_{\beta})=0$ for all $\alpha\neq \beta^c$ in $\hat{\cG}$ \cite[Theorem 5.7]{Wor87}, we get
$\gamma \in J$.

Now we just need to find $\psi\in W$ such that
(\ref{approx e:eq1}) holds for all $a\in A$.
Since $h$ is faithful, the {\GNS} representation $(\pi_{A}, H_{A})$ of $A$ associated to $h$ is faithful.
Thus convex combinations of vector states from $(\pi_{A}, H_{A})$ are
weak-$*$ dense in $S(A)$ \cite[Lemma T.5.9]{WO}. Note that $\cA$ is dense in $A$.
Therefore $W$ is weak-$*$ dense in $S(A)$.
Take $R\ge r_A$.
Since $\cD_R(A)$ is totally bounded by Proposition~\ref{criterion of Lip:prop},
we can find $\psi\in W$  such
that
\begin{eqnarray} \label{approx e:eq2}
|(\phi-\psi)(a)|\le \varepsilon
\end{eqnarray}
for all $a\in \cD_R(A)$.
By Proposition~\ref{criterion of Lip:prop} we have $\cE(A)=\cD_R(A)+\Re\cdot 1_A$.
Therefore (\ref{approx e:eq2}) holds for all $a\in \cE(A)$, from which
(\ref{approx e:eq1}) follows.
\end{proof}

The next lemma is an analogue of \cite[Lemmas 8.3 and 8.4]{Rieffel00} and \cite[Lemma 10.8]{Li10}.

\begin{lemma} \label{finite approx:lemma}
Let $B$ and $L_B$ be as in Lemma~\ref{radius:lemma}.
Assume that $A$ is co-amenable.
Let $\varepsilon>0$ and take $\psi$ and $J$ in Lemma~\ref{approx e:lemma}
for $\phi$ being the counit $e$.
Denote by $P_{\psi}$ the linear map $B\rightarrow B$ sending $b\in B$ to $\psi*b$.
Then $P_{\psi}(B)\subseteq B_{J}$ and
\begin{eqnarray} \label{finite approx:eq}
\pa P_{\psi}(b)\pa \le \pa b\pa,
\quad  \mbox{ and } \quad \pa b-P_{\psi}(b)\pa\le
\varepsilon L_{B}(b)
\end{eqnarray}
for all $b\in B_{\sa}$.
\end{lemma}
\begin{proof} Since $\psi$ vanishes on $A_{\gamma}$ for all $\gamma \in \hat{\cG}\setminus J$
and $\sigma(B_{\beta})\subseteq B_{\beta}\odot A_{\beta}$ for all $\beta\in \hat{\cG}$,
we have $P_{\psi}(B_{\beta})\subseteq B_{J}$ for all $\beta\in \hat{\cG}$.
Note that $B_{\beta}$ is finite dimensional and $\cB=\sum_{\beta \in \hat{\cG}}B_{\beta}$ is dense in
$B$. Thus $P_{\psi}(B)\subseteq B_{J}$.

For any $b\in B$ clearly $\pa P_{\psi}(b)\pa \le \pa b\pa$.
If $b\in B_{\sa}$, by Remark~\ref{injective:remark} we have
\begin{eqnarray*}
\pa b-P_{\psi}(b)\pa &=&\pa e*(b-P_{\psi}(b))\pa =\sup_{\varphi\in S(B)}| \varphi(e*(b-P_{\psi}(b)))| \\
                     &=& \sup_{\varphi\in S(B)} | e(b*\varphi)-\psi(b*\varphi)|
                        \overset{(\ref{approx e:eq1})}\le  \sup_{\varphi\in S(B)} \varepsilon L_{A}(b*\varphi) \\
                     &\overset{(\ref{induced Lip-norm:eq})}= & \varepsilon L_{B}(b).
\end{eqnarray*}
This finishes the proof of Lemma~\ref{finite approx:lemma}.
\end{proof}

We are ready to prove Theorem~\ref{induced Lip-norm:thm}.
\begin{proof}[Proof of Theorem~\ref{induced Lip-norm:thm}]
We verify the conditions in Proposition~\ref{criterion of Lip:prop}.
For any $b\in B$ and $\varphi$ in $S(B)$ we have
$b^**\varphi=(b*\varphi)^*$. Since $L_{A}$ satisfies the reality condition
(\ref{real:eq}), so does $L_{B}$.
For any $b\in \cB$, $\{b*\varphi:\varphi\in S(B)\}$ is bounded and contained in
a finite dimensional subspace of $\cA$ since $\sigma(\cB)\subseteq \cB\odot \cA$.
Then $L_{B}$ is finite on $\cB$ because of the regularity of $L_{A}$.
Clearly $L_{B}$ vanishes on $\Ce 1_{B}$.
By Lemma~\ref{radius:lemma} we have $\pa \cdot\pa^{\sim} \le 2r_{A}\tilde{L}_{B}$ on $(\tilde{B})_{\sa}$.
For any $\varepsilon>0$ let $P_{\psi}$ and $J$ be as in Lemma~\ref{finite approx:lemma}.
Then $P_{\psi}(\cD_1(B))$ is a bounded subset of the finite dimensional space
$B_{J}$. Thus $P_{\psi}(\cD_1(B))$ is totally bounded.
Since $\varepsilon>0$ is arbitrary, $\cD_1(B)$ is also totally bounded.
Therefore Theorem~\ref{induced Lip-norm:thm} follows from Proposition~\ref{criterion of Lip:prop}.
\end{proof}

Now we consider the invariance of a (possibly $+\infty$-valued) seminorm on $B$ with respect to
an action $\sigma$ of $\cG$. We consider first the case $\cG=G$ is a compact group.
For any action of $A=C(G)$ on $B$, there is a strongly continuous action $\alpha$ of $G$ on $B$
such that for any $b\in B$, the element $\sigma(b)\in B\otimes A=C(G, B)$ is given by $(\sigma(b))(x)=\alpha_x(b)$
for all $x\in G$.
If a seminorm $L_{B}$ on $B$ is lower semi-continuous, which is the case
if $L_{B}$ is defined via (\ref{Rieffel induce:eq}),
and is $\alpha$-invariant, then for any $\psi\in S(A)$ corresponding to a Borel probability measure
$\mu$ on $G$, we have
\begin{eqnarray*}
L_{B}(\psi*b)=L_{B}(\int_G \alpha_x(b)\, d\mu(x))\le L_{B}(b)
\end{eqnarray*}
for all $b\in B$. Conversely, if $L_{B}(\psi*b)\le L_{B}(b)$ for
all $b\in B$ and $\psi\in S(B)$, taking $\psi$ to be the evaluation at $x\in G$,
one sees immediately that $L_{B}$ is $\alpha$-invariant.
Note that the essential information about the quantum metric is the restriction of $L_B$ on $B_{\sa}$.
This leads to the following

\begin{definition} \label{invariant:def}
Let $A=C(\cG)$ be a compact quantum group.
We say that a (possibly $+\infty$-valued) seminorm  $L_{A}$ on $A$ is {\it right-invariant} ({\it left-invariant} resp.)
if
\begin{eqnarray*}
L_{A}(\psi*a)\le L_{A}(a)  \quad (L_{A}(a*\psi)\le L_{A}(a) \mbox{ resp. })
\end{eqnarray*}
for all $a\in A_{\sa}$ and $\psi\in S(A)$.
For an action $\sigma:B\rightarrow B\otimes A$ of  $\cG$ on a unital $C^*$-algebra $B$,
we say that a (possibly $+\infty$-valued) seminorm  $L_{B}$ on $B$
is {\it invariant}   if
\begin{eqnarray*}
L_{B}(\psi*b)\le L_{B}(b)
\end{eqnarray*}
for all $b\in B_{\sa}$ and $\psi\in S(A)$.
\end{definition}

\begin{proposition} \label{exist invariant:prop}
Let $L_{A}$ be a regular Lip-norm on $A$.
Define (possibly $+\infty$-valued) seminorms $L'_{A}$ and $L''_{A}$ on $A$ by
\begin{eqnarray*}
L'_{A}(a)=\sup_{\varphi\in S(A)} L_{A}(\varphi*a),
\end{eqnarray*}
and
\begin{eqnarray*}
L''_{A}(a)=\sup_{\varphi\in S(A)} L_{A}(a*\varphi)
\end{eqnarray*}
for $a\in A$.
Assume that $A$ has bounded counit.
Then $L'_{A}$ ($L''_{A}$ resp.) is a right-invariant (left-invariant resp.) regular
Lip-norm on $A$, and
$L'_{A}\ge L_{A}$ ($L''_{A}\ge L_{A}$ resp.).
If $L_{A}$ is left-invariant (right-invariant resp.), then so is $L'_{A}$ ($L''_{A}$ resp.).
\end{proposition}
\begin{proof}
An argument similar to that in the proof of Theorem~\ref{induced Lip-norm:thm}
shows that $L'_{A}$ satisfies the reality condition (\ref{real:eq}), vanishes on
$\Ce 1_{A}$, and is finite on $\cA$.
Taking $\varphi$ to be the counit we see that $L'_{A}\ge L_{A}$.
It follows immediately from  Proposition~\ref{criterion of Lip:prop}
that $L'_{A}$ is a regular Lip-norm on $A$.
For any $a\in A_{\sa}$ and $\psi\in S(A)$ we have
\begin{eqnarray*}
L'_{A}(\psi*a)=\sup_{\varphi\in S(A)}L_{A}(\varphi*(\psi*a))
=\sup_{\varphi\in S(A)}L_{A}((\varphi*\psi)*a)
\le L'_{A}(a),
\end{eqnarray*}
where $\varphi*\psi$ is the state on $A$ defined via $(\varphi*\psi)(a')=(\varphi \otimes \psi)(\Phi(a'))$
for $a'\in A$.
Therefore $L'_{A}$ is right-invariant.
Assume that $L_{A}$ is left-invariant. Then for any $a\in A_{\sa}$ and $\psi\in S(A)$ we have
\begin{eqnarray*}
L'_{A}(a*\psi)&=&\sup_{\varphi\in S(A)}L_{A}(\varphi*(a*\psi))
=\sup_{\varphi\in S(A)}L_{A}((\varphi*a)*\psi)\\
&\le & \sup_{\varphi\in S(A)}L_{A}(\varphi*a)
= L'_{A}(a).
\end{eqnarray*}
Thus $L'_{A}$ is also left-invariant.
The assertions about $L''_{A}$ are proved similarly.
\end{proof}

Using Remark~\ref{existence of regular:remark} and
applying the construction in Proposition~\ref{exist invariant:prop} twice,
we get

\begin{corollary} \label{biinvariant:coro}
Every separable compact quantum group with bounded counit
has a bi-invariant regular Lip-norm.
\end{corollary}

An argument similar to that in the proof of Proposition~\ref{exist invariant:prop} shows

\begin{proposition} \label{invariant:prop}
Let $\sigma$ be an action of $\cG$ on a unital $C^*$-algebra $B$.
If $L_{A}$ is a right-invariant regular Lip-norm on $A$, then $L_{B}$ defined via (\ref{induced Lip-norm:eq})
is invariant.
\end{proposition}

\section{Quantum distance}
\label{distance:sec}

In this section we introduce the quantum distance $\dist_{\rre}$
between ergodic
actions of $\cG$,
and prove Theorem~\ref{distance vs topology:thm}.

Throughout this section,
$A$ will be a co-amenable compact quantum group with a fixed regular Lip-norm $L_{A}$.
For any ergodic action $(B, \sigma)$ of $\cG$, we endow $B$ with the Lip-norm
$L_{B}$ in Theorem~\ref{induced Lip-norm:thm}.

In \cite{Kerr, KL, Li10, Li12, Rieffel00} several quantum
Gromov-Hausdorff distances are introduced, applying to quantum
metric spaces in various contexts as order-unit spaces, operator
systems, and $C^*$-algebras. They are all applicable to
$C^*$-algebraic compact quantum metric spaces, which we are
dealing with now. Among these distances, the unital version
$\dist_{\rnu}$ of the one introduced in \cite[Remark 5.5]{Li12} is
the strongest one, which we recall below from \cite[Section 5]{KL}. To simplify the
notation, for fixed unital $C^*$-algebras $B_1$ and $B_2$,
when we take infimum over unital $C^*$-algebras $C$ containing
both $B_1$ and $B_2$, we mean to take infimum over all unital
injective $*$-homomorphisms of $B_1$ and $B_2$ into some
unital $C^*$-algebra $C$. We denote by $\dist^{C}_{\rH}$ the
Hausdorff distance between subsets of $C$.
Recall that $\cE(B):=\{b\in B_{\sa}:L_B(b)\le 1\}$.
For any $C^*$-algebraic
compact quantum metric spaces $(B_1, L_{B_1})$ and $(B_2,
L_{B_2})$, the distance $\dist_{\rnu}(B_1, B_2)$
is  defined as
\begin{eqnarray*}
\dist_{\rnu}(B_1, B_2)&=&\inf \dist^C_{\rH}(\cE(B_1), \cE(B_2)),
\end{eqnarray*}
where the infimum is taken over all unital $C^*$-algebras $C$ containing $B_1$ and $B_2$.
Note that $\dist_{\rnu}(B_1, B_2)$ is always finite since
$\cD_R(\cB)$ is totally bounded and $\cE(B)=\cD_R(\cB)+\Re\cdot 1_B$ for any $R\ge r_B$ by
Proposition~\ref{criterion of Lip:prop}.
These distances become zero whenever there is a $*$-isomorphism $\varphi:B_1\rightarrow B_2$
preserving the Lip-norms on the self-adjoint parts. In particular, as the following example shows,
these distances may not distinguish
the actions when the Lip-norms $L_{B_i}$ come from ergodic actions of $\cG$ on $B_i$.

\begin{example}   \label{tori:ex}
Let $\nl'$ be a length function on the circle $S^1$.
Set $\nl$ to be the length function on the two-torus $\Te^2$
defined as $\nl(x, y)=\nl'(x)+\nl'(y)$ for $x, y\in S^1$.
Then $\nl(x, y)=\nl(x^{-1}, y)$ for all $(x, y)\in \Te^2$.
Let $\theta\in\Re$, and
let $B_{\theta}$ be the non-commutative two-torus generated by
unitaries $u_{\theta}$ and $v_{\theta}$ satisfying $u_{\theta}v_{\theta}=e^{2\pi i \theta}v_{\theta}u_{\theta}$.
Then $\Te^2$ has a strongly continuous action $\alpha_{\theta}$
on $B_{\theta}$ specified by  $\alpha_{\theta, (x, y)}(u_{\theta})=xu_{\theta}$ and $\alpha_{\theta, (x, y)}(v_{\theta})=yv_{\theta}$.
Consider the $*$-isomorphism $\psi:B_{\theta}\rightarrow B_{-\theta}$ determined
by $\psi(u_{\theta})=(u_{-\theta})^{-1}$ and $\psi(v_{\theta})=v_{-\theta}$. Then $\psi$ preserves the Lip-norms
defined via (\ref{Rieffel induce:eq})
for the actions $\alpha_{\theta}$ and $\alpha_{-\theta}$ of $\Te^2$,  and
hence $B_{\theta}$ and $B_{-\theta}$ have distances zero under all the quantum distances
defined in \cite{Kerr, KL, Li10, Li12, Rieffel00}.  However, when $0<\theta<1/2$, the actions
$(B_{\theta}, \alpha_{\theta})$ and
$(B_{-\theta}, \alpha_{-\theta})$ are not
isomorphic, as can be seen from the fact that $\Ce u_{\theta}=\{b\in B_{\theta}: \alpha_{\theta, (x, y)}(b)=xb \mbox{ for all } (x, y)\in \Te^2\}$
and $\Ce v_{\theta}=\{b\in B_{\theta}: \alpha_{\theta, (x, y)}(b)=yb \mbox{ for all } (x, y)\in \Te^2\}$.
\end{example}

\begin{notation} \label{graph:notation}
For any $C^*$-algebra $C$ we denote $C\oplus (C\otimes A)$ by $C^{\sharp}$.
For any action $\sigma:B\rightarrow B\otimes A$ of $\cG$ on a unital $C^*$-algebra $B$ and
any subset $\cX$ of $B$ we denote by $\cX_{\sigma}$ the graph
\begin{eqnarray*}
\{(b, \sigma(b))\in B^{\sharp}: b\in \cX\}
\end{eqnarray*}
of $\sigma|_{\cX}$.
\end{notation}


We are going to introduce a quantum distance between ergodic actions of $\cG$
to distinguish the actions.
Modifying the above definition of $\dist_{\rnu}$,  we just need
to add one term to take care of the actions:

\begin{definition} \label{distance:def}
Let $(B_1, \sigma_{1})$ and $(B_2, \sigma_{2})$ be ergodic actions of
$A$. We set
\begin{eqnarray*}
\dist_{\rre}(B_1, B_2)=\inf \dist^{C^{\sharp}}_{\rH}((\cE(B_1))_{\sigma_1},
(\cE(B_2))_{\sigma_2}),
\end{eqnarray*}
where the infimum is taken over all unital $C^*$-algebras $C$ containing
both $B_1$ and $B_2$.
\end{definition}

Clearly $\dist_{\rre}\ge \dist_{\rnu}$.
An argument similar to that in the proof of \cite[Theorem 3.15]{Li12} yields

\begin{proposition} \label{metric:prop}
The distance $\dist_{\rre}$
is a metric on $\EA(\cG)$.
\end{proposition}

We relate first continuous fields of ergodic actions of $\cG$ to the distance $\dist_{\rre}$.

\begin{proposition} \label{criterion:prop}
Suppose that $L_{A}$ is left-invariant.
Let $(\{(B_t, \sigma_t)\}_{t\in T}, B)$ be a continuous field of ergodic actions of $\cG$ over
a compact metric space $T$. Let $t_0\in T$.
If $\lim_{t\rightarrow t_0} \mul(B_t, \gamma)= \mul(B_{t_0}, \gamma)$
for all $\gamma \in \hat{\cG}$, then $\dist_{\rre}(B_t, B_{t_0})\rightarrow 0$ as $t\rightarrow t_0$.
\end{proposition}

To simplify the notation, we shall write $L_t$ for $L_{B_t}$ below.

\begin{lemma} \label{cont of Lip:lemma}
Let the notation be as in Proposition~\ref{criterion:prop}.
Let $J$ be a finite subset of $\hat{\cG}$, and let $b\in B$ such that $\pi_t(b)\in (B_t)_{J}$
for each $t\in T$. Then the function $t\mapsto L_{t}(\pi_t(b))$ is continuous on $T$.
\end{lemma}
\begin{proof}
Let $s\in T$.  To prove the continuity of $t\mapsto
L_{t}(\pi_t(b))$ at $t=s$, it suffices to show that for any
sequence $t_n\rightarrow s$ one has
$L_{t_n}(\pi_{t_n}(b))\rightarrow L_{s}(\pi_{s}(b))$.
By Remark~\ref{existence of regular:remark} each $B_t$ is separable.
Taking
restriction to the closure of this sequence, we may assume that
$B$ is separable. Since $A$ is co-amenable, any unital
$C^*$-algebra admitting an ergodic action of $A$ is nuclear
\cite{DLRZ}. Every separable continuous field of unital nuclear
$C^*$-algebras over a compact metric space can be subtrivialized
\cite[Theorem 3.2]{Blanchard}. Thus we can find a unital
$C^*$-algebra $C$ and unital embeddings $B_t\rightarrow C$
for all $t\in T$ such that (via identifying each $B_t$ with its
image in $C$) elements in $B$ are exactly those continuous
maps $T\rightarrow C$ whose images at each $t$ are in $B_t$.

Let $\varphi_{s}\in S(B_{s})$. Extend it to a state of $C$ and let $\varphi_t$ be the restriction
on $B_t$ for each $t\in T$. Then $\varphi_t\in S(B_t)$ for each $t\in T$ and
$\varphi_t(\pi_t(c))\rightarrow \varphi_{s}(\pi_{s}(c))$ as $t\rightarrow s$ for any $c\in B$.
Say,
\begin{eqnarray*}
\sigma_t(\pi_t(b))= \sum_{\gamma \in J}\sum_{ 1\le i, j\le d_{\gamma}}c_{\gamma ij}(t)\otimes u^{\gamma}_{ij}
\end{eqnarray*}
for all $t\in T$. Then clearly the sections $t\mapsto c_{\gamma ij}(t)$ are in $B$.
Thus $\pi_t(b)*\varphi_t$ converges to $\pi_{s}(b)*\varphi_{s}$ in $A_{J}$ as $t\rightarrow s$.
Since $A_{J}$ is finite dimensional, $L_{A}$ is continuous on $A_{J}$. Therefore
$L_{A}(\pi_t(b)*\varphi_t)$ converges to $L_{A}(\pi_{s}(b)*\varphi_{s})$ as $t\rightarrow s$.
Then it follows easily that the function $t\mapsto L_{t}(\pi_t(b))$ is lower semi-continuous at $s$.

Let $\varepsilon>0$. Take a sequence $t_1, t_2, \cdots$ in $T$ converging to $s$ such that
\begin{eqnarray*}
\varepsilon+L_{t_n}(\pi_{t_n}(b))\ge \limsup_{t\rightarrow s}L_t(\pi_t(b))
\end{eqnarray*}
for each $n\ge 1$. Take $\varphi_{t_n}\in S(B_{t_n})$ for each $n\ge 1$ such that
\begin{eqnarray*}
\varepsilon +L_{A}(\pi_{t_n}(b)*\varphi_{t_n})\ge L_{t_n}(\pi_{t_n}(b)).
\end{eqnarray*}
Since $B$ is separable, passing to a subsequence if necessary, we may assume that
$\varphi_{t_n}\circ \pi_{t_n}$ converges to some state $\psi$ of $B$ (in the weak-$*$ topology)
as $n\rightarrow \infty$.
Then $\psi= \varphi_{s}\circ \pi_{s}$ for some $\varphi_{s}\in S(B_{s})$ by the upper semi-continuity of 
the field $(\{B_t\}_{t\in T}, B)$.
We have $\varphi_{t_n}(\pi_{t_n}(c))\rightarrow \varphi_{s}(\pi_{s}(c))$ as $n\rightarrow \infty$ for any $c\in B$. As
in the second paragraph of the proof, $L_{A}(\pi_{t_n}(b)*\varphi_{t_n})$ converges to $L_{A}(\pi_{s}(b)*\varphi_{s})$
as $n\rightarrow \infty$. Therefore,
\begin{eqnarray*}
2\varepsilon+L_{s}(\pi_{s}(b))\ge 2\varepsilon+L_{A}(\pi_{s}(b)*\varphi_{s})\ge \limsup_{t\rightarrow t_0} L_t(\pi_t(b)).
\end{eqnarray*}
Thus  the function $t\mapsto L_{t}(\pi_t(b))$ is upper semi-continuous at $s$ and hence continuous at $s$.
\end{proof}

\begin{lemma} \label{uniform dense:lemma}
Let $V$ be a finite-dimensional vector space, and let $W$ be a linear subspace of $V$.
Let $T$ be a topological space. Let $\pa \cdot \pa_t$ be a norm on $V$ and
$L_t$ be a seminorm on $V$ vanishing exactly on $W$ for each $t\in T$ such that the functions
$t\mapsto \pa v\pa_t $ and $t\mapsto L_t(v)$ are upper semicontinuous and
continuous respectively on $T$ for every $v\in V$.
Let $t_0\in T$, and let $\varepsilon>0$. Then
\begin{eqnarray} \label{uniform dense:eq}
\dist^{\pa \cdot \pa_t}_{\rH}(\cE_{t_0}(V), \cE_t(V))\le \varepsilon
\end{eqnarray}
throughout some neighborhood $U$ of $t_0$, where $\cE_t(V)=\{v\in V:L_t(v)\le 1\}$.
\end{lemma}
\begin{proof}
Via considering $V/W$ we may assume that $W=\{0\}$.
For any $\delta >0$, a standard compactness argument shows that
\begin{eqnarray*}
  &\pa \cdot \pa_t& \le (1+\delta)\pa \cdot\pa_{t_0},\\
\frac{1}{1+\delta}L_{t_0}\le &L_t& \le (1+\delta)L_{t_0}
\end{eqnarray*}
throughout some neighborhood $U_{\delta}$ of $t_0$.
Then we can find some $R>0$ such that $\pa \cdot \pa_t\le R L_t(\cdot)$ throughout $U_1$.
Fix $\delta=R/\varepsilon$.
Let $t\in U_1\cap U_{\delta}$ and $v\in \cE_{t_0}(V)$. Then $v/(1+\delta)\in \cE_t(V)$, and
\begin{eqnarray*}
\pa v-v/(1+\delta)\pa_t=\frac{\delta}{1+\delta}\pa v\pa_t\le \delta \pa v\pa_{t_0}\le \delta R=\varepsilon.
\end{eqnarray*}
Similarly, for any
$t\in U_1\cap U_{\delta}$ and $v\in \cE_t(V)$, we have
$v/(1+\delta)\in \cE_{t_0}(V)$ and $\pa v-v/(1+\delta)\pa_t \le \varepsilon$.
This proves (\ref{uniform dense:eq}).
\end{proof}

We are ready to prove Proposition~\ref{criterion:prop}.

\begin{proof}[Proof of Proposition~\ref{criterion:prop}]
As in the first paragraph of the proof of Lemma~\ref{cont of Lip:lemma} we
may assume that there is a unital $C^*$-algebra $C$ containing each $B_t$
as a unital $C^*$-subalgebra and that elements in $B$ are exactly those continuous maps $T\rightarrow C$
whose images at each $t$ are in $B_t$.
Let $\varepsilon>0$.  Pick $\psi\in S(A)$ and $J\subseteq \hat{\cG}$
in Lemma~\ref{approx e:lemma} for $\phi$ being the counit.
We may assume that $\gamma_0\in J$ and $\gamma^c\in J$ for each $\gamma\in J$.
Then $1_{B_t}\in (B_t)_{J}$ and $((B_t)_{J})^*=(B_t)_{J}$.
By Proposition~\ref{invariant:prop} $L_t$ is invariant
for all $t\in T$.
By Lemma~\ref{finite approx:lemma}
we have
\begin{eqnarray} \label{criterion:eq1}
\dist^{C}_{\rH}(\cE(B_t),\, \cE((B_t)_{J}))\le
\varepsilon
\end{eqnarray}
for all $t\in T$, where $\cE((B_t)_{J}):=\cE(B_t)\cap
(B_t)_{J}$. Suppose that $\lim_{t\rightarrow t_0} \mul(B_t,
\gamma)= \mul(B_{t_0}, \gamma)$ for all $\gamma \in \hat{\cG}$.
By Lemma~\ref{lifting ON:lemma} there are a neighborhood $U$ of
$t_0$ and a linear isomorphism
$\varphi_t:(B_{t_0})_{J}\rightarrow (B_t)_{J}$ for each
$t\in U$ such that $\varphi_{t_0}=\id$, $\varphi_t((B_{t_0})_{\gamma})=(B_t)_{\gamma}$ for
each $\gamma \in J$ and $t\in U$, and the map $t\mapsto
\varphi_t(v)\in C$ is continuous over $U$ for all $v\in
(B_{t_0})_{J}$. Replacing $\varphi_t$ by $(\varphi_t+\varphi_t^*)(\varphi_t(1_{B_{t_0}})+\varphi_t(1_{B_{t_0}})^*)^{-1}$ and
shrinking $U$ if necessary,
we may assume
that $\varphi_t$ is unital
and Hermitian throughout $U$. By Lemma~\ref{cont of Lip:lemma} we know that
$\{\pa \cdot \pa_{C} \circ \varphi_t\}_{t\in U}$ and $\{L_t\circ
\varphi_t\}_{t\in U}$ are continuous families of norms and
seminorms on $(B_{t_0})_{J}$. By Lemma~\ref{uniform
dense:lemma}, shrinking $U$ if necessary, we have
\begin{eqnarray} \label{criterion:eq2}
\dist^{C}_{\rH}(\varphi_t(\cX),\, \cE((B_t)_{J}))<\varepsilon
\end{eqnarray}
throughout $U$, where $\cX=\cE((B_{t_0})_{J})$.
Putting (\ref{criterion:eq1}) and (\ref{criterion:eq2}) together, we get
\begin{eqnarray} \label{criterion:eq4}
\dist^{C}_{\rH}(\cE(B_t),\, \varphi_t(\cX))< 2\varepsilon
\end{eqnarray}
throughout $U$.
Note that
\begin{eqnarray*}
\dist^{C}_{\rH}(\cY,\, \cZ)=\dist^{C^{\sharp}}_{\rH}(\cY_{\sigma_t},\, \cZ_{\sigma_t})
\end{eqnarray*}
for any subsets $\cY, \cZ$ of $B_t$.
Thus
\begin{eqnarray} \label{criterion:eq5}
\dist^{C^{\sharp}}_{\rH}((\cE(B_t))_{\sigma_t}, \, (\varphi_t(\cX))_{\sigma_t})< 2\varepsilon
\end{eqnarray}
throughout $U$. By Lemma~\ref{tensor product:lemma} we may identify elements of $B\otimes A$ 
with the continuous maps $T\rightarrow C\otimes A$
whose images at each $t$ are in $B_t\otimes A$.
Since $\cD_R((B_{t_0})_{J})$ is totally bounded and $\cX=\cD_R((B_{t_0})_{J})+\Re\cdot 1_{B_{t_0}}$ for
any $R\ge 2r_A$ by Lemma~\ref{radius:lemma},
shrinking $U$ if necessary, we may assume that $\pa \sigma_t(\varphi_t(x))-\sigma_{t_0}(x)\pa_{C\otimes A}$,
$\pa \varphi_t(x)-x\pa_{C}<\varepsilon$ for
all $x\in \cX$ and $t\in U$.
Then
\begin{eqnarray} \label{criterion:eq6}
\dist^{C^{\sharp}}_{\rH}((\varphi_t(\cX))_{\sigma_t}, \, \cX_{\sigma_{t_0}})<\varepsilon
\end{eqnarray}
throughout $U$.  Putting (\ref{criterion:eq5}) and (\ref{criterion:eq6})
together, we get
\begin{eqnarray*}
\dist_{\rre}(B_t, B_{t_0})&\le &\dist^{C^{\sharp}}_{\rH}((\cE(B_t))_{\sigma_t}, \, (\cE(B_{t_0}))_{\sigma_{t_0}})
< 6\varepsilon
\end{eqnarray*}
throughout $U$. This finishes the proof of Proposition~\ref{criterion:prop}.
\end{proof}

\begin{remark} \label{cont of nu:remark}
Since $\dist_{\rre}\ge \dist_{\rnu}$ and
$\dist_{\rnu}$ is the strongest one among the quantum distances defined in \cite{Kerr, KL, Li10, Li12, Rieffel00},
Proposition~\ref{criterion:prop} also holds with $\dist_{\rre}$ replaced by any of them.
\end{remark}

We are ready to prove Theorem~\ref{distance vs topology:thm}.

\begin{proof}
[Proof of Theorem~\ref{distance vs topology:thm}]
By Proposition~\ref{metric:prop} $\dist_{\rre}$ is a metric on
$\EA(\cG)$.
By Theorem~\ref{separable:thm} and
Proposition~\ref{criterion:prop} the topology on $\EA(\cG)$ defined in
Definition~\ref{topology on EA:def} is stronger than that induced
by $\dist_{\rre}$. By Theorem~\ref{topology on EA:thm} the former is
compact. Thus these two topologies coincide.
\end{proof}




\begin{thebibliography}{999}

\bibitem{AP}
J.  Anderson and W. Paschke. The rotation algebra.
{\it Houston J. Math.}  {\bf 15}  (1989),  no. 1, 1--26.

\bibitem{Banica}
T. Banica. Representations of compact quantum groups and subfactors.
{\it J. Reine Angew. Math.}  {\bf 509}  (1999), 167--198. math.QA/9804015.

\bibitem{BMT}
E. B\'edos, G. J. Murphy, and L. Tuset. Co-amenability of compact quantum groups.
{\it J. Geom. Phys.}  {\bf 40}  (2001),  no. 2, 130--153. math.OA/0010248.

\bibitem{Blanchard}
\'E. Blanchard. Subtriviality of continuous fields of nuclear $C^*$-algebras.
{\it J. Reine Angew. Math.}  {\bf 489}  (1997), 133--149. math.OA/0012128.

\bibitem{Boca}
F. P. Boca. Ergodic actions of compact matrix pseudogroups on $C^*$-algebras.
In: {\it Recent Advances in Operator Algebras (Orl\'eans, 1992).  Ast\'erisque}  No. 232 (1995), 93--109.

\bibitem{Cuntz}
J. Cuntz.  Simple $C^*$-algebras generated by isometries.
{\it Comm. Math. Phys.}  {\bf 57}  (1977), no. 2, 173--185.

\bibitem{Dixmier}
J. Dixmier. {\it $C^*$-algebras.}
Translated from the French by Francis Jellett.
North-Holland Mathematical Library, Vol. 15. North-Holland Publishing Co., Amsterdam-New York-Oxford, 1977.

\bibitem{DLPS}
L. D\c abrowski, G. Landi, M. Paschke, and A. Sitarz. The spectral geometry of the equatorial Podle\'s sphere.
{\it C. R. Math. Acad. Sci. Paris}  {\bf 340}  (2005),  no. 11, 819--822.
math.QA/0408034.

\bibitem{DS}
L.  D\c abrowski and A. Sitarz. Dirac operator on the standard Podle\'s quantum sphere.
In: {\it Noncommutative Geometry and Quantum Groups (Warsaw, 2001)},  49--58, Banach Center Publ., 61,
Polish Acad. Sci., Warsaw, 2003.

\bibitem{DLRZ}
S. Doplicher, R. Longo, J. E. Roberts, and  L. Zsid\'o. A remark on quantum group actions and nuclearity.
Dedicated to Professor Huzihiro Araki on the occasion of his 70th birthday.
{\it Rev. Math. Phys.}  {\bf 14}  (2002),  no. 7-8, 787--796. math.OA/0204029.

\bibitem{DG}
M. Dupr\'e and R. M. Gillette. {\it Banach Bundles, Banach Modules and Automorphisms of $C^*$-algebras.}
Research Notes in Mathematics, 92.
Pitman (Advanced Publishing Program), Boston, MA, 1983.




\bibitem{EN}
R. Exel and C.-K. Ng. Approximation property of $C^*$-algebraic bundles.
{\it Math. Proc. Cambridge Philos. Soc.} {\bf 132} (2002), no. 3, 509--522. math.OA/9906070.



\bibitem{HMS}
P. M. Hajac, R. Matthes, and W. Szyma\'nski. Chern numbers for two families of noncommutative Hopf fibrations.
{\it C. R. Math. Acad. Sci. Paris}  {\bf 336}  (2003),  no. 11, 925--930. math.QA/0302256.

\bibitem{HLS}
R.  H\o egh-Krohn, M. R. Landstad, and E. St\o rmer. Compact ergodic groups of automorphisms.
{\it Ann. of Math. (2)}  {\bf 114}  (1981), no. 1, 75--86.

\bibitem{Kerr}
D. Kerr. Matricial quantum Gromov-Hausdorff distance.
{\it J. Funct. Anal.}  {\bf 205}  (2003),  no. 1, 132--167. math.OA/0207282.

\bibitem{KL}
D. Kerr and H. Li. On Gromov-Hausdorff convergence for operator
metric spaces.  {\it J. Operator Theory} to appear. math.OA/0411157.

\bibitem{KW}
E. Kirchberg and S. Wassermann. Operations on continuous bundles of $C^*$-algebras.
{\it Math. Ann.}  {\bf 303}  (1995),  no. 4, 677--697.

\bibitem{KNW}
Y. Konishi, M. Nagisa, and Y. Watatani. Some remarks on actions of compact matrix quantum groups on $C^*$-algebras.
{\it Pacific J. Math.} {\bf 153} (1992), no. 1, 119--127.

\bibitem{Lance}
E. C. Lance. {\it Hilbert $C\sp *$-modules. A Toolkit for Operator Algebraists.}
London Mathematical Society Lecture Note Series, 210. Cambridge University Press, Cambridge, 1995.

\bibitem{Landstad}
M. B. Landstad. Ergodic actions of nonabelian compact groups.
In: {\it Ideas and Methods in Mathematical Analysis, Stochastics, and Applications (Oslo, 1988)},  365--388,
Cambridge Univ. Press, Cambridge, 1992.

\bibitem{Li10}
H. Li. Order-unit quantum Gromov-Hausdorff distance.
{\it J. Funct. Anal.} {\bf 231} (2006), no. 2, 312--360.
math.OA/0312001.

\bibitem{Li12}
H. Li. $C^*$-algebraic quantum Gromov-Hausdorff
distance.
math.OA/0312003 v3.

\bibitem{MV}
A. Maes and A. Van Daele. Notes on compact quantum groups.
{\it Nieuw Arch. Wisk. (4)}  {\bf 16}  (1998),  no. 1-2, 73--112.

\bibitem{Marciniak}
M. Marciniak. Actions of compact quantum groups on $C^*$-algebras.
{\it Proc. Amer. Math. Soc.}  {\bf 126}  (1998),  no. 2, 607--616.

\bibitem{Nagy}
G. Nagy. On the Haar measure of the quantum $\SU(N)$ group.
{\it Comm. Math. Phys.}  {\bf 153}  (1993),  no. 2, 217--228.

\bibitem{OPT}
D. Olesen, G. K. Pedersen, and M. Takesaki. Ergodic actions of compact abelian groups.
{\it J. Operator Theory}  {\bf 3}  (1980), no. 2, 237--269.

\bibitem{Rieffel-Ozawa}
N. Ozawa and M. A. Rieffel. Hyperbolic group $C^*$-algebras and free-product $C^*$-algebras as compact quantum metric spaces.
{\it Canad. J. Math.} {\bf 57} (2005), no. 5, 1056--1079. math.OA/0302310.

\bibitem{Paolucci}
A. Paolucci. Coactions of Hopf algebras on Cuntz algebras and their fixed point algebras.
{\it Proc. Amer. Math. Soc.}  {\bf 125}  (1997),  no. 4, 1033--1042.

\bibitem{Podles87}
P. Podle\'s. Quantum spheres.
{\it Lett. Math. Phys.}  {\bf 14}  (1987),  no. 3, 193--202.

\bibitem{Podles95}
P. Podle\'s. Symmetries of quantum spaces. Subgroups and quotient spaces of quantum $\SU(2)$ and ${\rm SO}(3)$ groups.
{\it Comm. Math. Phys.}  {\bf 170}  (1995),  no. 1, 1--20. hep-th/9402069.

\bibitem{PW}
P. Podle\'s and S. L. Woronowicz. Quantum deformation of Lorentz group.
{\it Comm. Math. Phys.}  {\bf 130}  (1990),  no. 2, 381--431.



\bibitem{Rieffel89}
M. A. Rieffel. Continuous fields of $C^*$-algebras coming from group cocycles and actions.
{\it Math. Ann.}  {\bf 283}  (1989),  no. 4, 631--643.

\bibitem{Rieffel98b}
M. A. Rieffel. Metrics on states from actions of
compact groups.
{\it Doc. Math.} {\bf 3}  (1998), 215--229 (electronic).
math.OA/9807084.


\bibitem{Rieffel00}
M. A. Rieffel. Gromov-Hausdorff distance for
quantum  metric spaces.
{\it Mem. Amer. Math. Soc.}  {\bf 168}  (2004),  no. 796, 1--65.
 math.OA/0011063.


\bibitem{Rieffel01}
M. A. Rieffel. Matrix algebras converge to the
sphere for  quantum Gromov-Hausdorff distance.
{\it Mem. Amer. Math. Soc.}  {\bf 168}  (2004),  no. 796, 67--91.
math.OA/0108005.

\bibitem{Rieffel02}
M. A. Rieffel. Group $C^*$-algebras as compact quantum metric spaces.
{\it Doc. Math.}  {\bf 7}  (2002), 605--651 (electronic). math.OA/0205195.

\bibitem{Rieffel03}
M. A. Rieffel. Compact quantum metric spaces.
In:  {\it Operator Algebras, Quantization, and Noncommutative Geometry},  315--330,
Contemp. Math., 365, Amer. Math. Soc., Providence, RI, 2004.
math.OA/0308207.

\bibitem{VS}
L. L. Vaksman and Ya. S. Soibelman. An algebra of functions on the quantum group $\SU(2)$. (Russian)
{\it Funktsional. Anal. i Prilozhen.}  {\bf 22}  (1988),  no. 3, 1--14, 96;
translation in  {\it Funct. Anal. Appl.}  {\bf 22}  (1988),  no. 3, 170--181 (1989).

\bibitem{Wang95}
S. Wang. Tensor products and crossed products of compact quantum groups.
{\it Proc. London Math. Soc. (3)}  {\bf 71}  (1995),  no. 3, 695--720.

\bibitem{Wang98}
S. Wang. Quantum symmetry groups of finite spaces.
{\it Comm. Math. Phys.}  {\bf 195}  (1998),  no. 1, 195--211.

\bibitem{Wang99}
S. Wang. Ergodic actions of universal quantum groups on operator algebras.
{\it Comm. Math. Phys.}  {\bf 203}  (1999),  no. 2, 481--498. math.OA/9807093.

\bibitem{Wassermann89}
A. Wassermann. Ergodic actions of compact groups on operator algebras. I. General theory.
{\it Ann. of Math. (2)}  {\bf 130}  (1989),  no. 2, 273--319.

\bibitem{Wassermann}
A. Wassermann. Ergodic actions of compact groups on operator algebras. II. Classification of full multiplicity ergodic actions.
{\it Canad. J. Math.}  {\bf 40}  (1988),  no. 6, 1482--1527.

\bibitem{Wassermann88}
A. Wassermann, Ergodic actions of compact groups on operator algebras. III. Classification for $\SU(2)$.
{\it Invent. Math.}  {\bf 93}  (1988),  no. 2, 309--354.

\bibitem{WO}
N. E. Wegge-Olsen.
{\it $K$-theory and $C\sp *$-algebras. A Friendly Approach.}
Oxford Science Publications. The Clarendon Press, Oxford University Press, New York, 1993.

\bibitem{WorSU2}
S. L. Woronowicz. Twisted $\SU(2)$ group. An example of a noncommutative differential calculus.
{\it Publ. Res. Inst. Math. Sci.}  {\bf 23}  (1987),  no. 1, 117--181.

\bibitem{Wor87}
S. L. Woronowicz. Compact matrix pseudogroups.
{\it Comm. Math. Phys.}  {\bf 111}  (1987),  no. 4, 613--665.

\bibitem{Wor95}
S. L. Woronowicz. Compact quantum groups.
In: {\it Sym\'etries Quantiques (Les Houches, 1995)},  845--884, North-Holland, Amsterdam, 1998.

\end{thebibliography}
\end{document}